\documentclass[11pt,a4paper]{article}

\usepackage[utf8]{inputenc}
\usepackage[ngerman, english]{babel}

\usepackage{amsmath}
\usepackage{amsfonts}
\usepackage{amssymb}
\usepackage{dsfont}
\usepackage{amsthm}
\usepackage{fancyhdr}
\usepackage[textheight=600pt]{geometry}
\usepackage[hang,flushmargin]{footmisc}
\usepackage[style=ext-authoryear, citestyle=authoryear, articlein=false, backend=biber, uniquename=init, giveninits=true, maxbibnames=10, maxcitenames=2]{biblatex}
\usepackage{enumitem}

\usepackage{csquotes}
\usepackage{comment}
\usepackage{setspace}
\usepackage{graphicx}
\usepackage{ifthen}

\usepackage{hyperref}
\hypersetup{breaklinks=true}

\usepackage[noabbrev]{cleveref}
\usepackage{authblk}
\usepackage[dvipsnames]{xcolor}
\theoremstyle{plain}
\newtheorem{theorem}{Theorem}[section]
\newtheorem{lemma}[theorem]{Lemma}
\newtheorem{proposition}[theorem]{Proposition}
\newtheorem{corollary}[theorem]{Corollary}
\newtheorem{definition}[theorem]{Definition}

\theoremstyle{definition}
\newtheorem{remark}[theorem]{Remark}
\newtheorem{example}[theorem]{Example}

\crefname{equation}{relation}{relations}

\newcommand{\E}{{\, \mathbb E}}
\newcommand{\R}{\mathbb{R}}
\newcommand{\N}{\mathbb{N}}

\renewcommand{\P}{\mathbb{P}}
\newcommand{\PC}[2]{\mathbb{P} \left( #1 \| #2 \right)}
\newcommand{\EC}[2]{\mathbb{E} \left( #1 \| #2 \right)}

\newcommand{\supnorm}[1]{\| #1 \|_\infty}
\newcommand{\lpnorm}[2]{\| #1 \|_{L_{#2} }}
\newcommand{\norm}[1]{\| #1 \|_{2}}
\newcommand{\limN}{\underset{n \to \infty}{\lim}}
\newcommand{\plim}{\overset{p}{\longrightarrow}}

\newcommand{\Fh}{\widehat{\mathrm{F}}_n}
\newcommand{\Ft}{\mathrm{F}_n}
\newcommand{\G}{\widehat{G}_n^{CV}}
\newcommand{\Gp}{\widehat{G}_n^{CV+}}
\newcommand{\Tn}{T_n}
\newcommand{\Tnl}[1]{T_{n}^{\backslash #1}}
\newcommand{\Tncv}[1]{T_{n}^{\backslash K_{ #1}}}
\newcommand{\pred}{\hat{y}_{n+1}}
\newcommand{\predl}[1]{\hat{y}^{\backslash #1}_{n+1}}
\newcommand{\PI}[3]{
    \ifthenelse{\equal{#2}{0}}
    {PI_{#1}^{#3}}
    {PI_{#1}^{#3}(#2)}
}
\newcommand{\PCv}[2]{\PI{#1}{#2}{CV}}
\newcommand{\PCvp}[2]{\PI{#1}{#2}{CV+}}
\newcommand{\PJi}[2]{\PJo{#1}(#2)}
\newcommand{\PJo}[1]{PI_{#1}^J}
\newcommand{\PJpi}[2]{\PJpo{#1}(#2)}
\newcommand{\PJpo}[1]{PI_{#1}^{J+}}

\newcommand{\avg}[1]{\dfrac{1}{#1}\sum_{i=1}^{#1}}
\newcommand{\sn}{\sigma_{n,\mathcal{P}_n}}
\newcommand{\snl}{\sigma_{n-1,\mathcal{P}_n}}
\newcommand{\bhl}[1]{\hat{\beta}_{\lambda}^{#1}}

\newcommand{\ldnamEnd}{L\'{e}vy gauge}
\newcommand{\ldname}{\ldnamEnd~}
\newcommand{\gam}{\delta}
\newcommand{\ld}[1]{\mathcal{L}_{#1}}
\newcommand{\LD}{\ld{\delta}(F,G)}
\newcommand{\Q}[2]{\mathcal{Q}_{#1}^{#2}}
\newcommand{\fatir}{\text{ for all } t \in \R}
\newcommand{\CCend}{\textbf{CC}}
\newcommand{\CC}{\CCend~}

\renewcommand{\epsilon}{\varepsilon}

\newcommand{\Op}[1]{O_p\left(#1\right)}

\newcommand\smallO[1]{
  \mathchoice
    {{\scriptstyle {O}_p}}
    {{\scriptstyle {O}_p}}
    {{\scriptscriptstyle {O}_p}}
    {\scalebox{.7}{$\scriptscriptstyle {O}_p$}}
  {\left(#1\right)}}

\newenvironment{nalign}{
    \begin{equation}
    \begin{aligned}
}{
    \end{aligned}
    \end{equation}
    \ignorespacesafterend
}

\newlist{lemmenum}{enumerate}{1}
\setlist[lemmenum]{label=\alph*), ref=\thelemma~(\alph*)}
\crefalias{lemmenumi}{lemma} 

\newlist{thmenum}{enumerate}{1}
\setlist[thmenum]{label=\alph*), ref=\thethm~(\alph*)}
\crefalias{thmenumi}{theorem} 
\title{Uncertainty quantification via cross-validation and its variants under algorithmic stability}
\author{Nicolai Amann}
\author{Hannes Leeb}
\author{Lukas Steinberger\footnote{
    Lukas Steinberger was supported by the Austrian Science Fund (FWF):
    I 5484-N, as part of the Research Unit 5381 of the German Research Foundation.
}}
\affil{University of Vienna}

\def \version{arxive}

\begin{document}
    \selectlanguage{english}
    \maketitle
    \begin{abstract}
Recently, there has been substantial interest in statistical guarantees for cross-validation (CV) methods of uncertainty quantification in statistical learning (cf. \cite{barber2021predictive}, \cite{liang2023algorithmic}, \cite{steinberger2020conditional}). These guarantees should hold under minimal assumptions on the data generating process and conditional on the training data, because numerous predictions are usually computed based on one and the same training sample. We push this objective to the limit: 
We prove asymptotic conditional conservativeness of CV, that is, 
the probability of the actual coverage probability, conditional on the training data, undershooting its nominal level vanishes asymptotically, under minimal assumptions.
In particular, we impose a stability condition, require that the prediction error is stochastically bounded, and show that neither condition can be dropped in general. 
By way of an asymptotic equivalence result, we also show that the closely related CV+ method of \textcite{barber2021predictive} provides exactly the same conditional statistical guarantees as CV in large samples, thereby extending the range of applicability of CV+ to the high-dimensional regime. We conclude that, in view of its marginal coverage guarantee, CV+ does indeed improve over simple CV. For our proofs we introduce a new concept called \ldnamEnd, which can be of independent interest.
\end{abstract}
	\section{Introduction}\label{sec:introduction}
Given a prediction algorithm and training data, it is important to quantify the corresponding prediction risk by, e.g., providing valid prediction intervals or estimating the misclassification error.
As there is typically only one set of training data available, one is interested in consistent estimation of the risk quantities \emph{conditional} on the specific training data rather than having a marginal guarantee which averages over all possible (but not available) training sets. For this reason, we focus on uncertainty quantification conditional on the training data.

\subsection{Background}
In the classical asymptotic setting where the number of parameters $p$ is fixed while the number of observations $n$ tends to infinity, consistent estimation of the underlying parameters is often possible. This case allows to consistently quantify the risk by using the fitted errors to get an approximation of the prediction error's distribution. However, if $p$ grows with $n$, this approach typically no longer works (cf. \cite{mammen1996empirical}) and even the bootstrap will fail (cf. \cite{elkaroui2018can} and \cite{bickel1983bootstrapping}). 
Such a setting, including the high-dimensional case, is precisely what applicants of modern data science are often confronted with, e.g., in image classification or (deep) neural networks. 
One possible solution for consistent risk estimation in this challenging scenario is sample splitting.
However, it comes at the price of using only a fraction of the data to train the prediction algorithm, which can lead to overly large prediction intervals compared to other methods (cf. \cite{steinberger2020conditional}).
An alternative approach avoiding this drawback is based on $k$-fold cross-validation (CV), its special case, the Jackknife, alongside with its variants such as the recently introduced $k$-fold cross-validation+ (CV+) and its special case, the Jackknife+ (cf. \cite{barber2021predictive}). 

While the special case of $n$-fold cross-validation, the so-called Jackknife, dates back to the Fifties (cf. \cite{quenouille1956notes}), 
little is known about this method in asymptotic settings where the data-generating process, in particular the number of parameters, may depend on $n$.
One recent exception is \textcite{steinberger2020conditional}, who showed that prediction intervals based on CV 
are asymptotically conditionally valid, that is, their actual 
coverage probability conditional on the training data converges to its nominal level also in high-dimensional asymptotic frameworks.
However, besides the natural condition of algorithmic stability, \textcite{steinberger2020conditional} also rely on a continuity assumption on the distribution of the response variable $y_{n+1}$ conditional on the feature vector $x_{n+1}$, which excludes many important cases such as classification. 

For the Jackknife+ even less is known, with the notable exception of \textcite{liang2023algorithmic}, who showed that prediction intervals based on the Jackknife+ are conditionally conservative, that is,
the probability of the actual coverage probability conditional on the training data undershooting the nominal level vanishes asymptotically.
For their asymptotic results, they require the algorithm to be stable with respect to the exclusion of one data point and the distribution of the response-feature pair $(y_{n+1}, x_{n+1})$ to be fixed over $n$ (see Theorem~3.6 in \cite{liang2023algorithmic} and \Cref{sub:appRem_mStability} for a discussion). The latter condition prevents the dimension $p$ of $x_{n+1}$ to grow with $n$ and is therefore only applicable in the classical fixed-$p$, large-$n$ framework where typically also the bootstrap and the naive approach based on fitted errors give a precise approximation of the true prediction error's distribution.
In particular, the asymptotic results of \textcite{liang2023algorithmic} do not cover the common situation in modern data science where the number $p$ of parameters is not negligible compared to the number $n$ of observations.

\subsection{Our contributions}
The four main contributions of this paper are as follows:
First, we show that CV provides asymptotically conservative prediction intervals 
\emph{conditional} on the training data under algorithmic stability without posing any other restriction to the underlying model than stochastic boundedness of the prediction error (cf. \Cref{thm:jack_asymptoticValidPi}). Our framework allows the dimension $p$ of the regressor to grow, thus including the challenging high-dimensional asymptotic framework.

Our results contain the Jackknife as a special case and extend to the Jackknife+ and CV+.
For the sake of simplicity, we present our findings for the Jackknife and Jackknife+ in the main part of the paper and defer the statements for CV and CV+ to \Cref{sec:CV}.
While our main results focus on uncertainty quantification through prediction intervals, we extend 
these findings to an accurate estimation of risk functions based on the prediction error, such as the conditional misclassification error and the conditional mean-squared prediction error, under mild conditions on the loss function (cf. \Cref{sub:appExt_risk}).

Our \emph{second finding} is the asymptotic equivalence of the Jackknife and the Jackknife+ and the equivalence of CV and CV+ under algorithmic stability (cf. \Cref{thm:eqi_PacBoundAsymptotic}). Thus, all (our and other) existing and future results concerning one of these methods immediately carry over to their Jackknife or Jackknife+ variant given only algorithmic stability.
Even more, our results are based on finite sample error bounds showing that even for a fixed $n$ both Jackknife variants are similar to each other under algorithmic stability (see \Cref{sub:appRem_finiteSampleResults}). Consequently, future research may focus on one of the variants and use the finite sample similarity or asymptotic equivalence results to carry over their results to the other variants under algorithmic stability. 

In the \emph{third part} of our work, we investigate for an arbitrary algorithm the necessity of our assumptions by analyzing the stability condition and the stochastic boundedness of the prediction error: We do so by, firstly, proving that the expected length of any prediction interval containing the predictor (such as the symmetric version of the Jackknife) diverges to infinity unless the prediction error is stochastically bounded (cf. \Cref{prop:nec_expectedInfinityLength}).
Secondly, we show that the stability condition for the Jackknife essentially cannot be dropped for \emph{any} prediction algorithm whose prediction error has uniformly bounded $(2+\epsilon)$-th absolute moments (see \Cref{sub:necessityStab}).

Our \emph{fourth contribution} is the introduction of a new statistical distance measure that we call the \ldname (see \Cref{sec:ld}). It is an essential tool in our proofs and can be of independent interest.
In fact, its application is not restricted to CV, but can also be useful for the analysis of other methods of uncertainty quantification, such as conformal prediction (see \Cref{sec:appRemarks}).

\subsection{Organization of the paper}
In \Cref{sec:setting}, we introduce the setting of our analysis, the Jackknife, and related notation.
We present the asymptotic results for the Jackknife prediction intervals in \Cref{sec:jackknife}. In \Cref{sec:equivalence} we show that the same holds true for the Jackknife+ by proving an equivalence result between the two methods under algorithmic stability. In \Cref{sec:necessity} we discuss the necessity of our assumptions.

In fact, all our results also apply to CV and CV+ as long as the number $k$ of folds diverges. For the sake of clarity, we have decided to present only the special case $k=n$, that is, the Jackknife and the Jackknife+, in \Cref{sec:jackknife} and \Cref{sec:equivalence}, while the general statement can be found in \Cref{sec:CV}.

We present finite sample results for CV methods in \Cref{sub:appRem_finiteSampleResults} and several extensions, including the generalization of uncertainty quantification to other risk measures, in \Cref{sec:appRemarks}. 
Our proofs heavily rely on a new concept called \ldname which is introduced in \Cref{sec:ld}.
We defer all proofs to the last chapters of the%
\ifthenelse{\equal{\version}{arxive}}{
 Appendix.
}{
 online supplement.
}

\section{Setting}\label{sec:setting}
\subsection{The data}
Since we want to embed our finite sample results in an asymptotic regime that accommodates high-dimensional settings, 
our framework allows all quantities to depend on the sample size $n$:
For every $n \geq 2$ we assume that $\mathcal{P}_n$ is a Borel-measure on the space 
$\R \times \R^{p_n}$ and we are given training data $\Tn$ consisting of 
$n$ independent and identically distributed (i.i.d.) response-feature pairs $(y_{i}^{(n)}, x_{i}^{(n)})$ with $1 \leq i \leq n$ distributed according to $\mathcal{P}_n$, 
where $y_{i}^{(n)}$ is a real-valued random variable and $x_{i}^{(n)}$ is a random vector in $\R^{p_n}$.\footnote{
	Our arguments can easily be generalized to the case where we are given no regressors at all ($p_n=0$), to the case of randomized prediction algorithms, and cases where the feature vectors $x_i$ are contained in a more general space $\mathcal{X}$.
}  
In particular, the dimension $p_n$ may grow with $n$.
We are interested in predicting another random variable $y_{n+1}^{(n)}$, where
we are additionally given a new regressor $x_{n+1}^{(n)}$.
Furthermore, we assume that the new response-feature pair $(y_{n+1}^{(n)}, x_{n+1}^{(n)})$ is also distributed according to $\mathcal{P}_n$ and independent of the training data $\Tn = (y_{i}^{(n)}, x_{i}^{(n)})_{i=1}^n$. 
For the sake of readability, we will suppress the dependence on $n$ of the aforementioned expressions whenever it is clear from the context. Furthermore, we simply write $\P$ if the underlying probability measure is clear from the context.

\subsection{The prediction algorithm}
Throughout this paper, we take the prediction algorithm as given and construct prediction intervals and estimators for the prediction error using this algorithm. This includes both accurate and ill-suited algorithms.

A prediction algorithm is a procedure that uses the training data $\Tn$ together with a new feature vector $x_{n+1}$ to predict the response variable associated with that feature vector. 
To formalize this, we fix $p \in \N, n \in \N$, assume $\mathcal{A}_{p,n}$ to be a Borel measurable function from $\R^p \times \R^{(p+1)n}$ to $\R$ and call $\pred := \mathcal{A}_{p,n}(x_{n+1}, \Tn)$ the predictor for $y_{n+1}$ given $x_{n+1}$ and the training data $\Tn$.\footnote{
	If we want to allow for randomized predictors, we replace the domain of the prediction algorithm by $\R^p \times \R^{(p+1)n} \times \Omega$, where $\Omega$ is the sample space of a probability space $(\Omega, \mathcal{F}, P)$.
}

A prediction algorithm $\mathcal{A}_{p,n}$ is called symmetric if the predictor is not changed by a permutation of the training data, that is, for all $x_{n+1} \in \R^p$, any permutation $\pi: \{1, \ldots, n\} \to \{1, \ldots, n\}$ and training data 
$\Tn = (y_i, x_i)_{i=1}^n$ the prediction algorithm satisfies 
$\mathcal{A}_{p,n}(x_{n+1}, \Tn) = \mathcal{A}_{p,n}(x_{n+1}, (y_{\pi(i)}, x_{\pi(i)})_{i=1}^n)$. 
Being somewhat imprecise, we call a predictor $\pred$ symmetric if the underlying prediction algorithm $\mathcal{A}_{p,n}$ is symmetric for all $n \in \N$ and $p \in \N$.

\subsection{The goal}
Assume we have chosen a method to create a prediction interval $PI_{\alpha_1, \alpha_2}$ based on the training data $\Tn$ and $x_{n+1}$ with \emph{nominal} coverage probability of $\alpha_2 - \alpha_1$ for some unknown $y_{n+1}$. 
Our main goal is to analyze under which conditions these prediction intervals are conditionally \emph{conservative} in large samples, that is,
$$\PC{y_{n+1} \in PI_{\alpha_1, \alpha_2}}{\Tn} \geq (\alpha_2 - \alpha_1) - \smallO{1}.$$
Unless noted otherwise in the following, ``conditional'' always means ``conditional on the training data''.\footnote{
    For the analysis of the coverage probability of prediction intervals conditional on $x_{n+1}$ and its limitation we refer to 
    \textcite{foygel2021limits}.
}
In order to prevent this prediction interval to be overly large, we additionally want to analyze whether the prediction interval is even 
conditionally \emph{valid} in large samples, that is,
$$\PC{y_{n+1} \in PI_{\alpha_1, \alpha_2}}{\Tn} - (\alpha_2 - \alpha_1) = \smallO{1}.$$
Since any non-randomized prediction interval may fail to be valid if the distribution of $(y_{n+1}, x_{n+1})$ is discrete,\footnote{
	For example, assume the data $(y_{n+1}, x_{n+1})$ are Dirac distributed. Then, the coverage probability of any non-randomized prediction interval is either $0$ or $1$.
} 
we restrict the analysis of validity to the case of continuous distribution functions. We emphasize that for the general analysis of conservativeness, we do not pose this continuity assumption.
Thus, in our results we will distinguish between the \emph{general case} and the \emph{continuous case} by using the following definition:
\begin{definition}[\CC Assumption]\label{def:set_CC}
	Let $n \in \N$ and $\mathcal{P}_n^x$ denote the marginal distribution of $x_{1}^{(n)}$.
	We say the Continuous Case Assumption (in short \CCend) is fulfilled if the following two properties hold true.
	\begin{enumerate}[label=CC.\alph*)]
		\item For every $n \in \N$ and $\mathcal{P}_n^x$-almost every $x$ the random variable $y_{1}^{(n)}$ conditional on $x_{1}^{(n)} = x$ is absolutely continuous 
		      and the supremum norm $\supnorm{f_{y_{1}^{(n)} \| x_{1}^{(n)} = x}}$ of its density $f_{y_{1}^{(n)} \| x_{1}^{(n)} = x}$ is finite.
		\item The sequence of random variables $(\supnorm{f_{y_{1}^{(n)} \| x_{1}^{(n)}}})_{n \in \N}$ is stochastically bounded.
	\end{enumerate} 
\end{definition}

\subsection{Notation}
For a measurable function $f: \R \to \R$ we will denote its supremum by $\supnorm{f} := \sup_{x \in \R} |f(x)|$.
Furthermore, $e_i$ denotes the $i$-th canonical basis vector and $\norm{\cdot}$ the Euclidean norm.
For a set $S \subseteq \R^k$ we denote by $|S|$ its cardinality and by $\mathds{1}_S: \R^k \to \{0,1\}$ the indicator function on the set $S$. 
For a real number $x$, the expression $\lceil x \rceil$ describes the smallest integer larger than or equal to $x$.
The expression $\lim_{x \searrow t}$ will denote the limit from above (in contrast to the usual limit $\lim_{x \to t}$). 
For a distribution function $F$, we define $F(x-) := \lim_{\delta \searrow 0}F(x-\delta)$ to be the limit from the left at the point $x \in \R$.

For a sequence of random variables $X_n$, we will abbreviate convergence in probability to a random variable $X$ by $X_n \plim X$. 
Furthermore, for a positive sequence $(v_n)_{n \in \N}$ we will write $X_n \sim \smallO{v_n}$ if the sequence $(X_n v_n^{-1})_{n \in \N}$ converges to $0$ in probability and $X_n \sim \Op{v_n}$ if the sequence $(X_n v_n^{-1})_{n \in \N}$ is stochastically bounded. 
In particular, we write $X_n \sim \smallO{1}$ or $X_n \plim 0$ to describe convergence in probability of $(X_n)_{n \in \N}$ to $0$ and $X_n \sim \Op{1}$ to denote stochastic boundedness of $(X_n)_{n \in \N}$. We define the infimum of the empty set as $+\infty$ and the supremum of the empty set as $-\infty$.

For the definition of prediction sets, we need to define $\alpha$-quantiles of distribution functions. In order to avoid case distinctions in our proofs and statements, we extend the definition to all $\alpha \in \R$.
\begin{definition}[Quantiles]\label{def:set_quantiles}
	Let $F: \R \to [0,1]$ be a cumulative distribution function. 
	For $\alpha \in \R$ we define the $\alpha$-quantiles $\Q{\alpha}{F}$ of $F$ as
	\begin{align*}
		\Q{\alpha}{F} := \inf\{x \in \R: F(x) \geq \alpha\} \in [-\infty, +\infty]. 
	\end{align*}
    Recalling that we defined the infimum of the empty set $\inf\{x \in \R: x \in \emptyset\}$ as $+\infty$, 
	we have $\Q{\alpha}{F} = - \infty$ whenever $\alpha \leq 0$ and $\Q{\alpha}{F} = +\infty$ whenever $\alpha > 1$.
\end{definition}

We immediately conclude from \Cref{def:set_quantiles} that
\begin{align}\label{eq:ld_alphaQuantilesExceedAlpha}
	F(\Q{\alpha}{F}-) \leq \alpha \leq F(\Q{\alpha}{F}) 
\end{align}
holds true whenever $\Q{\alpha}{F}$ is finite.

\section{Asymptotic validity of the Jackknife}\label{sec:jackknife}
The Jackknife estimates the distribution function 
\begin{align*}
	\Ft(t) = \PC{y_{n+1} - \pred \leq t}{\Tn} 
\end{align*}
of the prediction error $y_{n+1} - \pred$  conditional on the training data. 
Let $1 \leq i \leq n$ and denote by $\Tnl{i}$ the training data where the $i$-th pair $(y_i, x_i)$ is removed. 
For $1 \leq j \leq n+1$ and $1 \leq i \leq n$ let $\hat{y}_j^{\backslash i}$ be the prediction for $y_j$ based on the regressor $x_j$ and the training data $\Tnl{i}$. 
Define the \emph{leave-one-out} residuals $\hat{u}_i = y_i - \hat{y}_i^{\backslash i}$.
The Jackknife estimates $\Ft$ via the empirical distribution function $\Fh$ of the leave-one-out residuals:
\begin{align*}
	\Fh(t) = \dfrac{1}{n} \sum_{i=1}^n \mathds{1}\{\hat{u}_i \leq t \}.
\end{align*}

Consider $\Q{\alpha}{\Fh}$, the $\alpha$-quantile for $\alpha \in \R$ based on the Jackknife, and note that
\begin{align*}
	\Q{\alpha}{\Fh} = 
	\begin{cases}
	\hat{u}_{( \lceil\alpha n\rceil )} & \text{ if } \alpha \in (0, 1] \\
	-\infty                            & \text{ if } \alpha = 0,       
	\end{cases}
\end{align*}
where $\hat{u}_{(1)} \leq \hat{u}_{(2)} \leq \ldots \leq \hat{u}_{(n)}$ are the ordered leave-one-out residuals. 
Let $0 \leq \alpha_1 \leq \alpha_2 \leq 1$ and $\gam \in \R$.
Define the $\gam$-distorted prediction intervals based on the Jackknife
with nominal coverage probability of $\alpha_2 - \alpha_1$ for $y_{n+1}$ as
\begin{align}\label{def:jack_PiJackknifeInflated}
	\PJi{\alpha_1, \alpha_2}{\gam} = \pred + [\Q{\alpha_1}{\Fh} - \gam, \Q{\alpha_2}{\Fh} + \gam], 
\end{align}
where we use the convention that $[a,b] = \emptyset$ if $a > b$. For convenience, we extend this definition to all $\alpha_1, \alpha_2 \in \R$. For $\gam > 0$, we will refer to $\PJi{\alpha_1, \alpha_2}{\gam}$ as inflated prediction intervals, and for $\gam < 0$, we call them shrunken prediction intervals. We will occasionally abbreviate the non-distorted prediction intervals $\PJi{\alpha_1, \alpha_2}{0}$ with $\PJo{\alpha_1, \alpha_2}$.

Since $\Tnl{i}$ is independent of $x_i$ and $y_i$, the leave-one-out approach is able to capture the inherent nature of independence of the training data $\Tn$ from the new observation $y_{n+1}$ and therefore avoids underestimating the out-of-sample prediction error. However, the approach may fail if the inclusion of one observation to the training data results in a highly different prediction. We will refer to the property that a new observation does not change the prediction too much as the \emph{out-of-sample stability} of a predictor and formalize this as follows:
\begin{definition}[Asymptotic out-of-sample stability]\label{def:jack_oosStability}
	A predictor $\pred$ is asymptotically out-of-sample stable with respect to its leave-one-out analogues $(\predl{1}, \ldots, \predl{n})$
	if
	\begin{align*}
		\limN \avg{n} \P(|\pred - \predl{i}| \geq \epsilon) = 0 \text{ for all } \epsilon > 0. 
	\end{align*}
	Unless noted otherwise, an asymptotically out-of-sample stable predictor will be called stable. 
\end{definition}
In particular, for a symmetric predictor the distribution of $\pred - \predl{i}$ does not change with $i$ and therefore the asymptotic out-of-sample stability condition reduces to the assumption that 
$\pred - \predl{n}$ converges to zero in probability.

\begin{theorem}[Asymptotically conservative/valid prediction intervals with the Jackknife]\label{thm:jack_asymptoticValidPi}
	Suppose the prediction error $y_{n+1} - \pred$ is stochastically bounded and the predictor is stable.
	Then the following statements hold true.
	\begin{enumerate}[label=(\roman*)]
		\item \label{it:jack_asymptoticValidPi_GenCase}\emph{General case:} For every $\gam > 0$ the $\gam$-inflated Jackknife prediction intervals 
		      are asymptotically uniformly conditionally conservative, that is,
		      \begin{align*}
		      	\limN \P \left( \inf_{0 \leq \alpha_1 \leq \alpha_2 \leq 1}                                   
		      	\left[ \PC{ y_{n+1} \in \PJi{\alpha_1, \alpha_2}{\gam} }{\Tn} - (\alpha_2 - \alpha_1) \right] 
		      	\geq -\epsilon \right) = 1 \text{ for all } \epsilon > 0.                                     
		      \end{align*}
		      Additionally, the conditional coverage probability of $\gam$-shrunken prediction intervals uniformly does not overshoot its nominal level asymptotically, that is,
		      \begin{align*}
		      	\limN \P \left( \sup_{0 \leq \alpha_1 \leq \alpha_2 \leq 1}                                    
		      	\left[ \PC{ y_{n+1} \in \PJi{\alpha_1, \alpha_2}{-\gam} }{\Tn} - (\alpha_2 - \alpha_1) \right] 
		      	\leq \epsilon \right) = 1 \text{ for all } \epsilon > 0.                                       
		      \end{align*}
		\item \emph{Continuous Case: }If, additionally, Assumption \CC given in \Cref{def:set_CC} is satisfied, the non-distorted prediction interval ($\gam = 0$) is asymptotically uniformly conditionally valid, that is,
		      \begin{align*}
		      	\limN \E \left( \sup_{0 \leq \alpha_1 \leq \alpha_2 \leq 1}                                             
		      	\left| \PC{ y_{n+1} \in \PJi{\alpha_1, \alpha_2}{0} }{\Tn} - (\alpha_2 - \alpha_1) \right| \right) = 0. 
		      \end{align*}
	\end{enumerate}
\end{theorem}

\begin{remark}
	\Cref{thm:jack_asymptoticValidPi} can be generalized to $k$-fold cross-validation as long as the number of folds $k$ diverges to infinity asymptotically.
	The general statement together with the definition of leave-one-fold-out stability
	can be found in \Cref{sec:CV}.
\end{remark}

\begin{remark}
	Since the two statements for the general case in \Cref{thm:jack_asymptoticValidPi} hold true for every $\gam > 0$ and every $\epsilon > 0$, there are null-sequences $(\gam_n)_{n \in \N}$ and $(\epsilon_n)_{n \in \N}$ such that these statements remain true with $\gam_n$ and $\epsilon_n$ replacing $\delta$ and $\epsilon$, respectively. In the typical setting where the prediction error does not converge to a point $x \in \R$, the length of a conservative prediction interval (for large enough $\alpha_2 - \alpha_1$) is nonzero asymptotically. Thus, in such a setting, the artificial inflation $\gam_n$ is negligible compared to the length of the original interval for large samples.
\end{remark}

\begin{remark}
	The values $\alpha_1$ and $\alpha_2$ may depend on $n \in \N$ and, more importantly, on the training data.
	Thus, a practitioner can fix a value $\alpha \in (0,1)$, choose $\alpha_1$ and $\alpha_2$ to be the shortest Jackknife prediction interval for the training data with $\alpha_2 - \alpha_1 = 1 - \alpha$, and still get the same asymptotic coverage guarantees. 
	Furthermore, this optimization can be efficiently computed since there are no more than $\lceil \alpha n \rceil + 1$ different prediction intervals of the form \eqref{def:jack_PiJackknifeInflated} with nominal coverage probability of $1-\alpha$.
\end{remark}

\Cref{thm:jack_asymptoticValidPi} generalizes Theorem $2.4$ of \textcite{steinberger2020conditional} in two ways: Our major novelty is the extension to the general, possibly non-continuous, case where the Jackknife provides asymptotically uniformly conditionally conservative prediction intervals for stable algorithms if we are additionally willing to inflate the prediction intervals by some (small) $\delta > 0$. 
Secondly, in the continuous case our statement provides a worst-case analysis by holding uniformly over all nominal levels $\alpha_2 - \alpha_1$ rather than pointwise for fixed $0 \leq \alpha_1 \leq \alpha_2 \leq 1$.

We stress that \Cref{thm:jack_asymptoticValidPi} does not cover non-distorted prediction intervals in the non-continuous case.
Nevertheless, the Jackknife prediction intervals are not overly large since any arbitrarily small shrinkage of the prediction interval prevents the actual conditional coverage probability from overshooting the nominal level asymptotically.\footnote{Additional considerations regarding an upper bound for the actual conditional coverage probability of inflated prediction intervals can be found in \Cref{sub:appRem_upperBoundsCovProb}.}

Furthermore, the statements in the first part of \Cref{thm:jack_asymptoticValidPi} also entail for $\gam > 0$ the (weaker) marginal coverage guarantees
\begin{align*}
	\limN \P \left( y_{n+1} \in \PJi{\alpha_1, \alpha_2}{-\gam} \right) \leq \alpha_2 - \alpha_1 \leq \limN \P \left( y_{n+1} \in \PJi{\alpha_1, \alpha_2}{\gam} \right). 
\end{align*}
In fact, this is closely related to Theorem~$5$ of \textcite{barber2021predictive}, who proved marginal guarantees for the special case of the symmetrized Jackknife for symmetric predictors fulfilling a similar stability condition. 

We emphasize that \Cref{thm:jack_asymptoticValidPi} is based on a finite sample inequality stated as a PAC-bound\footnote{
	Here, a probably approximately correct (PAC) bound denotes an inequality of the form $\P( (1-\alpha) - \PC{y_{n+1} \in PI}{\Tn} \geq \delta) \leq \epsilon$ for some (small) $\delta, \epsilon > 0$.
} which shows that even for a fixed $n$ the conditional coverage probability of Jackknife prediction intervals avoids undershooting its nominal level with high probability if the underlying algorithm fulfills a finite sample analogue of the stability condition (cf. \Cref{prop:appRem_finSamCV}). 

Furthermore, the Jackknife can also be used to consistently estimate other measures of uncertainty quantification such as the conditional misclassification error or the conditional mean-squared prediction error (cf. \Cref{sub:appExt_risk}).

\section{Asymptotic equivalence of the Jackknife and the Jackknife+}\label{sec:equivalence}
Before defining the Jackknife+ prediction intervals we would like to note that for every $\alpha \in [0,1]$ we can rewrite $\pred + \Q{\alpha}{\Fh}$ as $\Q{\alpha}{\widehat{G}_n}$ where $\widehat{G}_n(t) = \Fh(t - \pred) = \frac{1}{n} \sum_{i=1}^n \mathds{1}_{[\pred + \hat{u}_i, \infty)}(t)$ is the empirical distribution function of 
\begin{align}\label{eq:eqi_defVi}
    v_i := \pred + \hat{u}_i. 
\end{align}
This gives the following alternative representation of the Jackknife prediction intervals:
\begin{align}\label{eq:eqi_PiWithTheJackAlternative}
    \PJi{\alpha_1, \alpha_2}{\gam} = [\Q{\alpha_1}{\widehat{G}_n} - \gam, \Q{\alpha_2}{\widehat{G}_n} + \gam]. 
\end{align}

Now the Jackknife+ is defined similar to the Jackknife after replacing $v_i$ by 
\begin{align}\label{eq:eqi_defViPlus}
    v_i^+ := \predl{i} + \hat{u}_i 
\end{align}
in the definition above.
To put it in other words: For any $\gam \in \R$ the $\gam$-distorted prediction interval based on the Jackknife+ with nominal coverage probability of $\alpha_2 - \alpha_1$ is defined as
\begin{align}\label{eq:eqi_PiWithTheJackp}
    \PJpi{\alpha_1, \alpha_2}{\gam} = [\Q{\alpha_1}{\widehat{G}_n^{+}} - \gam, \Q{\alpha_2}{\widehat{G}_n^{+}} + \gam], 
\end{align}
where $\widehat{G}_n^{+}$ denotes the empirical distribution function of $(v_i^+)_{i =1}^n$. Like for the Jackknife, we occasionally abbreviate the non-distorted Jackknife+ prediction intervals $\PJpi{\alpha_1, \alpha_2}{0}$ with $\PJpo{\alpha_1, \alpha_2}$.
Note that we have decided to use the asymmetric version of the Jackknife+ as we are dealing with the asymmetric version of the Jackknife as well.\footnote{
    We would like to point out that our definition of the asymmetric Jackknife+ slightly differs from the one given in Appendix A in \textcite{barber2021predictive} due to rounding issues: Our prediction intervals with nominal level $\alpha_2 - \alpha_1$ coincide with the asymmetric prediction intervals of \textcite{barber2021predictive} with nominal level $\frac{\lceil \alpha_2 n \rceil - \lceil \alpha_1 n \rceil }{n+1}$.} 
In fact, all our results also hold true if we use the symmetrized versions of the Jackknife and the Jackknife+ (cf. \Cref{sub:appExt_symmetrizedVersions}).
Equipped with these definitions, we present the main result of this section:
\begin{theorem}[Asymptotic equivalence of the Jackknife and the Jackknife+]\label{thm:eqi_PacBoundAsymptotic}
    Assume the predictor is stable as in \Cref{def:jack_oosStability}.
    \begin{enumerate}[label=(\roman*)]
        \item \emph{General case:} 
              Then the following statements are equivalent.
              \begin{enumerate}[label=G\arabic*)]
                \item\label{it:eqi_PacBoundAsymptoticGenJack} 
                For all $\gam > 0$ and $\epsilon > 0$ we have
                \begin{align*}
                    \limN \P \left( \inf_{0 \leq \alpha_1 \leq \alpha_2 \leq 1} \PC{y_{n+1} \in \PJi{\alpha_1, \alpha_2}{\gam}}{\Tn} 
                    - (\alpha_2 - \alpha_1) \geq - \epsilon \right) = 1.                                                             
                \end{align*}
                \item\label{it:eqi_PacBoundAsymptoticGenJackp} 
                For all $\gam > 0$ and $\epsilon > 0$ we have
                \begin{align*}
                    \limN \P \left( \inf_{0 \leq \alpha_1 \leq \alpha_2 \leq 1} \PC{y_{n+1} \in \PJpi{\alpha_1, \alpha_2}{\gam}}{\Tn} 
                    - (\alpha_2 - \alpha_1) \geq - \epsilon \right) = 1.                                                              
                \end{align*}
              \end{enumerate}
        \item \emph{Continuous case:} If, in addition, Assumption \CC is satisfied, the following statements are equivalent.
              \begin{enumerate}[label=C\arabic*)]
                \item\label{it:eqi_PacBoundAsymptoticConJack}
                    $\limN \E(\sup_{0 \leq \alpha_1 \leq \alpha_2 \leq 1} |\PC{y_{n+1} \in \PJi{\alpha_1, \alpha_2}{0}}{\Tn} - (\alpha_2 - \alpha_1)|) = 0$.
                \item\label{it:eqi_PacBoundAsymptoticConJackp}
                    $\limN \E(\sup_{0 \leq \alpha_1 \leq \alpha_2 \leq 1} |\PC{y_{n+1} \in \PJpi{\alpha_1, \alpha_2}{0}}{\Tn} - (\alpha_2 - \alpha_1)|) = 0.$
              \end{enumerate}
    \end{enumerate}
\end{theorem}

\Cref{thm:eqi_PacBoundAsymptotic} shows that in the continuous case one prediction interval is valid if and only if the other one is, provided the predictor is stable.
In the general case, \Cref{thm:eqi_PacBoundAsymptotic} yields a similar statement for the prediction intervals being asymptotically conditionally conservative. In fact, also a similar equivalence result for the shrunken prediction intervals not overshooting their nominal level can be proved (cf. \Cref{prop:appJa_ShrunkenEqInflated} and \Cref{rem:appJa_ShrunkenEquivalence}). 

Combining \Cref{thm:jack_asymptoticValidPi} with \Cref{thm:eqi_PacBoundAsymptotic} yields the following corollary:
\begin{corollary}[Asymptotically conservative/valid prediction intervals with the Jackknife+]\label{cor:eqi_asymptoticValidPiPlus}
    Assume the prediction error $y_{n+1} - \pred$ is stochastically bounded and the predictor is stable. Then, the following statements hold true.
    \begin{enumerate}[label=(\roman*)]
        \item \emph{General case:} For every $\gam > 0$ the $\gam$-inflated Jackknife+ prediction intervals 
              are asymptotically uniformly conditionally conservative, that is,
              \begin{align*}
                \limN \P \left( \inf_{0 \leq \alpha_1 \leq \alpha_2 \leq 1} \PC{ y_{n+1} \in \PJpi{\alpha_1, \alpha_2}{\gam} }{\Tn} - (\alpha_2 - \alpha_1) \geq -\epsilon \right) = 1 \text{ for all } \epsilon > 0. 
              \end{align*}
              Additionally, the conditional coverage probability of $\gam$-shrunken prediction intervals uniformly does not overshoot its nominal level asymptotically, that is,
              \begin{align*}
                \limN \P \left( \sup_{0 \leq \alpha_1 \leq \alpha_2 \leq 1} \PC{ y_{n+1} \in \PJpi{\alpha_1, \alpha_2}{-\gam} }{\Tn} - (\alpha_2 - \alpha_1) \leq \epsilon \right) = 1 \text{ for all } \epsilon > 0. 
              \end{align*}
        \item \emph{Continuous Case: }If, in addition, Assumption \CC is satisfied, the non-distorted prediction interval ($\gam = 0$) is asymptotically uniformly conditionally valid, that is,
              \begin{align*}
                \limN \E \left( \sup_{0 \leq \alpha_1 \leq \alpha_2 \leq 1}                                              
                \left| \PC{ y_{n+1} \in \PJpi{\alpha_1, \alpha_2}{0} }{\Tn} - (\alpha_2 - \alpha_1) \right| \right) = 0. 
              \end{align*}
    \end{enumerate}
\end{corollary}

\Cref{cor:eqi_asymptoticValidPiPlus} shows that indeed the Jackknife+ provides asymptotically uniformly conditionally conservative prediction intervals in the general case if the prediction intervals are slightly inflated by some $\delta > 0$, while in the continuous case, the non-distorted Jackknife+ prediction intervals are even asymptotically uniformly conditionally valid. Thus, \Cref{cor:eqi_asymptoticValidPiPlus} generalizes Theorem $3.6$ of \textcite{liang2023algorithmic} in several ways: Firstly, we do not need to fix the distribution $\mathcal{P}_n$ of $(y_{n+1}, x_{n+1})$ over $n \in \N$ and therefore allow the dimension $p$ of $x_{n+1}$ to depend and grow with $n$.\footnote{
    A detailed discussion why Theorem $3.6$ of \cite{liang2023algorithmic} crucially relies on the assumption of a fixed distribution can be found in \Cref{sub:appRem_mStability}.} 
In particular, our results also contain the high-dimensional asymptotic framework as a special case which is a typical setting a statistician is confronted with, e.g., in image classification or deep neural networks.
Secondly, the additional assumption of a stochastically bounded prediction error allows us to prove the uniformity over all prediction intervals with $0 \leq \alpha_1 \leq \alpha_2 \leq 1$. 
Furthermore, both \Cref{thm:eqi_PacBoundAsymptotic} and \Cref{cor:eqi_asymptoticValidPiPlus} are special cases of a more general result concerning CV+ and its connection to the classical CV, which can be found in \Cref{sec:CV}. 

\begin{remark}
    The equivalence conditions in \Cref{thm:eqi_PacBoundAsymptotic} require the prediction intervals to be asymptotically conditionally conservative (or valid, in the continuous case) \emph{uniformly} over all nominal levels $\alpha_2 - \alpha_1$ for any combination $0 \leq \alpha_1 \leq \alpha_2 \leq 1$.
    This property is indeed fulfilled for the Jackknife if the predictor is stable and the prediction error is stochastically bounded.
    It might also be of interest to have such an equivalence result which holds pointwise for all $0 \leq \alpha_1 \leq \alpha_2 \leq 1$ in the sense that the infimum in \emph{both} statements is located outside of the outer probability in the general case and the supremum is located outside of the expectation in the continuous case.
    Indeed, with a slight modification of the proof, such a pointwise statement can also be made. We omit the details for the sake of brevity.
\end{remark}

\section{Necessity of the imposed conditions}\label{sec:necessity}
In \Cref{sec:jackknife} we require algorithmic stability and stochastic boundedness of the prediction error to prove the validity of prediction intervals based on the Jackknife.
In this section, we argue that essentially neither of these two conditions can be removed.
We stress that under these conditions also the Jackknife+ provides the same asymptotic conditional coverage guarantees as the original Jackknife.

\subsection{On stochastic boundedness}
We consider the assumption of stochastic boundedness of the prediction error to be no true restriction:
In fact, even if we only ask for asymptotically \emph{marginal} conservative prediction intervals based on $\pred$, the expected length of these prediction intervals diverges to $\infty$ along a subsequence $(n_k)_{n \in \N}$ unless the prediction error $y_{n+1} - \pred$ is stochastically bounded:

\begin{proposition}\label{prop:nec_expectedInfinityLength}
    Let $(PI_{\alpha_1, \alpha_2}^{n, p})_{0 \leq \alpha_1 \leq \alpha_2 \leq 1, n \in \N, p \in \N}$ be a method for creating prediction intervals\footnote{
        That is, for every fixed $0 \leq \alpha_1 \leq \alpha_2 \leq 1$ and $n,p \in \N$, 
        $PI_{\alpha_1, \alpha_2}^{n,p}(x_{n+1}, \Tn)$ is an interval of the form $[L,U]$, where $L$ and $U$ are measurable functions.
    }.
    Assume this method provides asymptotically \emph{marginal} conservative prediction intervals for every $0 \leq \alpha_1 \leq \alpha_2 \leq 1$, that is,
    \begin{align*}
        \liminf_{n \to \infty} \P(y_{n+1} \in PI_{\alpha_1, \alpha_2}^{n,p_n}(x_{n+1}, \Tn)) \geq \alpha_2 - \alpha_1. 
    \end{align*}
    Let $\pred$ be any predictor that is contained in the prediction interval almost surely. 
    If the prediction error $y_{n+1} - \pred$ is not stochastically bounded, then there is an $\epsilon > 0$ such that for all $\alpha_2 - \alpha_1 \geq 1 - \epsilon$
    the expected length of the prediction interval diverges to $\infty$ along a subsequence $(n_k)_{k \in \N}$.
\end{proposition}

\Cref{prop:nec_expectedInfinityLength} justifies the assumption of stochastic boundedness of the prediction error:
Since for the symmetric version of the Jackknife the predictor is always contained in the prediction interval,
\Cref{prop:nec_expectedInfinityLength} shows that the symmetrized Jackknife can only provide asymptotically conservative prediction intervals if its expected length diverges to infinity (along a subsequence $(n_k)_{k \in \N}$) or the prediction error is stochastically bounded.
In fact, with a slightly more refined technique a similar statement for the symmetric version of the Jackknife+ can be proved given the stability of the predictor.\footnote{Note that for the symmetrized Jackknife+ the predictor $\pred$ need not be contained in the corresponding prediction interval and therefore \Cref{prop:nec_expectedInfinityLength} cannot be applied directly.} We omit the details for the sake of brevity.

Note that this phenomenon is not limited to the (symmetric) Jackknife and its variants alone since \Cref{prop:nec_expectedInfinityLength} holds true for arbitrary methods of creating prediction intervals. Consequently, if one desires a prediction interval that includes $\pred$ to be of finite length, the associated prediction error has to be stochastically bounded.

\subsection{On the stability condition}\label{sub:necessityStab}
In this subsection, we discuss whether the Jackknife can provide conditionally conservative (or valid) prediction intervals in large samples even without the stability condition given in \Cref{def:jack_oosStability}. In such a case, the equivalence result typically no longer holds, and therefore it might be possible to find scenarios where the Jackknife indeed outperforms the Jackknife+ in terms of its conditional coverage guarantees. 

Before stating our results, we point out that there exist non-stable algorithms such that the Jackknife, alongside with the Jackknife+, do not provide conditionally conservative prediction intervals (cf. Theorem $2$ in \cite{barber2021predictive} and Theorem $3$ in \cite{bian2022training}).
However, this does not answer the question whether for \emph{every} non-stable algorithm these methods fail to give conditionally conservative prediction intervals. 
In this subsection, we argue that the stability assumption can essentially not be removed for any predictor with uniformly bounded $(2+\xi)$-th absolute moments, thereby complementing the counterexamples of \textcite{barber2021predictive} and \textcite{bian2022training}.

For this, we start with the following equivalent characterization of stability for prediction algorithms. For the sake of simplicity, we restrict our analysis to the case of symmetric predictors.
\begin{lemma}[Equivalent stability condition in $L_{2 + \xi}$]\label{lem:nec_equivalentStability}
    Let $\pred$ be a
    symmetric predictor and assume the $(2+\xi)$-th absolute moments of $\pred$ and $\predl{n}$ are uniformly bounded over $n$ for some $\xi > 0$.
    Furthermore, let $\sn^2$ and $\snl^2$ denote the variance of $\pred$ and $\predl{n}$, respectively.
    Then a prediction algorithm is stable if and only if the following two statements hold.
    \begin{enumerate}
        \item The variance stabilizes in the sense that $\limN \sn^2 - \snl^2 = 0$.
        \item The drift added by a new observation vanishes asymptotically in $L_2$, that is,
              \begin{align}\label{eq:nec_OneStepUpdateDrift}
                \limN \E \left( \left(\EC{\pred}{\Tnl{n}, x_{n+1}} - \predl{n}\right)^2 \right) = 0. 
              \end{align}
    \end{enumerate}
\end{lemma}

\Cref{lem:nec_equivalentStability} gives an alternative representation of the stability condition of \Cref{def:jack_oosStability} in $L_{2+\xi}$. 
Next, we will show that the two conditions derived in \Cref{lem:nec_equivalentStability} are not only sufficient but essentially also necessary for the Jackknife to be asymptotically uniformly conditionally conservative.

\begin{proposition}[Necessity of a stabilizing variance in $L_{2+\xi}$]\label{prop:nec_VanishingVariance}
    Let $\pred$ be a 
    symmetric predictor and assume the $(2+\xi)$-th absolute moments of $\pred$, $\predl{n}$ and $y_{n+1}$ are uniformly bounded over $n$ for some $\xi > 0$.
    Furthermore, assume the drift added by a new observation vanishes asymptotically in the sense of \cref{eq:nec_OneStepUpdateDrift}.
    Then the following two statements are equivalent.
    \begin{enumerate}
        \item The Jackknife method provides asymptotically uniformly conditionally conservative prediction intervals, that is,
              for all $\epsilon > 0$ and every $\gam > 0$ we have
              \begin{align*}
                \limN \P \left(\inf_{0 \leq \alpha_1 \leq \alpha_2 \leq 1} \PC{y_{n+1} \in \PJi{\alpha_1, \alpha_2}{\gam}}{\Tn} - (\alpha_2 - \alpha_1) \geq - \epsilon \right) = 1. 
              \end{align*}
        \item The variance stabilizes in the sense that $\limN \sn^2 - \snl^2 = 0$.
    \end{enumerate}
\end{proposition}

In the setting of \Cref{prop:nec_VanishingVariance}, the assumption of a stabilizing variance is a necessary and sufficient condition for the Jackknife prediction intervals to be asymptotically uniformly conditionally conservative.
For the assumption of a vanishing one-step-update drift, we can only provide a slightly weaker statement.
\begin{proposition}\label{prop:nec_OneStepUpdateDrift}
    Let $\pred$ be a symmetric predictor 
    and assume the $(2+\xi)$-th absolute moments of $\pred$, $\predl{n}$ and $y_{n+1}$ are uniformly bounded over $n$ for some $\xi > 0$.
    Assume the \emph{expectation} of the one-step-update drift does not vanish asymptotically, that is,
    \begin{align*}
        \limsup_{n \to \infty} |\E(\pred - \predl{n})| > 0. 
    \end{align*}
    Then, the Jackknife method \emph{does not} provide asymptotically uniformly conditionally conservative prediction intervals, that is, there exists an $\epsilon > 0$ and a $\gam > 0$ such that 
    \begin{align*}
        \liminf_{n \to \infty} \P \left(\inf_{0 \leq \alpha_1 \leq \alpha_2 \leq 1} \PC{y_{n+1} \in \PJi{\alpha_1, \alpha_2}{\gam}}{\Tn} - (\alpha_2 - \alpha_1) \geq - \epsilon \right) < 1. 
    \end{align*}
\end{proposition}

To sum it up, we can show that for symmetric predictors with uniformly bounded $(2+\xi)$-th absolute moments, the stability assumption can be decomposed into two conditions. If the $(2+\xi)$-th absolute moments of $y_{n+1}$ are also uniformly bounded, then these two conditions are sufficient to prove statistical coverage guarantees for the Jackknife. 
One of the two conditions, stabilizing variance, is also necessary, while the assumption of a vanishing one-step-update drift needs to be replaced by a slightly weaker condition to be necessary. Indeed, this leaves the question unanswered whether the Jackknife comes with coverage guarantees if the one-step-update drift does not vanish in $L_2$ while its expectation does.
Nevertheless, we consider the foregoing results as a major step forward towards answering the question whether, for a given algorithm, the stability condition is necessary for the Jackknife to provide conditionally conservative prediction intervals in large samples.

\subsection{Comparison of the Jackknife and the Jackknife+}
In view of the equivalence result given by \Cref{thm:eqi_PacBoundAsymptotic}, one could now ask which prediction interval is shorter. As \Cref{ex:appRem_bothLengthsPossibleJ} and \Cref{ex:appRem_bothLengthsPossibleJp} show, neither method is able to provide shorter prediction intervals than the other for all prediction algorithms. Thus, we restrict our comparison of the two methods to coverage guarantees.

\Cref{sub:necessityStab} shows that the stability condition essentially cannot be removed if we are interested in conditional coverage guarantees of the Jackknife. Since the equivalence result for the Jackknife and the Jackknife+ can only be guaranteed under algorithmic stability, this opens the question whether the Jackknife+ may still provide conditionally conservative prediction intervals in large samples if the underlying prediction algorithm is not stable.
While Theorem $3$ of \textcite{bian2022training} proves the existence of a non-stable algorithm for which the Jackknife+ does not fulfill this property, it remains an open question whether for every arbitrary algorithm the Jackknife+ fails to give conditionally conservative prediction intervals in large samples.

In fact, by a modification of our proofs, we can show that for symmetric predictors the Jackknife+ is still asymptotically uniformly conditionally conservative if the asymptotic out-of-sample stability condition is replaced by the slightly weaker asymptotic swap-stability condition, which is defined as  
\begin{align*}
    \predl{1} - \predl{2} \plim 0. 
\end{align*}
This result is in line with Theorem~4.4 of \textcite{liang2023algorithmic}, who showed that the symmetric Jackknife+ is (pointwise) asymptotically conditionally conservative if the distribution $\mathcal{P}_n$ and $\alpha_1, \alpha_2$ are fixed.

While there are pathological examples of algorithms that are asymptotically swap-stable but not asymptotically out-of-sample stable, evidence so far suggests that the difference between these two stability conditions is small in practice. In fact, \textcite{liang2023algorithmic} showed that any asymptotically swap-stable algorithm can be used to create an asymptotically out-of-sample stable algorithm which is based on the former (cf. Proposition $4.6$ therein).

Nevertheless, we argue that the Jackknife+ does indeed improve over the simple Jackknife for the following reason:
While under algorithmic stability the two approaches are equivalent in large samples as shown in \Cref{thm:eqi_PacBoundAsymptotic}, the Jackknife+ with nominal level $1-\alpha$ guarantees a universal marginal coverage probability of at least $1-2\alpha$ even for non-stable algorithms (cf. Theorem $1$ of \cite{barber2021predictive} for the symmetric case and Appendix A therein for the asymmetric version). In contrast, there exist non-stable algorithms where the marginal coverage probability of the original Jackknife equals $0$ (cf. Theorem $2$ in \cite{barber2021predictive}).
Since the marginal coverage guarantee for the Jackknife+ holds for all algorithms, we expect the Jackknife+ to be preferred over the Jackknife if no theoretical analysis of the underlying prediction algorithm and model is available.

\section{Summary}
In this work, we prove that prediction intervals based on CV are conditionally conservative in large samples if the underlying prediction algorithm is stable and the prediction error is stochastically bounded. If, additionally, the conditional distribution of $y_{n+1}$ given $x_{n+1}$ fulfills a continuity assumption, the CV method even provides conditionally valid prediction intervals in large samples, showing that these prediction intervals are not overly large. 

Furthermore, we present an equivalence result between CV and CV+ under algorithmic stability, showing that neither method is able to outperform the other one in terms of their actual conditional coverage probability. From this, we conclude that the CV+ also provides conditionally conservative prediction intervals in large samples under the same assumptions as original CV.

We discuss that these conditions can essentially not be removed for the classical Jackknife. Combining these findings with the fact that the Jackknife+ variants offer \emph{marginal} coverage guarantees even for non-stable algorithms, we conclude that CV+ indeed improves over the original CV approach. Furthermore, we show that the leave-one-out residuals can be used for a more general uncertainty quantification by allowing to consistently estimate other risk measures such as the conditional mean-squared prediction error or the conditional misclassification error.

Given the aforementioned negative results, we believe the Jackknife alongside with the Jackknife+ to be used and investigated mainly for stable algorithms, in which case they are asymptotically equivalent.
    \printbibliography[heading=bibintoc,title=References]
    \appendix{}
    \section{Cross-validation and CV+}\label{sec:CV}
\subsection{Definitions}
Let $(K_j)_{j=1}^k$ be a partition of $\{1 \ldots, n\}$, that is, $\bigcup_{j=1}^k K_j = \{1 \ldots, n\}$ and $K_i \cap K_j = \emptyset$ for all $i \neq j$. Additionally, we assume $K_j \neq \emptyset$ for all $1 \leq j \leq k$ and write $|K_j|$ to denote the cardinality of the set $K_j$. 
For our asymptotic statements we call $(\mathcal{K}_n)_{n \geq 2}$ a sequence of partitions if for every integer $n \geq 2$ $\mathcal{K}_n = (K_j^{(n)})_{j=1}^{k_n}$ is a partition of $\{1 \ldots, n\}$. Furthermore, we will implicitly assume $k_n \geq 2$ for all $n \in \N$. For the sake of readability, we will suppress the index $n$ in the partitions and in $k_n$ whenever it is clear from the context.

Define a leave-one-fold-out predictor for $y_i$ based on CV as
\begin{align*}
    \hat{y}_i^{\backslash K_j} = \mathcal{A}_{p, n-|K_j|}(x_i, \Tncv{j}),
\end{align*}
for $1 \leq i \leq n+1$ and $1 \leq j \leq k$, where $\Tncv{j}$ denotes the training data $\Tn$ without the observations of the $j$-th fold $K_j$.
Denoting the fold where the $i$-th observation is contained in with $K_{j(i)}$, we call $\hat{u}_i^{CV} = y_i - \hat{y}_i^{\backslash K_{j(i)}}$ the leave-one-fold-out errors.
Furthermore, let $v_i^{CV} = \pred + \hat{u}_i^{CV}$ and $v_i^{CV+} = \hat{y}_{n+1}^{\backslash K_{j(i)}} + \hat{u}_i^{CV}$ the CV and CV+ analogues of $v_i$ and $v_i^+$ defined in \eqref{eq:eqi_defVi} and \eqref{eq:eqi_defViPlus}, respectively. The corresponding distribution functions $G^{CV}$ and $G^{CV+}$ are defined as
\begin{align*}
    G^{CV}(t) &= \dfrac{1}{k} \sum_{j=1}^k \dfrac{1}{|K_j|} \sum_{i \in K_j} \mathds{1}\{v_i^{CV} \leq t\}
    \text{ and } \\
    G^{CV+}(t) &= \dfrac{1}{k} \sum_{j=1}^k \dfrac{1}{|K_j|} \sum_{i \in K_j} \mathds{1}\{v_i^{CV+} \leq t\}.
\end{align*}

For $\gam \in \R$ and $0 \leq \alpha_1 \leq \alpha_2 \leq 1$ we define the $\gam$-distorted prediction intervals with nominal coverage probability of $\alpha_2 - \alpha_1$ based on CV and CV+ as
\begin{nalign}\label{eq:appExt_defineKfoldCvPi}
    \PCv{\alpha_1, \alpha_2}{\gam} &= [\Q{\alpha_1}{G^{CV}} - \gam, \Q{\alpha_2}{G^{CV}} + \gam] \text{ and } \\
    \PCvp{\alpha_1, \alpha_2}{\gam} &= [\Q{\alpha_1}{G^{CV+}} - \gam, \Q{\alpha_2}{G^{CV+}} + \gam].
\end{nalign}
For convenience, we extend this definition to all $\alpha_1, \alpha_2 \in \R$ (cf. \Cref{def:set_quantiles}) and simply write $\PCv{\alpha_1, \alpha_2}{0}$ and $\PCvp{\alpha_1, \alpha_2}{0}$ to denote the non-distorted prediction intervals, that is, $\gam = 0$.
In particular, if $K^{(n)}_j = \{j\}$ for all $1 \leq j \leq n$ and all $n \in \N$, the above definitions coincide with the Jackknife and the Jackknife+, respectively.

\begin{definition}[Generalized definition of stability]\label{def:appExtensions_kFoldStability}
    We call a predictor asymptotically out-of-sample stable with respect to the sequence of partitions 
    $(\mathcal{K}_n)_{n \in \N}$ and the corresponding 
    leave-one-fold-out analogues $(\predl{K^{(n)}_1}, \ldots, \predl{K^{(n)}_{k_n}})$ if
    \begin{align*}
        \limN \dfrac{1}{k_n} \sum_{j=1}^{k_n} \P( | \pred - \predl{K^{(n)}_j} | \geq \epsilon) = 0 
        \text{ for all } \epsilon > 0.
    \end{align*}
    If the partition and its leave-one-fold-out analogues are clear from the context, we will abbreviate the definition by just stating that the predictor is stable. 
\end{definition}

\subsection{Finite sample results for CV and CV+}\label{sub:appRem_finiteSampleResults}
\begin{proposition}[Finite sample guarantees for the CV]\label{prop:appRem_finSamCV}
    Let $\gam > 0, \epsilon > 0$ and $\mu \in \R$.
    For every $L \geq 0$, we then have
    \begin{align*}
        \P &\left[ 
            \inf_{0 \leq \alpha_1 \leq \alpha_2 \leq 1} 
            \PC{y_{n+1} \in \PCv{\alpha_1, \alpha_2}{\gam}}{\Tn} - (\alpha_2 - \alpha_1)
            > - 2\epsilon \right] \\
        \geq&~ 1 - 
        \dfrac{2\P(|y_{n+1} - \pred - \mu| \geq L)}{\epsilon}
        - \dfrac{(8L+12\gam)}{k\gam \epsilon^2}\\
        &- \dfrac{4(5k+1)}{k^2 \gam\epsilon^2} \sum_{j=1}^k \E(\min(2L+3\gam, |\pred - \predl{K_j}|)).
    \end{align*}
    Furthermore,
    \begin{align*}
        &\P \left[ 
            \inf_{0 \leq \alpha_1 \leq \alpha_2 \leq 1} 
            \PC{y_{n+1} \in \PCv{\alpha_1, \alpha_2}{\gam}}{\Tn} - (\alpha_2 - \alpha_1)
            > - 2\epsilon \right] \\
        &\geq 1 - 
        \dfrac{1}{k\gam \epsilon^2} \E( |y_{n+1} - \pred - \mu|)
        - \dfrac{5k+1}{k^2\gam\epsilon^2} \sum_{j=1}^k \E(|\pred - \predl{K_j}|).
    \end{align*}
    The statements remain true if we replace the expression
    \begin{align*}
        \P \left[ 
            \inf_{0 \leq \alpha_1 \leq \alpha_2 \leq 1} 
            \PC{y_{n+1} \in \PCv{\alpha_1, \alpha_2}{\gam}}{\Tn} - (\alpha_2 - \alpha_1)
            > - 2\epsilon \right]
    \end{align*}
    by 
    \begin{align*}
        \P \left[ 
            \sup_{0 \leq \alpha_1 \leq \alpha_2 \leq 1} 
            \PC{y_{n+1} \in \PCv{\alpha_1, \alpha_2}{-\gam}}{\Tn} - (\alpha_2 - \alpha_1)
            <  2\epsilon \right].
    \end{align*}
\end{proposition}

In order to get finite sample statements for the non-distorted prediction intervals in the continuous case 
we can use \Cref{prop:appRem_finSamCV} together with
\begin{align*}
    &\left| \PC{y_{n+1} \in \PCv{\alpha_1, \alpha_2}{\gam}}{\Tn} 
    - \PC{y_{n+1} \in \PCv{\alpha_1, \alpha_2}{0}}{\Tn} \right| \\
    &\leq 
    \E( \min(1, 2|\gam| \supnorm{f_{y_{n+1}\|x_{n+1}}} ) ) \text{ a.s.}
\end{align*}
for every $\gam \in \R$ (see \cref{eq:appExtProofs_InflationBound}).

Recall that for $\alpha_1 > \alpha_2$ the prediction interval $\PCv{\alpha_1, \alpha_2}{0}$ is either the empty set (as one might have anticipated) or a singleton.
In order to avoid this case distinction and additional technicalities in some proofs, the following result is formulated for all $\alpha_1, \alpha_2 \in \R$ rather than for the more usual set $0 \leq \alpha_1 \leq \alpha_2 \leq 1$.
\begin{proposition}[Finite sample relation between CV and CV+]\label{prop:appRem_finSamEqui}
    For every $\epsilon > 0, \gam \geq 0$ and $\kappa \in \R$ we have
    \begin{align*}
        &\P \left[\inf_{\alpha_1, \alpha_2 \in \R} 
        \PC{y_{n+1} \in \PCv{\alpha_1 - \epsilon, \alpha_2 + \epsilon}{\kappa + \gam} }{\Tn} 
        - \PC{y_{n+1} \in \PCvp{\alpha_1, \alpha_2}{\kappa} }{\Tn} \leq -\epsilon \right] \\
        &\leq \dfrac{1}{k \epsilon^2} \sum_{j=1}^k \P(|\pred - \predl{K_j}| > \gam).
    \end{align*}
    Conversely, the following holds true 
    \begin{align*}
        &\P \left[\inf_{\alpha_1, \alpha_2 \in \R} 
        \PC{y_{n+1} \in \PCvp{\alpha_1 - \epsilon, \alpha_2 + \epsilon}{\kappa + \gam}}{\Tn} 
        - \PC{y_{n+1} \in \PCv{\alpha_1, \alpha_2}{\kappa}}{\Tn} \leq -\epsilon \right] \\
        &\leq \dfrac{1}{k \epsilon^2} \sum_{j=1}^k \P(|\pred - \predl{K_j}| > \gam).
    \end{align*}
\end{proposition}
We stress the fact that $\kappa$ may also be negative, which allows for a comparison of shrunken prediction intervals. For example, setting $\kappa = -2\gam$ gives
\begin{align*}
    &\P \left[\sup_{\alpha_1, \alpha_2 \in \R} 
        \PC{y_{n+1} \in \PCvp{\alpha_1, \alpha_2}{-2\gam} }{\Tn}
        - \PC{y_{n+1} \in \PCv{\alpha_1 - \epsilon, \alpha_2 + \epsilon}{-\gam} }{\Tn} 
        \geq \epsilon \right] \\
    &= \P \left[\inf_{\alpha_1, \alpha_2 \in \R} 
        \PC{y_{n+1} \in \PCv{\alpha_1 - \epsilon, \alpha_2 + \epsilon}{-\gam} }{\Tn} 
        - \PC{y_{n+1} \in \PCvp{\alpha_1, \alpha_2}{-2\gam} }{\Tn} \leq -\epsilon \right] \\
        &\leq \dfrac{1}{k \epsilon^2} \sum_{j=1}^k \P(|\pred - \predl{K_j}| > \gam).
\end{align*}

Now, finite sample statements for CV+ can be obtained by combining \Cref{prop:appRem_finSamCV} with \Cref{prop:appRem_finSamEqui}.

\subsection{Asymptotics for CV and CV+}
The finite sample statements in \Cref{prop:appRem_finSamCV} and \Cref{prop:appRem_finSamEqui} imply the following asymptotic results:
\begin{theorem}\label{thm:appExt_kFoldCvAlsoValid}
    Let $(\mathcal{K}_n)_{n \in \N}$ be a sequence of partitions, $0 \leq \alpha_1 \leq \alpha_2 \leq 1$, $\gam \in \R$ and define $\PCv{\alpha_1, \alpha_2}{\gam}$ and 
    $\PCvp{\alpha_1, \alpha_2}{\gam}$ as in \cref{eq:appExt_defineKfoldCvPi}. 
    Suppose the predictor is stable as in \Cref{def:appExtensions_kFoldStability}.
    \begin{thmenum}[label=\alph*)]
        \item \label{it:appExt_asymptoticEquivalence}\emph{Asymptotic equivalence:}
            The $k_n$-fold CV prediction intervals are asymptotically equivalent to their CV+ analogues in the sense that \Cref{thm:eqi_PacBoundAsymptotic} holds with $\PCv{\alpha_1, \alpha_2}{\cdot}$ and $\PCvp{\alpha_1, \alpha_2}{\cdot}$ replacing $\PJi{\alpha_1, \alpha_2}{\cdot}$ and $\PJpi{\alpha_1, \alpha_2}{\cdot}$.
        \item \label{it:appExt_CvValid}\emph{Asymptotic validity of $k_n$-fold CV:} 
            Assume, additionally, that the prediction error $y_{n+1} - \pred$ 
            is stochastically bounded and
            the number $k_n$ of folds diverges with $n$ to infinity.
            Then, \Cref{thm:jack_asymptoticValidPi} holds with $\PCv{\alpha_1, \alpha_2}{\cdot}$ replacing $\PJi{\alpha_1, \alpha_2}{\cdot}$.
    \end{thmenum}
\end{theorem}

The following result shows a close connection between inflated prediction intervals being asymptotically uniformly conditionally conservative and shrunken prediction being asymptotically uniformly conditionally anti-conservative.
\begin{proposition}\label{prop:appJa_ShrunkenEqInflated}
    The following two statements are equivalent.
    \begin{enumerate}[label = V\arabic*)]
        \item\label{it:appJa_Shrunken}
            For all $\gam > 0$, the $\gam$-shrunken prediction intervals based on CV are asymptotically uniformly conditionally anti-conservative,
            that is, for all $\epsilon > 0$ we have
            \begin{align*}
                \limN \P \left( \sup_{0 \leq \alpha_1 \leq \alpha_2 \leq 1} [\PC{y_{n+1} \in \PCv{\alpha_1, \alpha_2}{-\gam}}{\Tn} 
                - (\alpha_2 - \alpha_1)] \leq \epsilon \right) = 1.
            \end{align*}
        \item\label{it:appJa_Inflated}
            For all $\gam > 0$, the $\gam$-inflated prediction intervals based on CV are asymptotically uniformly conditionally conservative,
            that is, for all $\epsilon > 0$ we have
            \begin{align*}
                \limN \P \left( \inf_{0 \leq \alpha_1 \leq \alpha_2 \leq 1} 
                \left[ \PC{ y_{n+1} \in \PCv{\alpha_1, \alpha_2}{\gam} }{\Tn} - (\alpha_2 - \alpha_1) \right] 
                \geq -\epsilon \right) = 1.
            \end{align*}
    \end{enumerate}
    The same holds true if we replace $\PCv{\alpha_1, \alpha_2}{\cdot}$ with $\PCvp{\alpha_1, \alpha_2}{\cdot}$.
\end{proposition}

\begin{remark}\label{rem:appJa_ShrunkenEquivalence}
    Combining \Cref{prop:appJa_ShrunkenEqInflated} with the asymptotic equivalence result of \Cref{thm:appExt_kFoldCvAlsoValid} gives the following statement for shrunken prediction intervals:
    Assume the predictor is stable. Then $\gam$-shrunken prediction intervals based on CV are asymptotically uniformly conditionally anti-conservative for all $\gam > 0$ if and only if $\gam$-shrunken prediction intervals based on CV+ are (for all $\gam > 0$).
\end{remark}

Combining \Cref{prop:appJa_ShrunkenEqInflated} and the two statements of \Cref{thm:appExt_kFoldCvAlsoValid}, we immediately conclude for the CV+:
\begin{corollary}[Asymptotic conservativeness/validity of CV+]\label{cor:appExt_CvPlusValid}
    Let $(\mathcal{K}_n)_{n \in \N}$ be a sequence of partitions, $0 \leq \alpha_1 \leq \alpha_2 \leq 1$, $\gam \in \R$ and define
    $\PCvp{\alpha_1, \alpha_2}{\gam}$ as in \cref{eq:appExt_defineKfoldCvPi}. 
    Suppose the predictor is stable as in \Cref{def:appExtensions_kFoldStability}.
    If the prediction error $y_{n+1} - \pred$ is stochastically bounded and the number $k_n$ of folds diverges with $n$ to infinity, then \Cref{cor:eqi_asymptoticValidPiPlus} holds with $\PCvp{\alpha_1, \alpha_2}{\cdot}$ replacing $\PJpi{\alpha_1, \alpha_2}{\cdot}$.
\end{corollary}

In particular, \Cref{cor:appExt_CvPlusValid} shows that under the assumptions above, also CV+ is asymptotically conditionally conservative in general and asymptotically conditionally valid in the continuous case, that is, if Assumption \CC is fulfilled. 

This complements Theorem $4$ of \textcite{bian2022training} who showed that the actual conditional coverage probability of the CV+ method with nominal coverage probability of $1 - \alpha$ does not undershoot the level $1 - 2 \alpha$ asymptotically as long as each fold has the same size $n/k$ diverging to infinity \emph{without} posing any stability assumption.
We are able to show that indeed the CV+ method does not undershoot the target level $1-\alpha$ (instead of $1-2\alpha$) asymptotically under the assumptions above. Furthermore, our results also show that the CV+ method does not overshoot the nominal level if the continuity assumption \CC is fulfilled. 

However, the difference between the two results is that we require the \emph{number} of folds $k$ to diverge, while \textcite{bian2022training} require the \emph{size} $n/k$ (for $k \geq 2$) to do so while additionally assuming an equal fold size. In particular, Theorem $4$ in \textcite{bian2022training} is only applicable for a sample size $n$ which is divisible by $k \geq 2$, while we do not have this restriction. The fact that the CV+ method indeed provides reasonable prediction intervals regardless of the size $k$ and the ratio $n/k$ is no surprise in view of Theorem $4$ of \textcite{barber2021predictive}, who showed that the \emph{marginal} coverage probability of the CV+ method does not undershoot its nominal level in both scenarios.

\subsection{Upper bounds for the conditional coverage probability}\label{sub:appRem_upperBoundsCovProb}
Recall that only in the continuous case we are able to guarantee to meet the prescribed level $\alpha_2 - \alpha_1$ asymptotically.
Indeed, it may not be possible to achieve every desired coverage probability in the non-continuous case with non-randomized prediction intervals.
In the continuous case, the additional assumption on the boundedness of $\supnorm{f_{y_{n+1} \| x_{n+1}}}$ is needed to ensure that the distribution is ``sufficiently'' smooth. Without that assumption, a continuous distribution can be arbitrarily close to a non-continuous one (with respect to the $L_1$ distance or the \ldnamEnd), which can lead to the same problems as if the underlying distribution actually were non-continuous. 
However, in practice, it might not be known whether the distribution satisfies this continuity assumption. To solve the problem, the following statement can be useful in practice: 
For all $\gam > 0$ and $\epsilon > 0$, we have\footnote{
    \Cref{eq:appRem_computBoundsInflatedPi} can be seen as follows: 
    By a similar argument as in \Cref{cor:ld_generalPiInequality}, the probability in \cref{eq:appRem_computBoundsInflatedPi} can be bounded from below by $\P(2\ld{\gam}(\Fh, \Ft) \leq \epsilon)$, where $\ld{\gam}$ denotes the \ldname introduced in \Cref{sec:ld}. Then, under the assumptions of \Cref{thm:jack_asymptoticValidPi}, we have $\limN \P(2\ld{\gam}(\Fh, \Ft) \leq \epsilon) = 1$ for all $\gam$, $\epsilon > 0$.
}
\begin{nalign}\label{eq:appRem_computBoundsInflatedPi}
    &\P \left( \Fh(\Q{\alpha_2}{\Fh} + 2 \gam) - \Fh( (\Q{\alpha_1}{\Fh} - 2 \gam)-) + \epsilon \geq 
        \PC{ y_{n+1} \in \PJi{\alpha_1, \alpha_2}{\gam} }{\Tn} \geq
        \alpha_2 - \alpha_1 - \epsilon     
    \right)\\
    &\underset{n \to \infty}{\longrightarrow} 1.
\end{nalign}
The crucial point of \cref{eq:appRem_computBoundsInflatedPi} is that 
$\Fh(\Q{\alpha_2}{\Fh} + 2 \gam) - \Fh(( \Q{\alpha_1}{\Fh} - 2 \gam)-)$ is computable by the statistician. Hence, in applications one can compare the value $\Fh(\Q{\alpha_2}{\Fh} + 2 \gam) - \Fh(( \Q{\alpha_1}{\Fh} - 2 \gam)-)$ with $\alpha_2 - \alpha_1$. If the two values are close to each other, one can hope to meet the prescribed target well. Otherwise, there are many leave-one-out residuals close to $\hat{u}_{(\lceil n \alpha_i \rceil )}$ (for $i \in \{1, 2\}$), which could indicate a high density or even a point mass of the prediction error's conditional distribution located at the boundaries of the prediction interval. In that case, the prediction interval has potentially a larger conditional coverage probability than desired at the expense of a larger prediction interval by $2 \gam$.

\subsection{Comparison of the interval length}
Assume the predictor is stable. Then \Cref{thm:eqi_PacBoundAsymptotic} shows that the conditional coverage probability of the prediction intervals based on the Jackknife and the Jackknife+ are closely related to each other. 
In the special case where the continuity assumption \CC holds, we see that the conditional coverage probability of one prediction interval reaches its nominal level asymptotically (for all $0 \leq \alpha_1 \leq \alpha_2 \leq 1$) \emph{if, and only if} the prediction interval of the other method does.
This is the case whenever we additionally assume the prediction error to be bounded in probability (see \Cref{thm:jack_asymptoticValidPi}) or the distribution of $(y_{n+1}, x_{n+1})$ to be fixed with $n$ (cf. Theorem $3.6$ in \cite{liang2023algorithmic}). Thus, in those cases, we would like to use the one with shorter prediction intervals. As the following two examples show, no method provides (asymptotically valid) shorter prediction intervals \emph{for all algorithms} than the other one with respect to the Lebesgue measure $\lambda$.

\begin{example}\label{ex:appRem_bothLengthsPossibleJ}
    Let $\mathcal{A}_{p,n}(x_{n+1}, (y_i,x_i)_{i=1}^n) = \max_{1 \leq i \leq n} (y_i)$.
    Then, for every $n \geq 2$, $p \in \N$, every distribution $\mathcal{P}_n$ of $(y_{n+1}, x_{n+1})$ and
    every $0 \leq \alpha_1 \leq \alpha_2 \leq 1$ we have
    $\lambda \left( \PJpo{\alpha_1, \alpha_2} \right) \leq \lambda \left( \PJo{\alpha_1, \alpha_2} \right)$ almost surely.
    Furthermore, the algorithm is stable since $\frac{1}{n} \sum_{i=1}^n \P(|\pred - \predl{i}| > 0) \leq \frac{1}{n}.$
    If the distribution of $y_{n+1}$ is continuous, then for every $0 \leq \alpha_1 \leq \frac{n-1}{n} < \alpha_2 \leq 1$ 
    we even have $\lambda(\PJpo{\alpha_1, \alpha_2}) < \lambda(\PJo{\alpha_1, \alpha_2})$ almost surely.
\end{example}

\begin{example}\label{ex:appRem_bothLengthsPossibleJp}
    Let $\mathcal{A}_{p,n}(x_{n+1}, (y_i,x_i)_{i=1}^n) = -\max_{1 \leq i \leq n} (y_i)$.
    Then, for every $n \geq 2$, $p \in \N$, every distribution $\mathcal{P}_n$ of $(y_{n+1}, x_{n+1})$ and
    every $0 \leq \alpha_1 \leq \alpha_2 \leq 1$ we have
    $\lambda \left( \PJo{\alpha_1, \alpha_2} \right) \leq \lambda \left( \PJpo{\alpha_1, \alpha_2} \right)$ almost surely.
    Furthermore, the algorithm is stable since $\frac{1}{n} \sum_{i=1}^n \P(|\pred - \predl{i}| > 0) \leq \frac{1}{n}.$
    If the distribution of $y_{n+1}$ is continuous, then for every $0 \leq \alpha_1 \leq \frac{n-1}{n} < \alpha_2 \leq 1$ 
    we even have $\lambda(\PJo{\alpha_1, \alpha_2}) < \lambda(\PJpo{\alpha_1, \alpha_2})$ almost surely.
\end{example}
In the setting of \Cref{ex:appRem_bothLengthsPossibleJ} and \Cref{ex:appRem_bothLengthsPossibleJp}, both methods provide asymptotically valid prediction intervals if, additionally, 
the sequence of distributions $\mathcal{P}_n$ is fixed over $n$ 
or if for every $n \in \N$ the response $y_{n+1}$ is bounded almost surely by a constant $C$ independent of $n$.
However, in \Cref{ex:appRem_bothLengthsPossibleJ} the Jackknife+ provides shorter prediction intervals while in \Cref{ex:appRem_bothLengthsPossibleJp} the Jackknife does.

\section{Extensions and further remarks}\label{sec:appRemarks}
Our proofs mainly rely on the newly introduced concept of the \ldnamEnd, which is presented in \Cref{sec:ld} in detail.
Nevertheless, the usability of the \ldname is not restricted to prediction intervals based on the Jackknife only. To illustrate this, we will present in this section several extensions to other prediction intervals, including the naive approach and full-conformal prediction, and to other risk measures. Before doing so, we start with the symmetric Jackknife variants:

\subsection{The symmetrized Jackknife variants}\label{sub:appExt_symmetrizedVersions}
We define the symmetric versions of the Jackknife and the Jackknife+ as follows:
Let $s_i := \pred + |\hat{u}_i|$ for $1 \leq i \leq n$ and let $G_s$ denote the empirical distribution function (ecdf) corresponding to $s = (s_1, \ldots, s_n)$. Furthermore, for $1 \leq i \leq n$ let $s_i^+ := \predl{i} + |\hat{u}_i|$ and $ G_s^{+}$ denote the Jackknife+ analogues of $s_i$ and $G_s$.
Now the symmetric prediction intervals based on the Jackknife are defined by replacing the quantiles $\Q{\alpha}{G}$ of $G$ in \cref{eq:eqi_PiWithTheJackAlternative} by those of $G_s$. Similarly, we replace the quantiles of $G^+$ by those of $G_s^+$ in \cref{eq:eqi_PiWithTheJackp} for the Jackknife+.

\begin{proposition}\label{prop:appExt_symmetrizedJackknife}
    \Cref{thm:jack_asymptoticValidPi}, \Cref{thm:eqi_PacBoundAsymptotic} and \Cref{cor:eqi_asymptoticValidPiPlus} still hold if we replace the original Jackknife and Jackknife+ by their symmetric variants.
\end{proposition}

In fact, the arguments in the proof of \Cref{prop:appExt_symmetrizedJackknife} can also be applied mutatis mutandis to the general case of cross-validation to show that the results of \Cref{thm:appExt_kFoldCvAlsoValid} extend to their symmetrized versions.

\subsection{Naive approach on fitted errors}\label{sub:appExt_naive}
\begin{proposition}\label{prop:appExt_equivalenceForFittedValues}
    Let $w_i = \pred + y_i - \hat{y}_i$ for $1 \leq i \leq n$ and $H$ be the empirical distribution function of $(w_i)_{i=1}^n$.
    Define a $\gam$-distorted prediction interval based on the \emph{fitted values} for all $\gam \in \R$, 
    $0 \leq \alpha_1 \leq \alpha_2 \leq 1$ as
    \begin{align*}
        PI^{fv}_{\alpha_1, \alpha_2}(\gam) = [\Q{\alpha_1}{H}-\gam, \Q{\alpha_2}{H} + \gam].
    \end{align*}
    If for every $\epsilon > 0$ we have $\limN \frac{1}{n} \sum_{i=1}^n \P(|\hat{y}_i - \hat{y}_i^{\backslash i}| > \epsilon) = 0$, 
    then the prediction intervals based on the fitted values are asymptotically equivalent to prediction intervals based on the Jackknife in the sense that \Cref{thm:eqi_PacBoundAsymptotic} holds with $PI^{fv}_{\alpha_1, \alpha_2}(\gam)$ replacing $\PJpi{\alpha_1, \alpha_2}{\gam}$. 
\end{proposition}

For simplification, we decided to state the asymptotic version only. In view of \Cref{thm:jack_asymptoticValidPi} this also reveals that the \emph{conditional} coverage probability of prediction intervals based on fitted values does not significantly undershoot its nominal level if both $\frac{1}{n} \sum_{i=1}^n \P(|\hat{y}_i - \hat{y}_i^{\backslash i}| > \gam)$ and $\frac{1}{n} \sum_{i=1}^n \P(|\pred - \predl{i}| > \gam)$ are small.
Thus, \Cref{prop:appExt_equivalenceForFittedValues} generalizes Theorem $6$ of \textcite{barber2021predictive}, who showed the same for the \emph{marginal} coverage probabilities of these intervals.
The fact that functions of the fitted values can be related to their Jackknife analogues if $\hat{y}_i$ is close to $\hat{y}_i^{\backslash i}$ is in line with the literature (see the proof of Theorem $5$ in \cite{devroye1979distribution} or Remark $10$ in \cite{bousquet2002stability}).

\begin{remark}[Low-dimensional case]\label{rem:appExt_fittedValuesPseudoTarget}
    Assume the predictor is symmetric and let the distribution of $(y_{n+1}, x_{n+1})$ be fixed over $n$. In particular, this implies the dimension $p$ of $x_{n+1}$ to be fixed. Furthermore, assume that both the fitted value $\hat{y}_1$ and the leave-one-out prediction $\hat{y}_1^{\backslash 1}$ converge in probability to a pseudo-target $f(x_1)$, where $f: \R^p \to \R$ is a measurable function. 
    Then the assumption of \Cref{prop:appExt_equivalenceForFittedValues} is fulfilled and the prediction interval based on fitted values is equivalent to the Jackknife.
    Even more, the algorithm is stable, and the prediction error is bounded in probability. Thus, the Jackknife, the Jackknife+, and the naive approach based on fitted values provide conditionally conservative prediction intervals in large samples. However, predictors typically do not converge to a pseudo-target in high-dimensional asymptotic regimes where $p$ is not negligible to $n$.
\end{remark}

\subsection{Full-conformal prediction}\label{sub:appExt_fullConformal}
Indeed, the \ldname allows us to prove that prediction intervals are asymptotically conditionally conservative (or valid, in the continuous case) if the difference between the distribution of the conditional prediction error and the distribution the prediction interval is based on vanishes. 
If we replace the out-of-sample stability condition by an in-sample stability condition, the latter can be shown to hold true for full-conformal prediction. Thus, our approach also allows to prove that prediction sets based on the full-conformal method are asymptotically uniformly conditionally conservative and, under an additional continuity assumption, even uniformly conditionally valid in large samples. We defer the presentation and the proof of this result to a future publication.

\subsection{Extension to other risk measures}\label{sub:appExt_risk}
In a classification setting, the usage of prediction intervals to quantify uncertainty is debatable, and one is often more interested in consistently estimating the misclassification error \emph{conditional} on the available training data. For this reason, we present an extension of uncertainty quantification to the consistent estimation of more general loss functions of the prediction error. The main result is given by the following statement:
\begin{proposition}\label{prop:ext_lossFunctions}
    Assume the predictor is stable as in \Cref{def:jack_oosStability} and the prediction error $y_{n+1}-\pred$ is stochastically bounded. 
    Let $\ell: \R \to \R$ be a non-decreasing loss-function such that $\ell(|y_{n+1} - \pred|)$ and $\ell(|y_{n+1} - \predl{n}|)$ are uniformly integrable.\footnote{
        Note that the leave-one-out prediction error $y_{n+1} - \predl{i}$ has the same distribution for each $1 \leq i \leq n$. In contrast, the distribution of the stability term $\pred - \predl{i}$ may depend on $i$ unless the predictor is symmetric.
    }
    Then $\EC{\ell(|y_{n+1}-\pred|)}{\Tn}$ can be bounded by the $\epsilon$-shifted leave-one-out plug-in estimator for any $\epsilon > 0$, that is,
    \begin{align*}
        \limN \P \left( \avg{n} \ell(|\hat{u}_i| - \epsilon)) - \epsilon \leq \EC{\ell(|y_{n+1}-\pred|)}{\Tn} \leq \avg{n} \ell(|\hat{u}_i| + \epsilon) + \epsilon \right) = 1.
    \end{align*}
\end{proposition}

\Cref{prop:ext_lossFunctions} should be compared to a related result of \textcite{bousquet2002stability} (cf. Lemma~9 therein).
While \Cref{prop:ext_lossFunctions} comes in a rather technical form, we demonstrate its significance by applying it to the consistent estimation of the (conditional) mean-squared prediction error and the misclassification error.
\begin{corollary}\label{cor:ext_MseConsistency}
    Assume the predictor is stable and the $(2+\xi)$-th absolute moments of the prediction error and its leave-one-out analogue are uniformly bounded over $n$ for some $\xi > 0$. Then the conditional mean-squared prediction error can be consistently estimated by the empirical second moments of the leave-one-out residuals, that is,
    \begin{align*}
        \limN \E \left| \EC{(y_{n+1}-\pred)^2}{\Tn} - \avg{n}  \hat{u}_i^2 \right| = 0.
    \end{align*}
\end{corollary}
Thus, when comparing different predictors in terms of their (conditional) mean-squared prediction error, the second moments of the leave-one-out residuals give a consistent estimate if the underlying algorithms are stable. 
A similar picture arises if we want to estimate the misclassification error of classification algorithms. For this, we assume that the classes are encoded in integers and the predictor outputs only these classes. While the absolute value of the prediction error has no obvious interpretation here, $y_{n+1}$ is misclassified by $\pred$ if and only if the prediction error is non-zero. In fact, under a stability condition, the fraction of non-zero leave-one-out residuals allows for a consistent estimation of the true misclassification error conditional on the training data, as the following result shows.

\begin{corollary}\label{cor:appExt_misclassificationError}
        Let $\mathcal{S} = \{1, \ldots, K\} \subseteq \N$ be a set of finitely many integers encoding the classes such that $y_{n+1}$, $\pred$ and $\predl{n}$ are contained in $\mathcal{S}$ almost surely for every $n \in \N$.
        If a symmetric predictor fulfills $\limN \P(\pred \neq \predl{n}) = 0$, then the fraction of non-zero leave-one-out residuals is a consistent estimator of the conditional misclassification error, that is,
        \begin{align*}
            \limN \E \left| \PC{y_{n+1} \neq \pred}{\Tn} - \avg{n} \mathds{1}\{\hat{u}_i \neq 0\} \right| = 0.
        \end{align*}
\end{corollary}

\subsection{On m-stability}\label{sub:appRem_mStability}
Theorem~$3.6$ of \textcite{liang2023algorithmic} shows that indeed the Jackknife+ provides asymptotically conditionally conservative (or even valid, in the continuous case) prediction intervals if a symmetric algorithm is stable. However, this result crucially relies on the additional assumption that the distribution of $(y_{n+1}, x_{n+1})$ is fixed with $n$ and therefore prohibits any framework where the dimension $p$ of $x_{n+1}$ depends on $n$ including the high-dimensional asymptotic framework. 
In this subsection, we demonstrate that the strategy of the proof of Theorem~$3.6$ in \textcite{liang2023algorithmic} cannot be adapted to accommodate the high-dimensional framework without further argumentation.
Our \Cref{cor:eqi_asymptoticValidPiPlus}, however, by a different proof strategy, shows that the conclusion of Theorem~$3.6$ in \textcite{liang2023algorithmic} is still true, even if the data generating distribution is allowed to depend on $n$, as long as the prediction error is stochastically bounded.

To show why \textcite{liang2023algorithmic} need to rely on fixing the distribution, we need the concept of $m$-stability introduced therein:
Let $m \in \N, n \geq 2$ and let $(y_i^{(n)}, x_i^{(n)})$ be i.i.d. response-feature pairs with $1 \leq i \leq n+m$ distributed according to $\mathcal{P}_n$. Furthermore, let $\mathcal{T}_{n+m-1}^{(n)} = (y_i^{(n)}, x_i^{(n)})_{i=1}^{n+m-1}$ denote the augmented training data.
Recalling that $\Tnl{n} = (y_i^{(n)}, x_i^{(n)})_{i=1}^{n-1}$, 
$m$-stability $\beta_{m,n-1}^{out}(\mathcal{A},\mathcal{P}_n)$ of a prediction algorithm $\mathcal{A}$ is defined as
\begin{align*}
    \beta_{m,n-1}^{out}(\mathcal{A},\mathcal{P}_n) 
        :&= \E
        (|\mathcal{A}_{p_n,n+m-1}(x_{n+m}^{(n)}, \mathcal{T}_{n+m-1}^{(n)}) - \mathcal{A}_{p_n,n-1}(x_{n+m}^{(n)}, \Tnl{n})|).
\end{align*}
In particular, $1$-stability $\beta_{1,n-1}^{out}(\mathcal{A},\mathcal{P}_n)$ coincides with $\E(|\pred - \predl{n}|)$.

Theorem~$3.6$ of \textcite{liang2023algorithmic} shows that, for a fixed distribution $\mathcal{P}$, the Jackknife+ provides asymptotically conditionally conservative prediction intervals provided that $\beta_{1,n-1}^{out}(\mathcal{A},\mathcal{P}) \to 0$ as $n \to \infty$. 
The proof of Theorem~$3.6$ of \textcite{liang2023algorithmic} relies on two ingredients:
(a) A finite sample inequality, Theorem~$3.2$ in that reference, showing that the actual conditional coverage probability is close to the nominal level if 
$\beta_{m,n-1}^{out}(\mathcal{A},\mathcal{P})$ is small and both $m$ and $n$ are large (this inequality should be compared to our \Cref{prop:appRem_finSamCV}).
And (b) the observation that, for any $m$, $\beta_{m,n-1}^{out}(\mathcal{A},\mathcal{P}) \to 0$ if both $\beta_{1,n-1}^{out}(\mathcal{A},\mathcal{P}) \to 0$ \emph{and} $\mathcal{P}$ is fixed (cf. Lemma~$5.2$ in \cite{liang2023algorithmic}).
The following result shows that the last implication no longer holds if $\mathcal{P}$ is allowed to depend on $n$. In particular, in the finite sample bound of Theorem~$3.2$ of \textcite{liang2023algorithmic}, $\beta_{m,n-1}^{out}(\mathcal{A},\mathcal{P}_n)$ can not be controlled via $\beta_{1,n-1}^{out}(\mathcal{A},\mathcal{P}_n)$ in general if $P_n$ depends on $n$.

\begin{lemma}\label{lem:appRem_counterExMStability}
    Let $(p_n)_{n \in \N}$ be any sequence of integers and $L > 0$.
    Then, there is a sequence of distributions $\mathcal{P}_n^x$ of a $p_n$-dimensional random vector $x^{(n)}$ and an algorithm $\mathcal{A}$, such that
    for every sequence of distributions $\mathcal{P}_n^y$ of $y$ and any coupling $\mathcal{P}_n$ of $\mathcal{P}_n^y$ and $\mathcal{P}_n^x$
    \begin{align*}
        \beta_{m,n-1}^{out}(\mathcal{A},\mathcal{P}_n) 
        :&= \E
        (|\mathcal{A}_{p,n+m-1}(x_{n+m}^{(n)}, \mathcal{T}_{n+m-1}^{(n)}) - \mathcal{A}_{p,n-1}(x_{n+m}^{(n)}, \Tnl{n})|)
        = L\mathds{1}\{m \geq 2\}
    \end{align*}
    holds for all $n \in \N$ and $m \in \N$.
    In particular, the algorithm $\mathcal{A}$ is stable 
    but \emph{not} $m$-stable for $m \geq 2$. 
    Moreover, we have
    \begin{align*}
        \widetilde{\beta}^{out}_{m,n-1}(\mathcal{A},\mathcal{P}_n)(L)
        :&= \P
        \left( \left|\mathcal{A}_{p,n+m-1}(x_{n+m}^{(n)}, \mathcal{T}_{n+m-1}^{(n)}) - \mathcal{A}_{p,n-1}(x_{n+m}^{(n)}, \Tnl{n}) \right|
        \geq L \right) \\
        &= \mathds{1}\{m \geq 2\}.
    \end{align*}
    Furthermore, the predictor $\mathcal{A}_{p,n}(x_{n+1}^{(n)}, \Tn)$ is stochastically bounded and the prediction error is stochastically bounded whenever the sequence of distributions $\mathcal{P}_n^y$ of $y_{n+1}$ is tight.
\end{lemma}

Note that in \Cref{lem:appRem_counterExMStability} $\beta_{m,n-1}^{out}(\mathcal{A},\mathcal{P}_n)$ can be arbitrarily large for $m \geq 2$ and therefore Theorem $3.2$ of \textcite{liang2023algorithmic} only gives trivial bounds for the conditional coverage probability, although the algorithm of \Cref{lem:appRem_counterExMStability} is stable. 

One can hope that at least for some algorithms and distributions the $m$-stability can indeed be bounded from above by a multiple of the $1$-stability. In that case one would typically expect the $m$-stability to be much larger than $1$-stability.\footnote{One might even expect $m$-stability to be monotonically increasing in $m$ in practice (see Section $5.1$ in \cite{liang2023algorithmic} for a corresponding statement).}
As the following example shows, the $m$-stability can be $\sqrt{m}$ times larger than $1$-stability for the Ridge even if the distribution is Gaussian.

\begin{lemma}\label{lem:appRem_mStabilityGrowsWithSqM}
    Let $\rho \in [0,1)$, $\lambda > 0, \sigma > 0$ and $M > 0$.
    For every $n \in \N$, let $1 \leq m_n \leq n$ and $1 \leq p_n \leq \rho m_n$. Furthermore, we assume the following triangular array setting:
    For every $n \in \N$ there exists a $\beta_n \in \R^{p_n}$ with $\norm{\beta_n} \leq M$, such that $y_{0}^{(n)} = \beta_n'x_{0}^{(n)} + u_{0}^{(n)}$, where $u_{0}^{(n)} \sim \mathcal{N}(0, \sigma^2)$ and $x_{0}^{(n)} \sim \mathcal{N}(0, I_{p_n})$. Furthermore, assume that $(y_i^{(n)}, x_i^{(n)})$ are independent and identically distributed like $(y_0^{(n)}, x_0^{(n)})$ for all $1 \leq i \leq m_n + n + 1$.
    Define the Ridge estimator as 
    \begin{align*}
        \hat{\beta}^n_\lambda = (X_n'X_n + \lambda n I_{p_n})^{-1}X_n'Y_n,
    \end{align*}
    where $X_n = (x_1^{(n)}, \ldots, x_n^{(n)})'$ is the (random) matrix whose rows consist of $x_i^{(n)}$ 
    and $Y_n = (y_1^{(n)}, \ldots, y_n^{(n)})'$.
    Define the predictor $\mathcal{A}_{p_n, n}(x, \Tn)$ associated to the Ridge estimator as $x'\hat{\beta}^n_\lambda$ for $x \in \R^{p_n}$.
    Then, there exists a constant $C = C(\rho, \lambda, M, \sigma) > 0$ independent of $n$ and $m_n$ such that
    \begin{align*}
        \liminf_{n \to \infty} \dfrac{\beta_{m_n,n}^{out}(\mathcal{A}, \mathcal{P}_n)}{\sqrt{m_n} \beta_{1,n}^{out}(\mathcal{A}, \mathcal{P}_n)} \geq C.
    \end{align*}
\end{lemma}

Thus, although the $m$-stability can be linked to $1$-stability for the Ridge, it scales (at least) with rate $\sqrt{m}$. 
In fact, this complements Proposition 5.4 in \textcite{liang2023algorithmic}, who gave upper bounds for the $m$-stability of the Ridge in a slightly different framework.

\section{The \ldname}\label{sec:ld}
Our coverage guarantees for prediction intervals crucially rely on the new concept of the \ldnamEnd, which will be introduced in this section. 
\ifthenelse{\equal{\version}{arxive}}{%
    We here present the results only and defer the corresponding proofs to \Cref{sec:appLd}.}
{The proofs of most results are trivial and hence omitted with the exception of \Cref{lem:addDis_totalvariationHelpLemma}, \Cref{lem:dis_totalVariation}, and \Cref{lem:appJa_boundingLdSquared}.}

\subsection{Definition and basic properties}
\begin{definition}[\ldnamEnd]\label{def:ld_levyDivergence}
    Let $F$ and $G$ be cumulative distribution functions (cdfs) and $\gam \geq 0$. 
    Define the \ldname between $F$ and $G$ with tolerance parameter $\gam$ as
    \begin{align}\label{eq:ld_levyDivergenceDefinition}
        \LD = \inf\{\epsilon \geq 0: F(t- \gam) - \epsilon \leq G(t) \leq F(t+\gam) + \epsilon \fatir\}.
    \end{align}
\end{definition}

In comparison, the L\'{e}vy metric $L(F,G)$ between $F$ and $G$ is defined as 
    \begin{align*}
        L(F,G) = \inf\{\epsilon \geq 0: F(t- \epsilon) - \epsilon \leq G(t) \leq F(t+\epsilon) + \epsilon \fatir\}.
    \end{align*}
Furthermore, for $\gam = 0$ the \ldname coincides with the Kolmogorov distance between $F$ and $G$.

\begin{lemma}\label{lem:ld_basicProperties}
    Let $F$ and $G$ be cdfs and $\gam \geq 0$. 
    \leavevmode
    \begin{lemmenum}
        \item \label{it:ld_bpInfimumAttained}
            The infimum in \cref{eq:ld_levyDivergenceDefinition} is attained. In particular,
            \begin{align}\label{eq:ld_ldAttainsMinimum}
                F(t- \gam) - \LD \leq G(t) \leq F(t+\gam) + \LD \fatir.
            \end{align}
        \item \label{it:ld_bpSymmetry}
            The \ldname is symmetric, that is, $\ld{\gam}(F,G) = \ld{\gam}(G,F)$.
        \item \label{it:ld_bpMonoCont}
            The function $\gam \mapsto \LD$ is nonincreasing and continuous from the right on $[0, \infty)$. 
            In particular, we have $0 \leq \LD \leq \ld{0}(F,G) = \supnorm{F-G} \leq 1$.
        \item \label{it:ld_bpAltDef}
            The \ldname can be equivalently written as 
            \begin{align*}
                \LD = \sup_{t \in \R} \max( F(t) - G(t+\gam), G(t) - F(t+\gam)).
            \end{align*}
        \item \label{it:ld_bpTriangle}
            Let $H$ be another cdf and $\gam_1, \gam_2 \geq 0$. 
            Then
            \begin{align*}
                \ld{\gam_1 + \gam_2}(F,H) \leq \ld{\gam_1}(F,G) + \ld{\gam_2}(G,H).
            \end{align*}
        \item \label{it:ld_bpLevy}
            The \ldname is connected to the L\'{e}vy metric through
            \begin{align}\label{eq:ld_connectionToLevy}
                \min(\gam, \LD) \leq L(F,G) \leq \max(\gam, \LD) \text{ for every } \gam \geq 0.
            \end{align}
        \item \label{it:ld_bpScaling}
            Let $c > 0$ and let $F(c \cdot)$ and $G(c \cdot)$ 
            denote the functions $t \mapsto F(c t)$ and $t \mapsto G(c t)$, respectively.
            Then the \ldname satisfies
            \begin{align*}
                \LD =\ld{\frac{\gam}{c}}(F(c \cdot), G(c \cdot)),
            \end{align*}
    \end{lemmenum}
\end{lemma}

The \ldname does not satisfy the triangle inequality and is therefore not a metric in general (unless $\gam = 0$). 
To see this, we consider an arbitrary $\gam > 0$ and Dirac distributions at the points $0$, $\gam$ and $2\gam$ denoted by $D_0$, $D_\gam$ and $D_{2\gam}$, respectively. We then have $\ld{\gam}(D_0, D_\gam) = 0 = \ld{\gam}(D_{\gam}, D_{2\gam})$ but $\ld{\gam}(D_0, D_{2\gam}) = 1$.
We would like to emphasize that in view of the foregoing example, the inequality provided by \Cref{it:ld_bpTriangle} is sharp.

From \cref{eq:ld_connectionToLevy} we can immediately see that the \ldname metricizes weak convergence:
\begin{corollary}\label{cor:ld_ldWeakConvergence}
    Let $F$ and $F_n$ be cdfs for $n \in \N$.
    Then $F_n$ converges weakly to $F$ if and only if $\limN \ld{\gam}(F,F_n) = 0$ for all \emph{positive} $\gam > 0$.
\end{corollary}

\subsection{Connection to prediction intervals}
Changing our perspective from distribution functions to their corresponding quantiles, we get the following result, which yields the key argument in our proofs.
\begin{proposition}\label{prop:ld_quantileInequality}
    Let $F$ and $G$ be cdfs and $\gam \geq 0$. 
    We then have
    \begin{align}\label{eq:ld_quantileInequality}
        \Q{\alpha-\LD}{F} - \gam \leq \Q{\alpha}{G} \leq \Q{\alpha + \LD}{F} + \gam \text{ for all } \alpha \in \R.
    \end{align}
\end{proposition}
\Cref{prop:ld_quantileInequality} provides a connection between the quantiles of two distribution functions based on their \ldnamEnd. 
As the following result shows, the inequality in \eqref{eq:ld_quantileInequality} is indeed sharp.
\begin{lemma}\label{lem:ld_ldDefinitionViaQuantiles}
    Let $F$ and $G$ be distribution functions and $\epsilon \geq 0, \gam \geq 0$.
    Then
    \begin{align*}
        \LD \leq \epsilon \Longleftrightarrow \left[\forall \alpha \in \R: \Q{\alpha-\epsilon}{F} - \gam \leq \Q{\alpha}{G} \leq \Q{\alpha + \epsilon}{F} + \gam \right].
    \end{align*}\label{eq:ld_ldDefinitionViaQuantiles}
    In particular, the \ldname can alternatively be defined as
    \begin{align*}
        \LD = \inf\{\epsilon \geq 0: \Q{\alpha-\epsilon}{F} - \gam \leq \Q{\alpha}{G} \leq \Q{\alpha + \epsilon}{F} + \gam \text{ for all } \alpha \in \R \}.
    \end{align*}
\end{lemma}

In view of \Cref{prop:ld_quantileInequality}, one direction of the implication in \Cref{lem:ld_ldDefinitionViaQuantiles} is trivial. However, the interesting statement is the equivalence of the statements.
In order to draw the connection between quantiles and prediction intervals, 
we notice that \Cref{prop:ld_quantileInequality} implies that for every $\alpha_1, \alpha_2 \in \R$ and every $\gam \geq 0$ we have
\begin{align*}
    [\Q{\alpha_1}{F} - \gam, \Q{\alpha_2}{F} + \gam] \supseteq [\Q{\alpha_1 + \LD}{G}, \Q{\alpha_2 - \LD}{G}].
\end{align*}
Iterating this argument gives the following result.
\begin{corollary}\label{cor:ld_generalPiInequality}
    Let $F$ and $G$ be cdfs and $\gam \geq 0$. 
    Then, for every real-valued random variable $X$ and every $\alpha_1, \alpha_2, \mu \in \R$ we have
    \begin{align*}
        &\P_X \left( X \in \mu + [\Q{\alpha_1}{F} - \gam, \Q{\alpha_2}{F} + \gam] \right) \\
        &\geq \P_X \left( X \in \mu + [\Q{\alpha_1 + \LD}{G}, \Q{\alpha_2 - \LD}{G}] \right)\\
        &\geq \P_X \left( X \in \mu + [\Q{\alpha_1 + 2\LD}{F} + \gam, \Q{\alpha_2 - 2\LD}{F} - \gam] \right),
    \end{align*}
    where we use the convention that $[a,b] = \emptyset$ whenever $a > b$.
    Since $X$ only takes values in $\R$, the statements above hold true if we replace $[-\infty, b]$ by $(-\infty, b]$ and $[a, \infty]$ by $[a, \infty)$ for any $a, b \in [-\infty, \infty]$.
\end{corollary}

To put it in other words, \Cref{cor:ld_generalPiInequality} allows us to sandwich the coverage probability of a prediction interval based on a distribution function $G$ by the coverage probabilities of prediction intervals based on another distribution function $F$. Indeed, if the actual coverage probabilities of prediction intervals based on $F$ are close to its nominal level, this transfers to the prediction interval based on $G$ if $\LD$ is small.

\subsection{Connection to expectation}
In this subsection, we draw a connection between the \ldname and expected values. 
\begin{lemma}\label{lem:addDis_totalvariationHelpLemma}
    Let $f: \R \to \R$ be a non-decreasing function and let $X$ and $Y$ be real random variables. 
    Then, for every $\gam \geq 0$ and every $t \in \R$ we have
    \begin{align}\label{eq:addA_totalvariation1}
        \P( f(X) > t) - \P( f(Y + \gam) > t ) \leq \ld{\gam}(F_X, F_Y),
    \end{align}
    where $F_X$ and $F_Y$ denote the distribution functions of $X$ and $Y$, respectively.
\end{lemma}
\ifthenelse{\equal{\version}{arxive}}{}{
    \begin{proof}[Proof of \Cref{lem:addDis_totalvariationHelpLemma}]
	We start with the definition
	\begin{align*}
		f^{-1}(t) := \sup \{a \in \R: f(a) \leq t\} \in [-\infty, \infty], 
	\end{align*}
	where we use the convention $\sup \emptyset = - \infty$,
	and distinguish four cases:
	
	Case $f^{-1}(t) = -\infty$: Here, $f(x) > t$ holds true for all $x \in \R$. Thus, $\P( f(X) > t) = 1 = \P( f(Y + \gam) > t )$ holds true. Now the claim follows as $\ld{\gam}(F_X, F_Y) \geq 0$.
	
	Case  $f^{-1}(t) = \infty$: Since $f$ is non-decreasing, $f(x) \leq t$ holds true for all $x \in \R$. Thus, we have $\P( f(X) > t) = 0 = \P( f(Y + \gam) > t )$, which proves \cref{eq:addA_totalvariation1}.
	
	Case $|f^{-1}(t)| < \infty$ and $f(f^{-1}(t)) \leq t$: We claim $a > f^{-1}(t) \Leftrightarrow f(a) > t$ holds true or, equivalently, $a \leq f^{-1}(t) \Leftrightarrow f(a) \leq t$.
	To prove this, we assume $a \leq f^{-1}(t)$. By monotonicity, we have $f(a) \leq f(f^{-1}(t))$, where the right-hand side is not larger than $t$ in this case.
	For the reverse direction, we start with an $a$ fulfilling $f(a) \leq t$. Thus, $a \in \{ x \in \R: f(x) \leq t\}$ and therefore $a \leq \sup \{ x \in \R: f(x) \leq t\} = f^{-1}(t)$.
	Hence, we have
	\begin{align*}
		\P( f(X) > t) - \P( f(Y + \gam) > t ) & = \P( X > f^{-1}(t) ) - \P(Y > f^{-1}(t) - \gam)                     \\
		                                      & = F_Y(f^{-1}(t) - \gam) - F_X( f^{-1}(t) ) \leq \ld{\gam}(F_X, F_Y), 
	\end{align*}
	where the last inequality is given by \Cref{it:ld_bpInfimumAttained}.
	
	Case $|f^{-1}(t)| < \infty$ and $f(f^{-1}(t)) > t$: We claim that $a \geq f^{-1}(t) \Leftrightarrow f(a) > t$.
	By monotonicity, $a \geq f^{-1}(t)$ implies $f(a) \geq f(f^{-1}(t))$, which in this case is strictly larger than $t$.
	For the other direction, we start with an $a$ fulfilling $f(a) > t$.
	Thus, $a$ is not contained in the set $\{x \in \R: f(x) \leq t\}$. Recalling that $f$ is non-decreasing, we conclude $a \geq \sup \{ x \in \R: f(x) \leq t\} = f^{-1}(t)$.
	Hence, we have
	\begin{align*}
		\P( f(X) > t) - \P( f(Y + \gam) > t ) & = \P( X \geq f^{-1}(t) ) - \P(Y \geq f^{-1}(t) - \gam)                                                    \\
		                                      & = \lim_{\epsilon \searrow 0} \left[ F_Y(f^{-1}(t) - \gam - \epsilon) - F_X( f^{-1}(t) - \epsilon) \right] \\
		                                      & \leq \ld{\gam}(F_X, F_Y),                                                                                 
	\end{align*}  
	where, again, the last inequality is given by \Cref{it:ld_bpInfimumAttained}.
\end{proof}

}

\begin{lemma}\label{lem:dis_totalVariation}
    Let $f: \R \to [M_1, M_2]$ be a non-decreasing function for some $M_1, M_2 \in \R$, $X$ and $Y$ be real random variables and 
    let $F_X$ and $F_Y$ denote the distribution functions of $X$ and $Y$, respectively. 
    We then have for every $\gam \geq 0$
    \begin{align}
        \E( f(X) ) &\geq \E( f(Y - \gam) ) - (M_2 - M_1) \ld{\gam}(F_X, F_Y) \text{ as well as } \label{eq:dis_totalvariationLower1}\\
        \E( f(X) ) &\leq \E( f(Y + \gam) ) + (M_2 - M_1) \ld{\gam}(F_X, F_Y). \label{eq:dis_totalvariationUpper1}
    \end{align}
    If, additionally, $f$ is Lipschitz continuous with constant $L$ this yields the inequality
    \begin{align}\label{eq:dis_totalvariationLipschitz}
         \left| \E( f(X) ) - \E( f(Y) ) \right| \leq L \gam + (M_2 - M_1) \ld{\gam}(F_X, F_Y).
    \end{align}
    Furthermore, if $g: \R \to \R$ is a function of finite total variation $V_{-K}^K(g)$ on any interval $[-K, K]$ with $K > 0$, we also get
    \begin{align}\label{eq:dis_totalKoksma}
        \left| \E( g(X) ) - \E( g(Y) ) \right| \leq \sup_{K > 0} V_{-K}^K(g) \ld{0}(F_X, F_Y)
    \end{align}
    whenever $\sup_{K > 0} V_{-K}^K(g)$ is finite.
\end{lemma}
\ifthenelse{\equal{\version}{arxive}}{}{
    \begin{proof}[Proof of \Cref{lem:dis_totalVariation}]
	We start showing \cref{eq:dis_totalvariationUpper1}.
	W.l.o.g. we assume that $M_1 = 0$, as otherwise we could shift $f$ by $-M_1$, and let $M$ denote $M_2 - M_1$. 
	As $f$ is bounded and nonnegative, we get
	\begin{align*}
		\E( f(X) ) - \E( f(Y + \gam) ) & = \int_0^M \P( f(X) > t) - \P( f(Y + \gam) > t) dt            \\
		                               & \leq \int_0^M \ld{\gam}(F_X, F_Y) dt = M \ld{\gam}(F_X, F_Y), 
	\end{align*}
	where the inequality in the preceding display follows from \Cref{lem:addDis_totalvariationHelpLemma}.
	With the same arguments as before (or by applying the results above to $\tilde{X} = Y - \gam$ and $\tilde{Y} = X - \gam$) we also get $\E( f(X) ) - \E( f(Y - \gam) ) \geq -M \ld{\gam}(F_X, F_Y)$, proving \cref{eq:dis_totalvariationLower1}.
	
	\Cref{eq:dis_totalvariationLipschitz} can be directly derived from \cref{eq:dis_totalvariationLower1} and \cref{eq:dis_totalvariationUpper1} as the Lipschitz-continuity implies
	$\E(f(Y-\gam)) \geq \E(f(Y)) - L \gam$ and $\E(f(Y+\gam)) \leq \E(f(Y)) + L \gam$.
	
	In order to prove \cref{eq:dis_totalKoksma} we would like to point out that we already showed $|E(f(X)) - \E(f(Y))| \leq M \ld{0}(F_X, F_Y)  = M \supnorm{F_X - F_Y}$ to hold true for $f$ as above.
	If $\sup_{L > 0}V_{-L}^L(g) < \infty$, we can find for every $\epsilon > 0$ a $K > 0$, such that the total variation $\sup_{L > 0}V_{-L}^L(g) \leq V_{-K}^K(g) + \epsilon$.
	From this it follows that we have $|g(x) - g(K)| \leq \epsilon$ for every $x \geq K$ and $|g(x) - g(-K)| \leq \epsilon$ for every $x \leq -K$.
	We define a function $g_{K}(x): \R \to \R$ as
	\begin{align*}
		g_{K}(x) = \begin{cases}
		g(x)  & \text{ if } x \in [-K, K], \\
		g(-K) & \text{ if } x < -K \text{ and }       \\
		g(K)  & \text{ if } x > K.        
		\end{cases}
	\end{align*}
	We then have $\supnorm{g_{K} - g} \leq \epsilon$, implying $|\E(g(X)) - \E(g_{K}(X))| \leq \epsilon$ as well as 
	$|\E(g(Y)) - \E(g_{K}(Y))| \leq \epsilon$. Moreover, we have $V_{-K}^K(g_{K}) = V_{-K}^K(g) \leq \sup_{L > 0}V_{-L}^L(g)$. 
	    
	Now, as the total variation of $g_{K}$ on $[-K, K]$ is finite, we can split it into two non-decreasing functions $g_K^{+}$ and $g_K^{-}$, such that $g_{K}(x) = g_K^{+}(x) - g_K^{-}(x)$ for all $x$ in $[-K, K]$ (cf.~Section~31 of Chapter~6 in \cite{billingsley1995probability}). To be more precise, we define
	\begin{align*}
		g_K^+(x) & = \frac{1}{2}(V_{-K}^{x}(g_K) + g_K(x)) \\
		g_K^-(x) & = \frac{1}{2}(V_{-K}^{x}(g_K) - g_K(x)) 
	\end{align*}
	where $V_{-K}^{x}(g_K)$ denotes the total variation of $g_K$ on the interval $[-K,x]$.
	We extend the definition of $g_K^+$ and $g_K^-$ to all $x \in \R$ in the sense that
	$g_K^{+}(x) = g_K^{+}(K)$ for all $x > K$, $g_K^{+}(x) = g_K^{+}(-K)$ for all $x < -K$ (and $g_K^{-}$ is defined outside of $[-K,K]$ analogously).
	Thus, the equality $g_{K}(x) = g_K^{+}(x) - g_K^{-}(x)$ even holds for \emph{all} $x \in \R$. 
	As the total variation of $g_{K}$ is finite, the functions $g_K^{+}$ and $g_K^{-}$ are also bounded. 
	Furthermore, we have $V_{-K}^K(g_{K}) = g_K^{+}(K) - g_K^{+}(-K) + g_K^{-}(K) - g_K^{-}(-K)$.
	By the monotonicity of $g_K^{+}$ and $g_K^{-}$ the latter term can be expressed equivalently as $\sup_{x \in \R}g_K^{+}(x) - \inf_{x \in \R}g_K^{+}(x) + \sup_{x \in \R}g_K^{-}(x) - \inf_{x \in \R}g_K^{-}(x)$.
	Applying \cref{eq:dis_totalvariationLower1} and \cref{eq:dis_totalvariationUpper1} with $\gam = 0$ to the functions $g_K^{+}$ and $g_K^{-}$ we conclude
	\begin{align*}
		\left| \E(g_{K}(X)) - \E(g_{K}(Y)) \right| 
		  & \leq \left| \E(g_K^{+}(X)) - \E(g_K^{+}(Y)) \right| + \left| \E(g_K^{-}(X)) - \E(g_K^{-}(Y)) \right|                                                  \\
		  & \leq \supnorm{F_X - F_Y} \left( \sup_{x \in \R}g_K^{+}(x) - \inf_{x \in \R}g_K^{+}(x) + \sup_{x \in \R}g_K^{-}(x) - \inf_{x \in \R}g_K^{-}(x) \right) \\
		  & = V_{-K}^K(g_{K}) \supnorm{F_X - F_Y}.                                                                                                                
	\end{align*}
	Putting the pieces together, we end up with
	\begin{align*}
		\left| \E(g(X)) - \E(g(Y)) \right| 
		  & \leq \left| \E(g_{K}(X)) - \E(g_{K}(Y)) \right| + 2 \epsilon    \\
		  & \leq V_{-K}^K(g_{K}) \supnorm{F_X - F_Y} + 2 \epsilon           \\
		  & \leq \sup_{L > 0} V_{-L}^L(g) \supnorm{F_X - F_Y} + 2 \epsilon. 
	\end{align*}
	As $\epsilon > 0$ can be made arbitrarily small, \cref{eq:dis_totalKoksma} holds true.
\end{proof}
}

\Cref{lem:dis_totalVariation} allows us to approximate $\E(f(X))$ by $\E(f(Y \pm \gam))$ for certain functions $f$ if $\ld{\gam}(F_X, F_Y)$ is small. In fact, \cref{eq:dis_totalKoksma} is related to Koksma's inequality (cf. \cite{kuipers1974uniform}), but allows \emph{both} real random variables to follow arbitrary distributions.

\subsection{Bounds}
Since the \ldname between two distribution functions may be hard to compute, the following results provide upper bounds.
\begin{lemma}\label{lem:ld_WassersteinTypeBound}
    Let $F$ and $G$ be cdfs and $\gam \geq 0$.
    Let $\Gamma(\mu_F,\mu_G)$ be the set of all couplings between the measures $\mu_F$ and $\mu_G$ corresponding to the distribution functions $F$ and $G$.
    Then
    \begin{align}\label{eq:ld_WassersteinTypeBound}
        \LD &\leq \inf_{\substack{X \sim \mu_F, Y \sim \mu_G,\\ \P \in \Gamma(\mu_F,\mu_G)}} \P( |X-Y| > \gam).
    \end{align}
    In the case $\gam > 0$ this implies
    \begin{align*}
        \LD \leq \gam^{-p} \inf_{\substack{X \sim \mu_F, Y \sim \mu_G,\\ \P \in \Gamma(\mu_F,\mu_G)}} \E(|X-Y|^p),
    \end{align*}
    for any $p > 0$. In particular, for $p \geq 1$ this yields the inequality
    $\LD \leq \gam^{-p}\mathcal{W}_p(F,G)^p$,
    where $\mathcal{W}_p(F,G)$ denotes the $p$-th Wasserstein-distance between $F$ and $G$.
\end{lemma}

To sum it up, we can control the \ldname between two distribution functions if we can control the corresponding Wasserstein-distance. 
Moreover, since in the case of real-valued random variables the $L_1$ distance between two distribution functions coincides with the Wasserstein distance $\mathcal{W}_1$ (cf. \cite{panaretos2019statistical}), \Cref{lem:ld_WassersteinTypeBound} even yields a statement concerning the $L_1$ distance. In fact, it is not difficult to show that this extends to all $p \in (0, \infty)$: For every $\gam > 0$ and $p \in (0,\infty)$ we have
$\LD \leq \gam^{-1/p} \lpnorm{F-G}{p}$.
One particular result we will need in our proofs is the following lemma containing the foregoing statement for the special case $p = 2$.
\begin{lemma}\label{lem:appJa_boundingLdSquared}
    Let $F$ and $G$ be cdfs, $\gam > 0$, $K \geq 0$ and $\mu \in \R$. We then have the following bounds for the \ldname between $F$ and $G$.
    \begin{align*}
        &\LD \leq 1 - F(\mu + K) + F(\mu - K) + \left[ \dfrac{1}{\gam} \int_{[\mu - K-\gam, \mu + K+2\gam]} |F(x) - G(x)|^2 d\lambda(x) \right]^{\frac{1}{2}},\\
        &\LD^2 \leq \dfrac{1}{\gam} \int_{\R} |F(x) - G(x)|^2 d\lambda(x).
    \end{align*}
\end{lemma}

\ifthenelse{\equal{\version}{arxive}}{}{
    \begin{proof}[Proof of \Cref{lem:appJa_boundingLdSquared}]
    It is easy to see that it suffices to show the Lemma with $\mu = 0$ and define $F_{\mu}(t) = F(t+\mu)$ and $G_\mu(t) = G(t+\mu)$ to prove the general case.
	We start showing that 
	\begin{align}\label{eq:appJa_boundingLdSquaredEq1}
		\LD \leq 1 - F(K) + F(-K) + \sup_{|t| \leq K + \gam} \inf_{x \in [t, t + \gam]} |F(x) - G(x)|. 
	\end{align}
	We notice that for every $t \in \R$ and every $x \in [t, t + \gam]$ we have
	\begin{align*}
		\max(G(t) - F(t + \gam), F(t) - G(t + \gam))
		  & \leq \max(G(x) - F(x), F(x) - G(x)) \\
		  & = |G(x) - F(x)|.                    
	\end{align*}
	In particular, this yields
	\begin{align*}
		\ld{\gam}(F, G) 
		&= \sup_{t \in \R} \max(G(t) - F(t + \gam), F(t) - G(t + \gam))\\
		&\leq \sup_{t \in \R} \inf_{x \in [t, t+\gam]} |G(x) - F(x)|.
	\end{align*}
	We then have
	\begin{align*}
		\ld{\gam}(F, G) 
		  & \leq \sup_{t \in \R} \inf_{x \in [t, t+\gam]} |F(x) - G(x)|     \\
		\leq \max \bigg(
		  & \sup_{t < -K-\gam} \inf_{x \in [t, t+\gam]} |F(x) - G(x)|,      
		\sup_{t > K+\gam} \inf_{x \in [t, t+\gam]} |F(x) - G(x)|,\\
		  & \sup_{|t| \leq K + \gam} \inf_{x \in [t, t+\gam]} |F(x) - G(x)| 
		\bigg).
	\end{align*}
	
	We start bounding the first expression in the maximum:
	\begin{align*}
		  & \sup_{t < -K-\gam} \inf_{x \in [t, t+\gam]} |F(x) - G(x)| 
		\leq \sup_{t < -K-\gam} \inf_{x \in [t, t+\gam]} \max(F(x), G(x)) \\
		  & \leq \inf_{x \in [-K-\gam, -K]} \max(F(x), G(x))          
		\leq \inf_{x \in [-K-\gam, -K]} F(x) + |F(x) - G(x)|\\
		  & \leq F(-K) + \inf_{x \in [-K-\gam, -K]} |F(x) - G(x)|     
		\leq F(-K) + \sup_{|t| \leq K + \gam} \inf_{x \in [t, t+\gam]} |F(x) - G(x)|.
	\end{align*}
	Using the equality $|F(x) - G(x)| = |(1 - F(x)) - (1 - G(x))|$ we can show with a similar argument that 
	\begin{align*}
		\sup_{t > K+\gam} \inf_{x \in [t, t+\gam]} |F(x) - G(x)|                         
		\leq 1 - F(K) + \sup_{|t| \leq K + \gam} \inf_{x \in [t, t+\gam]} |F(x) - G(x)|. 
	\end{align*}
	Putting the pieces together, we end up with
	\begin{align*}
		\ld{\gam}(F, G) 
		  \leq&~ \max \bigg(                                                                       
		F(-K) + \sup_{|t| \leq K + \gam} \inf_{x \in [t, t+\gam]} |F(x) - G(x)|,\\
		  &\quad \quad \sup_{|t| \leq K + \gam} \inf_{x \in [t, t+\gam]} |F(x) - G(x)|,                          
		1 - F(K) + \sup_{|t| \leq K + \gam} \inf_{x \in [t, t+\gam]} |F(x) - G(x)|\bigg)\\
		  \leq&~ F(-K) + 1 - F(K) + \sup_{|t| \leq K + \gam} \inf_{x \in [t, t+\gam]} |F(x) - G(x)|), 
	\end{align*}
	which proves \cref{eq:appJa_boundingLdSquaredEq1}.
	In the next step, we will show that 
	\begin{align}\label{eq:appJa_boundingLdSquaredEq2}
		\left[ \sup_{|t| \leq K + \gam} \inf_{x \in [t, t+\gam]} |F(x) - G(x)|) \right]^2 
		\leq \dfrac{1}{\gam} \int_{[-K-\gam, K+2\gam]} |F(x) - G(x)|^2 d\lambda(x).       
	\end{align}
	In order to see this, we proceed as follows:
	\begin{align*}
		\left[ \sup_{|t| \leq K + \gam} \inf_{x \in [t, t+\gam]} |F(x) - G(x)| \right]^2
		  & = \sup_{|t| \leq K + \gam} \inf_{x \in [t, t+\gam]} |F(x) - G(x)|^2                          \\
		  & \leq \sup_{|t| \leq K + \gam} \dfrac{1}{\gam} \int_{[t, t+\gam]} |F(x) - G(x)|^2 d\lambda(x) \\
		  & \leq \dfrac{1}{\gam} \int_{[-K-\gam, K+2\gam]} |F(x) - G(x)|^2 d\lambda(x),                  
	\end{align*}
	which shows \cref{eq:appJa_boundingLdSquaredEq2}.
	    
	Thus, the first inequality of \Cref{lem:appJa_boundingLdSquared} is a direct combination of
	\cref{eq:appJa_boundingLdSquaredEq1} with \cref{eq:appJa_boundingLdSquaredEq2}.
	
	For the last inequality in \Cref{lem:appJa_boundingLdSquared}, we use the first inequality therein and let $K \to \infty$.
\end{proof}
}

Additionally, we will need upper bounds for the \ldname between two empirical distribution functions several times. 
For this, we present the following result, which is a consequence of 
\cref{eq:ld_WassersteinTypeBound}.
\begin{corollary}\label{cor:ld_BoundForEcdfs}
    Let $n \in \N$ and $a_i, b_i \in \R$ and for $1 \leq i \leq n$. 
    Furthermore, let $p_i \geq 0$ for all $1 \leq i \leq n$ such that $\sum_{i=1}^n p_i = 1$.
    Let $F_a$ and $F_b$ denote the weighted distribution functions of $a = (a_i)_{i=1}^n$ and $b = (b_i)_{i=1}^n$ with respect to $(p_i)_{i=1}^n$, that is,
    \begin{align*}
        F_a(t) = \sum_{i=1}^n p_i \mathds{1}\{a_i \leq t\}
        \text{ and }
        F_b(t) = \sum_{i=1}^n p_i \mathds{1}\{b_i \leq t\}.
    \end{align*}
    We then have for every $\gam \geq 0$
    \begin{align*}
        \ld{\gam}(F_a, F_b) \leq \sum_{i=1}^n p_i \mathds{1}\{|a_i - b_i | > \gam\}.
    \end{align*}
\end{corollary}

\section{Proofs for \Cref{sec:CV}}\label{sec:appEqi}
Since the results for the Jackknife and the Jackknife+ are special cases of the results for the CV and CV+, we present the proofs of \Cref{sec:CV} here. Then, the statements of \Cref{thm:jack_asymptoticValidPi} and  \Cref{thm:eqi_PacBoundAsymptotic} can be derived from \Cref{thm:appExt_kFoldCvAlsoValid} by setting $k_n = n$.

\subsection{Auxiliary Lemmas}
For our proofs, we need some auxiliary lemmas. They are formulated for an arbitrary Borel set $A$, but we will apply them only to the case where $A$ is either a bounded interval or the whole real line. Thus, in the latter case we have $\lambda(A) = \lambda(\R) = \infty$ and therefore $\min(\lambda(A),s) = s$ for all $s \in \R$.

\begin{lemma}\label{lem:appJa_expIntIndFunctions}
    Let $A \in \mathcal{B}(\R)$ be a Borel set and $X$ and $Y$ be random variables.
    We then have
    \begin{align*}
        \int_A \E_{X,Y} \left( |\mathds{1}\{X \leq t\} - \mathds{1}\{Y \leq t\}| \right) d\lambda(t)
        \leq \E_{X,Y}\left( \min(\lambda(A), |X-Y| )\right),
    \end{align*}
    where $\lambda$ denotes the Lebesgue-measure.
\end{lemma}
\begin{proof}
    Using Tonelli's theorem, we get
    \begin{align*}
        &\int_A \E_{X,Y} \left( |\mathds{1}\{X \leq t\} - \mathds{1}\{Y \leq t\}| \right) d\lambda(t)
        = \E_{X,Y} \left( \int_A  |\mathds{1}\{X \leq t\} - \mathds{1}\{Y \leq t\}| d\lambda(t)\right) \\
        &\leq \E_{X,Y} \left( \min(\lambda(A), \int_A \mathds{1}\{\min(X,Y) \leq t < \max(X,Y)\} d\lambda(t) \right)\\
        &\leq \E_{X,Y} \left( \min(\lambda(A), \max(X,Y) - \min(X,Y) ) \right).
    \end{align*}
    To finish the proof, we use the fact that $\max(X,Y) - \min(X,Y) = |Y - X|$.
\end{proof}

The following result is essentially contained in \textcite{steinberger2020conditional} and is adapted to our specific situation.
\begin{lemma}\label{lem:appJa_modifyingSteinLeebsLemma}
    Let $\Ft(t) = \PC{y_{n+1} - \pred \leq t}{\Tn}$,
    $(K_j)_{j=1}^k$ a partition of $\{1, \ldots, n\}$ into nonempty sets and define
    \begin{align*}
        \Fh(t) = \frac{1}{k} \sum_{j=1}^k \dfrac{1}{|K_j|} \sum_{i \in K_j} 
        \mathds{1}\{\hat{u}^{CV}_i \leq t\}.
    \end{align*}
    Then
    \begin{align*}
        \E(|\Fh(x) - \Ft(x)|^2)
        \leq&~ 
        \frac{5}{k} \sum_{j=1}^k \E( |\mathds{1}\{y_{n+1} - \pred \leq x\} - \mathds{1}\{y_{n+1} - \predl{K_j} \leq x \}|) \\
        &+ \dfrac{1}{2k}\dfrac{1}{k(k-1)} \sum_{l_1 \neq l_2} \dfrac{1}{|K_{l_1}| |K_{l_2}|} 
        \sum_{i_1 \in K_{l_1}} \sum_{i_2 \in K_{l_2}} \E(|\mathds{1}\{\hat{u}_{i_1}\leq x\} - \mathds{1}\{\hat{u}_{i_2}\leq x\}|).
    \end{align*}
\end{lemma}
\begin{proof}
    In order to find an upper bound for $\E(|\Fh(x) - \Ft(x)|^2)$ we need 
    a small modification of Lemma~C.3 in \textcite{steinberger2020conditional}, which itself is a consequence of Lemma~C.1 therein. 
    As in Lemma~C.3 of \textcite{steinberger2020conditional}, our loss function is $L(f,(x,y)) = \mathds{1}\{y - f(x)\leq s\}$ and therefore bounded by $C=1$.
    We modify Lemma~C.1 by \emph{not} applying Lemma~D.1 in \textcite{steinberger2020conditional} to get an upper bound, but rather use the inequalities derived for $\E(R^2_{CV})$ in Lemma~C.1 directly.
    To be more precise, we do not bound the expression
    \begin{align*}
        \Xi(x) := \dfrac{1}{k}\dfrac{1}{k(k-1)} \sum_{l_1 \neq l_2} 
        \E \left[ \left( \dfrac{1}{|K_{l_1}|} \sum_{i_1 \in K_{l_1}} 
        \mathds{1}\{\hat{u}_{i_1}\leq x\} \right) 
        \left(1 - \dfrac{1}{|K_{l_2}|} \sum_{i_2 \in K_{l_2}}  \mathds{1}\{\hat{u}_{i_2}\leq x\} \right) \right]
    \end{align*}
    by $\dfrac{1}{4(k-1)}$ and, instead, get the following inequality from Lemma~C.1 in \textcite{steinberger2020conditional}
    \begin{align}\label{eq:appJa_modifyingSteinLeebsLemmaEq1}
        \E(|\Fh(x) - \Ft(x)|^2)
        &\leq \Xi(x) + \frac{5}{k} \sum_{j=1}^k \E( |\mathds{1}\{y_{n+1} - \pred \leq x\} - \mathds{1}\{y_{n+1} - \predl{K_j} \leq x \}|).
    \end{align}

    Abbreviating $\mathds{1}\{\hat{u}_{a}\leq x\} \left(1 - \mathds{1}\{\hat{u}_{b}\leq x\} \right)$ with $h_{a, b}(x)$, 
    we can rewrite $\Xi(x)$ as
    \begin{align*}
        \dfrac{1}{k}\dfrac{1}{k(k-1)} \sum_{l_1 \neq l_2} \dfrac{1}{|K_{l_1}| |K_{l_2}|} 
        \sum_{i_1 \in K_{l_1}} \sum_{i_2 \in K_{l_2}}
        \E(h_{i_1, i_2}(x)).
    \end{align*}
    It is easy to see that
    \begin{align*}
        \Xi(x) = \dfrac{1}{2k}\dfrac{1}{k(k-1)} \sum_{l_1 \neq l_2} \dfrac{1}{|K_{l_1}| |K_{l_2}|} 
        \sum_{i_1 \in K_{l_1}} \sum_{i_2 \in K_{l_2}} \E(h_{i_1, i_2}(x) + h_{i_2, i_1}(x)).
    \end{align*}
    To finish the proof, note that $h_{i_1, i_2}(x) + h_{i_2, i_1}(x) = |\mathds{1}\{\hat{u}_{i_1}\leq x\} - \mathds{1}\{\hat{u}_{i_2}\leq x\}|$ and use \cref{eq:appJa_modifyingSteinLeebsLemmaEq1}.
\end{proof}

\begin{lemma}\label{lem:appJa_boundingLdMiddlePart}
    Let $A \in \mathcal{B}(\R)$ be a Borel set and define $\Ft$ and $\Fh$ as in \Cref{lem:appJa_modifyingSteinLeebsLemma}.
    We then have
    \begin{align*}
        \E \left( \int_{A} |\Fh(x) - \Ft(x)|^2 d\lambda(x) \right) 
        \leq&~ \dfrac{1}{k} \E(\min(\lambda(A), |y_{n+1} - \pred|)) \\
        &+ \frac{5k+1}{k^2} \sum_{j=1}^k \E( \min(\lambda(A), |\pred - \predl{K_j}|)).
    \end{align*}
\end{lemma}
\begin{proof}
    Using Tonelli's theorem, we have
    \begin{align*}
        \E \left( \int_{A} |\Fh(x) - \Ft(x)|^2 d\lambda(x) \right) 
        = \int_{A} \E(|\Fh(x) - \Ft(x)|^2) d\lambda(x).
    \end{align*}

    Applying \Cref{lem:appJa_modifyingSteinLeebsLemma} yields
    \begin{align*}
        \int_{A}& \E(|\Fh(x) - \Ft(x)|^2) d\lambda(x) \\
        \leq&~ \frac{5}{k} \sum_{j=1}^k \int_{A}\E( |\mathds{1}\{y_{n+1} - \pred \leq x\} - \mathds{1}\{y_{n+1} - \predl{K_j} \leq x \}|) d\lambda(x)\\
        &+ \dfrac{1}{2k}\dfrac{1}{k(k-1)} \sum_{l_1 \neq l_2} \dfrac{1}{|K_{l_1}| |K_{l_2}|} 
        \sum_{i_1 \in K_{l_1}} \sum_{i_2 \in K_{l_2}} \int_{A} \E(|\mathds{1}\{\hat{u}_{i_1}\leq x\} - \mathds{1}\{\hat{u}_{i_2}\leq x\}|) d\lambda(x).
    \end{align*}

    We now treat the two terms in the upper bound separately.
    First, we apply \Cref{lem:appJa_expIntIndFunctions} to get
    \begin{align*}
        \int_A \E(|\mathds{1}\{\hat{u}_{i_1}\leq x\} - \mathds{1}\{\hat{u}_{i_2}\leq x\}|) d\lambda(x) 
        &\leq \E( \min(\lambda(A), |\hat{u}_{i_1} - \hat{u}_{i_2} | )\\
        &\leq \E( \min(\lambda(A), |\hat{u}_{i_1}|)) 
        + \E( \min(\lambda(A), |\hat{u}_{i_2}|)).
    \end{align*}
    It follows that
    \begin{align*}
        &\dfrac{1}{2k}\dfrac{1}{k(k-1)} \sum_{l_1 \neq l_2} \dfrac{1}{|K_{l_1}| |K_{l_2}|} 
        \sum_{i_1 \in K_{l_1}} \sum_{i_2 \in K_{l_2}} \int_{A} \E(|\mathds{1}\{\hat{u}_{i_1}\leq x\} - \mathds{1}\{\hat{u}_{i_2}\leq x\}|) d\lambda(x) \\
        &\leq 
        \dfrac{1}{2k}\dfrac{1}{k(k-1)} \sum_{l_1 \neq l_2} \dfrac{1}{|K_{l_1}| |K_{l_2}|} 
        \sum_{i_1 \in K_{l_1}} \sum_{i_2 \in K_{l_2}} \E( \min(\lambda(A), |\hat{u}_{i_1}|)) 
        + \E( \min(\lambda(A), |\hat{u}_{i_2}|)) \\
        &= \dfrac{1}{k^2} \sum_{j=1}^k \dfrac{1}{|K_j|} \sum_{i \in K_j} 
        \E( \min(\lambda(A), |\hat{u}_{i}|))
        = \dfrac{1}{k^2} \sum_{j=1}^k \E( \min(\lambda(A), |y_{n+1} - \predl{K_j}|)),
    \end{align*}
    where the last equality follows from the fact that the distribution of $\hat{u}_{i}$ coincides with that of $y_{n+1} - \predl{K_j}$ for each $i \in K_j$.
    Using the triangle inequality, the last expression in the preceding display can be bounded from above by
    \begin{align*}
        &\dfrac{1}{k^2} \sum_{j=1}^k \E( \min(\lambda(A), |y_{n+1} - \pred|)) + \E( \min(\lambda(A), |\pred - \predl{K_j}|)) \\
        &= \dfrac{1}{k} \E( \min(\lambda(A), |y_{n+1} - \pred|)) + \dfrac{1}{k^2} \sum_{j=1}^k \E( \min(\lambda(A), |\pred - \predl{K_j}|)).
    \end{align*}
    To finish the proof, we once again apply \Cref{lem:appJa_expIntIndFunctions} to show that
    \begin{nalign}\label{eq:appJa_boundingLdMiddlePartEq1}
        &\int_A \frac{5}{k} \sum_{j=1}^k \E( |\mathds{1}\{y_{n+1} - \pred \leq x\} - \mathds{1}\{y_{n+1} - \predl{K_j} \leq x \}|) d\lambda(x) \\
        &\leq \frac{5}{k} \sum_{j=1}^k \E( \min(\lambda(A), |(y_{n+1} - \pred) - (y_{n+1} - \predl{K_j})|)) \\
        &= \frac{5}{k} \sum_{j=1}^k \E( \min(\lambda(A), |\pred - \predl{K_j}|)).
    \end{nalign}
\end{proof}

\begin{lemma}\label{lem:appJa_ultimateBoundJackknife}
    Let $\gam > 0$, $\epsilon > 0$, $\mu \in \R$, 
    and $\Fh$ and $\Ft$ as in \Cref{lem:appJa_modifyingSteinLeebsLemma}.
    We then have for every $L \geq 0$ that
    \begin{align*}
        \P( \ld{\gam}(\Fh, \Ft) \geq \epsilon)
        \leq&~ \dfrac{2\P(|y_{n+1} - \pred - \mu| \geq L)}{\epsilon}
        + \dfrac{(8L+12\gam)}{k\gam \epsilon^2}\\
        &+ \dfrac{4(5k+1)}{k^2\gam\epsilon^2} \sum_{j=1}^k \E(\min(2L+3\gam, |\pred - \predl{K_j}|)).
    \end{align*}
    Furthermore,
    \begin{align*}
        \P( \ld{\gam}(\Fh, \Ft) \geq \epsilon) 
        \leq&~ \dfrac{1}{k\gam \epsilon^2} \E( |y_{n+1} - \pred - \mu|)
        + \dfrac{5k+1}{k^2 \gam\epsilon^2} \sum_{j=1}^k \E(|\pred - \predl{K_j}|).
    \end{align*}
\end{lemma}
\begin{proof}
    We may assume $\mu = 0$ without loss of generality as otherwise we can define a new algorithm which is shifted by $\mu$.
    By \Cref{lem:appJa_boundingLdSquared} and Markov's inequality we have
    \begin{align*}
        \P(& \ld{\gam}(\Fh, \Ft) \geq \epsilon) 
        \leq \P\left(1 - \Ft(L) + \Ft(- L) \geq \frac{\epsilon}{2}\right) \\
        &+ \P \left( \left[ \frac{1}{\gam} \int_{[- L-\gam, L+2\gam]} |\Fh(x) - \Ft(x)|^2 d\lambda(x) \right]^{\frac{1}{2}} 
        \geq \frac{\epsilon}{2} \right)\\
        \leq&~ \dfrac{2 \E(1 - \Ft(L) + \Ft(- L))}{\epsilon} + \dfrac{4}{\gam \epsilon^2} \E \left(\int_{[- L-\gam, L+2\gam]} |\Fh(x) - \Ft(x)|^2 d\lambda(x) \right).
    \end{align*}
    We now bound the terms separately:
    \begin{align*}
        &\E(1 - \Ft(L) + \Ft(-L))  \\
        &= \E( \PC{y_{n+1} - \pred \leq -L}{\Tn}
        + \PC{y_{n+1} - \pred > L}{\Tn}) \\
        &\leq \P(|y_{n+1} - \pred| \geq L).
    \end{align*}
    For the second term we use \Cref{lem:appJa_boundingLdMiddlePart} with $A = [- L-\gam, L+2\gam]$ to get
    \begin{align*}
        &\dfrac{4}{\gam \epsilon^2} \E \left(\int_{[- L-\gam, L+2\gam]} |\Fh(x) - \Ft(x)|^2 d\lambda(x) \right) \\
        &\leq \dfrac{4}{k \gam \epsilon^2} \E(\min(2L+3\gam, |y_{n+1} - \pred|))
        + \dfrac{4(5k+1)}{k^2 \gam \epsilon^2} \sum_{j=1}^k \E( \min(2L+3\gam, |\pred - \predl{K_j}|)),
    \end{align*}
    which proves the first statement.
    
    For the second statement, we can use \Cref{lem:appJa_boundingLdSquared} again to get
    \begin{align*}
        \P( \ld{\gam}(\Fh, \Ft) \geq \epsilon)
        &\leq \dfrac{\E(\ld{\gam}(\Fh, \Ft)^2)}{\epsilon^2}
        \leq \dfrac{1}{\gam \epsilon^2} \int_\R \E( |\Fh(x) - \Ft(x)|^2 ) d\lambda(x).
    \end{align*}
    Applying \Cref{lem:appJa_boundingLdMiddlePart} with $A = \R$  yields the upper bound
    \begin{align*}
        \dfrac{1}{\gam \epsilon^2} \int_\R \E( |\Fh(x) - \Ft(x)|^2 ) d\lambda(x)
        \leq&~ \dfrac{1}{k \gam \epsilon^2} \E(|y_{n+1} - \pred|) \\
        &+ \dfrac{(5k+1)}{k^2 \gam \epsilon^2} \sum_{j=1}^k \E( |\pred - \predl{K_j}|),
    \end{align*}
    which was exactly what we had to prove.
\end{proof}

\begin{lemma}\label{lem:appEqi_CondCoverageEquivalence}
    Let $\epsilon > 0$, $\gam \geq 0$, $\kappa \geq 0$ and $\alpha_1, \alpha_2 \in \R$. We then have almost surely that
    \begin{nalign}\label{eq:appEqi_CondCoverageEquivalenceEq1}
        \PC{y_{n+1} \in \PCvp{\alpha_1 - \epsilon, \alpha_2 + \epsilon}{\gam + \kappa}}{\Tn}
        \geq&~ \PC{y_{n+1} \in \PCv{\alpha_1, \alpha_2}{\kappa}}{\Tn} \\
        &- \dfrac{1}{k \epsilon} \sum_{j=1}^k \PC{ | \pred - \predl{K_j} | > \gam}{\Tn}
    \end{nalign}  
    and the same holds true if we switch the roles of CV and CV+.
\end{lemma}
\begin{proof}
    We start noticing that the results for the \ldname in \Cref{sec:ld} 
    hold true for \emph{deterministic} distribution functions. 
    Therefore, we have to condition on both the training data $\Tn$ and $x_{n+1}$ in order to apply these results.
    For fixed training data $\Tn$ and $x_{n+1}$ we then know from \Cref{lem:ld_ldDefinitionViaQuantiles} with $\G$ and $\Gp$ replacing $F$ and $G$, respectively, that either 
    $\PCvp{\alpha_1, \alpha_2}{\kappa} \subseteq \PCv{\alpha_1 - \epsilon, \alpha_2 + \epsilon}{\gam + \kappa}$ or 
    $\ld{\gam}(\G, \Gp) > \epsilon$ holds true.
    Using this and the tower property of the conditional expectation, it is easy to see that almost surely
    \begin{align*}
        &\PC{y_{n+1} \in \PCv{\alpha_1 - \epsilon, \alpha_2 + \epsilon}{\kappa + \gam}}{\Tn} \\
        &= \EC{ \PC{y_{n+1} \in \PCv{\alpha_1 - \epsilon, \alpha_2 + \epsilon}{\kappa + \gam}}{\Tn, x_{n+1}} }{\Tn}\\
        &\geq \EC{ \PC{y_{n+1} \in \PCvp{\alpha_1, \alpha_2}{\kappa}}{\Tn, x_{n+1}} - \mathds{1}\{ \ld{\gam}(\G, \Gp) > \epsilon \} }{\Tn}\\
        &\geq \PC{y_{n+1} \in \PCvp{\alpha_1, \alpha_2}{\kappa}}{\Tn} - \PC{ \ld{\gam}(\G, \Gp) > \epsilon }{\Tn}.
    \end{align*}
    Applying \Cref{cor:ld_BoundForEcdfs} with $p_i = \frac{1}{k |K_{j(i)}|}$ we get the following bound almost surely
    \begin{align*}
        \ld{\gam}(\G, \Gp) 
        &\leq \sum_{i=1}^n \frac{1}{k |K_{j(i)}|} \mathds{1}\{ | v_i - v_i^{+} | > \gam\}
        = \sum_{i=1}^n \frac{1}{k |K_{j(i)}|} \mathds{1}\{ | \pred - \predl{K_{j(i)}} | > \gam\} \\
        &= \frac{1}{k} \sum_{j=1}^k
        \mathds{1}\{ | \pred - \predl{K_{j}} | > \gam\},
    \end{align*}
    which yields
    \begin{align*}
        \PC{ \ld{\gam}(\G, \Gp) > \epsilon }{\Tn}
        &\leq \EC{ \dfrac{\ld{\gam}(\G, \Gp)}{ \epsilon }}{\Tn} \\
        &\leq \frac{1}{k \epsilon} \sum_{j=1}^k \PC{ | \pred - \predl{K_j} | > \gam }{\Tn}
    \end{align*}
    almost surely.
    The second claim can be derived by exchanging the roles of CV and CV+.
\end{proof}

\subsection{Proof of the main results}
\begin{proof}[Proof of \Cref{prop:appRem_finSamCV}]
    Fix training data $\Tn$ and let 
    $\Fh$ and $\Ft$ be defined as in \Cref{lem:appJa_modifyingSteinLeebsLemma}.
    In particular, for fixed training data, $\Ft$, $\Fh$ and their corresponding quantiles are deterministic.
    Moreover, also $\ld{\gam}(\Ft, \Fh)$ is fixed.
    Take any $0 \leq \alpha_1 \leq \alpha_2 \leq 1$ and define the (infeasible) oracle prediction interval as 
    \begin{align*}
        \widetilde{PI}_{\alpha_1 + \ld{\gam}(\Ft, \Fh), \alpha_2 - \ld{\gam}(\Ft, \Fh)} = 
    \pred + [\Q{\alpha_1 + \ld{\gam}(\Ft, \Fh)}{\Ft}, \Q{\alpha_2 - \ld{\gam}(\Ft, \Fh)}{\Ft}].
    \end{align*}
    We then have
    \begin{align*}
        &\PC{y_{n+1} \in \widetilde{PI}_{\alpha_1 + \ld{\gam}(\Ft, \Fh), \alpha_2 - \ld{\gam}(\Ft, \Fh)}}{\Tn} \\
        &= \PC{y_{n+1} - \pred \in [\Q{\alpha_1 + \ld{\gam}(\Ft, \Fh)}{\Ft}, \Q{\alpha_2 - \ld{\gam}(\Ft, \Fh)}{\Ft}]}{\Tn}.
    \end{align*}
    
    If $\alpha_2 - \epsilon < \alpha_1 + \epsilon$ we have $\alpha_2 - \alpha_1 - 2 \epsilon < 0$ and we can therefore restrict the infimum in the statement of \Cref{prop:appRem_finSamCV} to the set $\{0 \leq \alpha_1 \leq \alpha_2 \leq 1: \alpha_2 \geq \alpha_1 + 2 \epsilon\}$.    
    If $\alpha_2 - \epsilon \geq \alpha_1 + \epsilon$ holds true, we conclude $0 \leq \alpha_1 \leq \alpha_1 + \epsilon \leq \alpha_2 - \epsilon \leq \alpha_2 \leq 1$. Thus, we also have $0 \leq \alpha_1 + \ld{\gam}(\Ft, \Fh) \leq \alpha_2 - \ld{\gam}(\Ft, \Fh) \leq 1$ whenever $\ld{\gam}(\Ft, \Fh) \leq \epsilon$.
    From the fact that $\Ft(\Q{\alpha}{\Ft}) \geq \alpha \geq \Ft(\Q{\alpha}{\Ft}-)$ for every $\alpha \in [0,1]$, we conclude
    \begin{align*}
        \PC{y_{n+1} \in \widetilde{PI}_{\alpha_1 + \ld{\gam}(\Ft, \Fh), \alpha_2 - \ld{\gam}(\Ft, \Fh)}}{\Tn} \geq \alpha_2 - \alpha_1 - 2 \ld{\gam}(\Fh, \Ft)
    \end{align*}
    whenever $\ld{\gam}(\Ft, \Fh) \leq \epsilon$.
    Thus, \Cref{cor:ld_generalPiInequality} yields
    \begin{align*}
        \PC{y_{n+1} \in \PCv{\alpha_1, \alpha_2}{\gam}}{\Tn}
        &= \PC{y_{n+1} - \pred \in [\Q{\alpha_1}{\Fh} - \gam, \Q{\alpha_2}{\Fh} + \gam]}{\Tn} \\
        &\geq \PC{y_{n+1} \in \widetilde{PI}_{\alpha_1 + \ld{\gam}(\Ft, \Fh), \alpha_2 - \ld{\gam}(\Ft, \Fh)}}{\Tn} \\
        &\geq \alpha_2 - \alpha_1 - 2 \ld{\gam}(\Fh, \Ft) - \mathds{1}\{\ld{\gam}(\Fh, \Ft) \geq \epsilon \} \text{ a.s.}
    \end{align*}
    where the last inequality follows from the preceding display and the fact that the lower bound is trivial if $\ld{\gam}(\Fh, \Ft) \geq \epsilon$.

    In particular, we have shown that $\PC{y_{n+1} \in \PCv{\alpha_1, \alpha_2}{\gam}}{\Tn} - (\alpha_2 - \alpha_1)$ can be bounded from below by 
    $-2\ld{\gam}(\Fh, \Ft) - \mathds{1}\{\ld{\gam}(\Fh, \Ft) \geq \epsilon \}$. We stress the fact that the lower bound is independent of the choice of $\alpha_1$ and $\alpha_2$ and therefore holds uniformly. 
    Before proceeding, we emphasize that the set $\{\Q{\alpha}{\Fh}: \alpha \in \R\}$ of all possible quantiles of $\Fh$ coincides with $\{\hat{u}_i: 1 \leq i \leq n\} \cup \{-\infty, \infty\}$. Thus, we can express $\inf_{0 \leq \alpha_1 \leq \alpha_2 \leq 1} \PC{y_{n+1} \in \PCv{\alpha_1, \alpha_2}{\gam}}{\Tn} - (\alpha_2 - \alpha_1)$ as the minimum of a \emph{finite} set.\footnote{
        It is not hard to show that the infimum equals $\min\left(0, \min_{0 \leq i_1 < i_2 \leq n} \Ft(\Q{\frac{i_2}{n}}{\Fh}+\delta) - \Ft((\Q{\frac{i_1+1}{n}}{\Fh}-\delta)-) - \frac{i_2 - i_1}{n}\right)$.
    }
    Hence, the infimum is indeed measurable, and from the inequality derived above, we conclude that even 
    \begin{align*}
        \inf_{0 \leq \alpha_1 \leq \alpha_2 \leq 1} \PC{y_{n+1} \in \PCv{\alpha_1, \alpha_2}{\gam}}{\Tn} - (\alpha_2 - \alpha_1) \geq -2 \ld{\gam}(\Fh, \Ft) - \mathds{1}\{\ld{\gam}(\Fh, \Ft) \geq \epsilon \}
    \end{align*}
    holds almost surely.
    This shows
    \begin{align*}
        \P \left[ \inf_{0 \leq \alpha_1 \leq \alpha_2 \leq 1} \PC{y_{n+1} \in \PCv{\alpha_1, \alpha_2}{\gam}}{\Tn} - (\alpha_2 - \alpha_1) \leq -2 \epsilon \right]
        \leq \P( \ld{\gam}(\Fh, \Ft) \geq \epsilon),
    \end{align*}
    where the right-hand side of the preceding display can be bounded from above by \Cref{lem:appJa_ultimateBoundJackknife}.

    For the shrunken prediction intervals, we start with a positive sequence $(\gam_m)_{m \in \N}$ converging to $\gam$ with
    $\gam_m < \gam$ for all $m \in \N$. We then have
    \begin{align*}
        \PC{y_{n+1} \in \PCv{\alpha_1, \alpha_2}{-\gam}}{\Tn} 
        &\leq \PC{y_{n+1} - \pred \in (\Q{\alpha_1}{\Fh} + \gam_m, \Q{\alpha_2}{\Fh} - \gam_m)}{\Tn} \\
        &\leq \PC{y_{n+1} - \pred \in (\Q{\alpha_1 - \ld{\gam_m}(\Fh, \Ft)}{\Ft}, \Q{\alpha_2 + \ld{\gam_m}(\Fh, \Ft)}{\Ft})}{\Tn},
    \end{align*}
    where we would like to point out that the last two prediction intervals are open sets and the second inequality in the preceding display follows from \Cref{prop:ld_quantileInequality}.
    Recalling \cref{eq:ld_alphaQuantilesExceedAlpha}, we conclude
    \begin{align*}
        \PC{y_{n+1} - \pred \in (\Q{\alpha_1 - \ld{\gam_m}(\Fh, \Ft)}{\Ft}, \Q{\alpha_2 + \ld{\gam_m}(\Fh, \Ft)}{\Ft})}{\Tn}
        \leq \alpha_2 - \alpha_1 + 2 \ld{\gam_m}(\Fh, \Ft)
    \end{align*}
    almost surely (even if $\alpha_2 + \ld{\gam_m}(\Fh, \Ft) > 1$ or $\alpha_1 -\ld{\gam_m}(\Fh, \Ft) < 0$).
    Thus, with the same arguments as before we conclude
    \begin{align*}
        \P \left( \sup_{0 \leq \alpha_1 \leq \alpha_2 \leq 1} \PC{y_{n+1} \in \PCv{\alpha_1, \alpha_2}{-\gam}}{\Tn}  
        - (\alpha_2 - \alpha_1) \geq 2 \epsilon \right)
        \leq \P( \ld{\gam_m}(\Fh, \Ft) \geq \epsilon ).
    \end{align*}
    Now, the right-hand side of the preceding display can be bounded from above by \Cref{lem:appJa_ultimateBoundJackknife} as in the first part of the proof. In fact, this only shows the statement of \Cref{prop:appRem_finSamCV} if we replace $\gam$ in the bounds by $\gam_m$. However, the bounds of \Cref{prop:appRem_finSamCV} are continuous in $\gam_m$ and therefore the statement for $\gam$ can be derived by taking the limit $m \to \infty$.
\end{proof}

\begin{proof}[Proof of \Cref{prop:appRem_finSamEqui}]
    We define
    \begin{align*}
        c(\alpha_1, \alpha_2, \Tn) &= \PC{y_{n+1} \in \PCv{\alpha_1 - \epsilon, \alpha_2 + \epsilon}{\kappa + \gam}}{\Tn},\\
        c^{+}(\alpha_1, \alpha_2, \Tn) &= \PC{y_{n+1} \in \PCvp{\alpha_1, \alpha_2}{\kappa}}{\Tn}\\
        \eta(\Tn, \gam) &= \frac{1}{k} \sum_{j=1}^k \PC{ |\pred - \predl{K_j}| > \gam}{\Tn}.
    \end{align*}
    Then \Cref{lem:appEqi_CondCoverageEquivalence} yields
    \begin{align*}
        c(\alpha_1, \alpha_2, \Tn) - c^{+}(\alpha_1, \alpha_2, \Tn) \geq - \dfrac{\eta(\Tn, \gam)}{\epsilon} \text{ a.s.}
    \end{align*}
    which holds uniformly for all $\alpha_1, \alpha_2 \in \R$. 
    Since the quantiles of $\G$ and $\Gp$ only take a finite number of values (including $\{-\infty, \infty\}$) the expression $\inf_{\alpha_1, \alpha_2 \in \R} c(\alpha_1, \alpha_2, \Tn) - c^{+}(\alpha_1, \alpha_2, \Tn)$ can indeed be rewritten as the minimum over a finite set.    
    Thus, the inequality from above also holds true almost surely for the infimum.
    \begin{align*}
        &\P\left( \inf_{\alpha_1, \alpha_2 \in \R} c(\alpha_1, \alpha_2, \Tn) - c^{+}(\alpha_1, \alpha_2, \Tn) \leq -\epsilon \right) \\
        &\leq \P( \eta(\Tn, \gam) \geq \epsilon^2 )\leq \dfrac{ \E( \eta(\Tn, \gam)) }{\epsilon^2}
        = \dfrac{1}{k\epsilon^2} \sum_{j=1}^k \P( |\pred - \predl{K_j}| > \gam).
    \end{align*}
    The second claim can be derived by exchanging the roles of CV and CV+.
\end{proof}

\begin{proof}[Proof of \Cref{thm:jack_asymptoticValidPi} and \Cref{thm:appExt_kFoldCvAlsoValid} \ref{it:appExt_CvValid}]
    We start with the general case:
    Fix $\gam > 0$, and $\epsilon$ and $\zeta$ fulfilling $0 < 2 \zeta < \epsilon$.
    Define $\eta_{n,\zeta} = \frac{1}{k_n} \sum_{j=1}^{k_n} \P( |\pred - \predl{K_j^{(n)}}| > \zeta)$. 
    Now, let $(L_n)_{n \in \N}$ be any sequence fulfilling $\limN L_n = \infty$, $L_n = \smallO{(\eta_{n,\zeta} + 1/k_n)^{-1}}$. 
    Applying the first statement in \Cref{prop:appRem_finSamCV} with $L = L_n$, $k = k_n$, $\mu = 0$ we see that both
    $\frac{8L_n + 12 \gam}{k_n \gam \zeta^2}$ and $\P(|y_{n+1} - \pred|>L_n)$ vanish asymptotically, where we used the fact that $y_{n+1} - \pred$ is bounded in probability by assumption.
    Furthermore, the asymptotic stability implies 
    \begin{align*}
        &\limsup_{n \to \infty} \frac{1}{k_n} \sum_{j=1}^{k_n} \E( \min(2L_n + 3\gam, |\pred - \predl{K_j}|)) \\
        &\leq \limsup_{n \to \infty} (2L_n + 3\gam) \frac{1}{k_n} \sum_{j=1}^{k_n}  \P( |\pred - \predl{K_j}| > \zeta) + \zeta
        = \limsup_{n \to \infty} (2L_n + 3\gam) \eta_{n,\zeta} + \zeta
        = \zeta.
    \end{align*}
    Thus, we have
    \begin{align*}
        \limsup_{n \to \infty}\P \left[ \inf_{0 \leq \alpha_1 \leq \alpha_2 \leq 1} \PC{y_{n+1} \in \PCv{\alpha_1, \alpha_2}{\gam}}{\Tn} - (\alpha_2 - \alpha_1) \leq  -2\zeta \right] \leq \zeta.
    \end{align*}
    Since $\zeta$ can be made arbitrarily small, the first statement of the general case is proven. The statement concerning the shrunken prediction intervals can be proven analogously since the upper bound in \Cref{prop:appRem_finSamCV} does not change if we replace $\gam$ by $-\gam$ and the infimum by the supremum.

    For the continuous case, we start with the observation that for every $\gam > 0$ we have
    \begin{align*}
        &\left| \PC{y_{n+1} \in \PCv{\alpha_1, \alpha_2}{\gam}}{\Tn}
        - \PC{y_{n+1} \in \PCv{\alpha_1, \alpha_2}{-\gam}}{\Tn} \right| \\
        &\leq \PC{y_{n+1} \in \pred + \Q{\alpha_1}{\Fh} + [-\gam, \gam)}{\Tn}
        + \PC{y_{n+1} \in \pred + \Q{\alpha_2}{\Fh} + (-\gam, \gam]}{\Tn}.
    \end{align*}
    Furthermore, we have for $i \in \{1,2\}$ almost surely that
    \begin{nalign}\label{eq:appExtProofs_InflationBound}
        \PC{y_{n+1} \in \pred + \Q{\alpha_i}{\Fh} + [-\gam, \gam]}{\Tn}
        &= \E\left( \PC{y_{n+1} \in \pred + \Q{\alpha_i}{\Fh} + [-\gam, \gam]}{\Tn, x_{n+1}} \right)\\
        &\leq \E\left( \min(1, 2\gam \supnorm{f_{y_{n+1}\|x_{n+1}}}) \right).
    \end{nalign}
    To sum it up, we conclude 
    \begin{align*}
        &\left| \PC{y_{n+1} \in \PCv{\alpha_1, \alpha_2}{\gam}}{\Tn}
        - \PC{y_{n+1} \in \PCv{\alpha_1, \alpha_2}{-\gam}}{\Tn} \right|\\
        &\leq 2 \E\left( \min(1, 2\gam \supnorm{f_{y_{n+1}\|x_{n+1}}}) \right) \text{ a.s.}
    \end{align*}
    Let $\epsilon > 0$ be arbitrary. Then, Assumption \CC implies the existence of a (fixed) $\gam_{\epsilon} > 0$,
    such that $\sup_{n \in \N} 2\E(\min(1, 2\gam_{\epsilon} \supnorm{f_{y_{n+1}\|x_{n+1}}})) \leq \epsilon$.
    To see this, recall that by Assumption \CC there is a $K_\epsilon > 0$ such that $\sup_{n \in \N} \P(\supnorm{f_{y_{n+1}\|x_{n+1}}} > K_\epsilon) \leq \epsilon/4$. Defining $\gam_\epsilon$ as $\epsilon/(8K_\epsilon)$, we have
    \begin{align}\label{eq:appEqi_AssCCimpliesEq1}
        \sup_{n \in \N} 2\E(\min(1, 2\gam_{\epsilon} \supnorm{f_{y_{n+1}\|x_{n+1}}}))
        \leq 2\sup_{n \in \N}\P\left( 2\gam_{\epsilon} \supnorm{f_{y_{n+1}\|x_{n+1}}} > \frac{\epsilon}{4}\right) + \frac{\epsilon}{2}
        \leq \epsilon.
    \end{align}
    We would like to emphasize that $\epsilon$ and $\gam_\epsilon$ are held fixed, while the distribution of $y_{n+1}$ conditional on $x_{n+1}$ may depend on $n$.
    For this, we have almost surely
    \begin{align*}
        &\PC{y_{n+1} \in \PCv{\alpha_1, \alpha_2}{\gam_{\epsilon}}}{\Tn} - \epsilon
        \leq
        \PC{y_{n+1} \in \PCv{\alpha_1, \alpha_2}{-\gam_{\epsilon}}}{\Tn}
        \leq \PC{y_{n+1} \in \PCv{\alpha_1, \alpha_2}{0}}{\Tn}\\
        &\leq \PC{y_{n+1} \in \PCv{\alpha_1, \alpha_2}{\gam_{\epsilon}}}{\Tn}
        \leq \PC{y_{n+1} \in \PCv{\alpha_1, \alpha_2}{-\gam_{\epsilon}}}{\Tn}
        + \epsilon.
    \end{align*}
    Applying the results for the general case we conclude
    \begin{align*}
        \P& \left[ \sup_{0 \leq \alpha_1 \leq \alpha_2 \leq 1} \left| \PC{y_{n+1} \in \PCv{\alpha_1, \alpha_2}{0}}{\Tn} - (\alpha_2 - \alpha_1) \right| \geq 2 \epsilon \right]\\
        \leq&~ \P \left[ \inf_{0 \leq \alpha_1 \leq \alpha_2 \leq 1} \PC{y_{n+1} \in \PCv{\alpha_1, \alpha_2}{\gam_{\epsilon}}}{\Tn} - (\alpha_2 - \alpha_1) \leq -\epsilon \right]\\
        &+ \P \left[ \sup_{0 \leq \alpha_1 \leq \alpha_2 \leq 1}\PC{y_{n+1} \in \PCv{\alpha_1, \alpha_2}{-\gam_{\epsilon}}}{\Tn} -  (\alpha_2 - \alpha_1) \geq \epsilon \right]
        \underset{n \to \infty}{\longrightarrow} 0.
    \end{align*}
    Since $\epsilon > 0$ was arbitrary, the result follows.    
\end{proof}

\begin{proof}[Proof of \Cref{thm:eqi_PacBoundAsymptotic} and 
    \Cref{thm:appExt_kFoldCvAlsoValid} \ref{it:appExt_asymptoticEquivalence}:]
    We start showing that \ref{it:eqi_PacBoundAsymptoticGenJackp} implies \ref{it:eqi_PacBoundAsymptoticGenJack}.
    To see this, we start with the observation that
    \begin{align*}
        &\inf_{0 \leq \alpha_1 \leq \alpha_2 \leq 1} \PC{y_{n+1} \in \PCv{\alpha_1, \alpha_2}{\gam}}{\Tn} - (\alpha_2 - \alpha_1) \\
        &\geq \inf_{0 \leq \alpha_1 \leq \alpha_2 \leq 1} \PC{y_{n+1} \in \PCv{\alpha_1, \alpha_2}{\gam}}{\Tn} - \PC{y_{n+1} \in \PCvp{\alpha_1+\epsilon/4, \alpha_2-\epsilon/4}{\gam/2}}{\Tn} \\
        &\indent+ \inf_{0 \leq \alpha_1 \leq \alpha_2 \leq 1} \PC{y_{n+1} \in \PCvp{\alpha_1+\epsilon/4, \alpha_2-\epsilon/4}{\gam/2}}{\Tn} - (\alpha_2 - \alpha_1 - \frac{\epsilon}{2}) - \frac{\epsilon}{2}.
    \end{align*}
    Thus, we have
    \begin{align*}
        \P&~ \left[ \inf_{0 \leq \alpha_1 \leq \alpha_2 \leq 1} \PC{y_{n+1} \in \PCv{\alpha_1, \alpha_2}{\gam}}{\Tn} - (\alpha_2 - \alpha_1) \leq - \epsilon \right] \\
        \leq&~ \P \left[ \inf_{0 \leq \alpha_1 \leq \alpha_2 \leq 1} \PC{y_{n+1} \in \PCv{\alpha_1, \alpha_2}{\gam}}{\Tn} - \PC{y_{n+1} \in \PCvp{\alpha_1+\epsilon/4, \alpha_2-\epsilon/4}{\gam/2}}{\Tn} \leq - \epsilon/4 \right] \\
        &+ \P \left[ \inf_{0 \leq \alpha_1 \leq \alpha_2 \leq 1} \PC{y_{n+1} \in \PCvp{\alpha_1+\epsilon/4, \alpha_2-\epsilon/4}{\gam/2}}{\Tn} - (\alpha_2 - \alpha_1 - \frac{\epsilon}{2}) \leq - \epsilon/4 \right].
    \end{align*}
    We treat the two expressions in the upper bound separately:
    Starting with the second expression, we conclude that for  $\alpha_2 - \alpha_1 < \epsilon/2$ the term $\PC{y_{n+1} \in \PCvp{\alpha_1+\epsilon/4, \alpha_2-\epsilon/4}{\gam/2}}{\Tn} - (\alpha_2 - \alpha_1 - \frac{\epsilon}{2})$ is always nonnegative. Defining $\tilde{\alpha}_1 = \alpha_1 + \epsilon/4, \tilde{\alpha}_2 = \alpha_2 - \epsilon/4$ we conclude that $\alpha_2 - \alpha_1 \geq \epsilon/2$ yields $\tilde{\alpha}_1 \leq \tilde{\alpha}_2$ and we therefore have
    \begin{align*}
        &\P \left[ \inf_{0 \leq \alpha_1 \leq \alpha_2 \leq 1} \PC{y_{n+1} \in \PCvp{\alpha_1+\epsilon/4, \alpha_2-\epsilon/4}{\gam/2}}{\Tn} - (\alpha_2 - \alpha_1 - \frac{\epsilon}{2}) \leq - \epsilon/4 \right] \\
        &= \P \left[ \inf_{\epsilon/4 \leq \tilde{\alpha}_1 \leq \tilde{\alpha}_2 \leq 1-\epsilon/4} 
        \PC{y_{n+1} \in \PCvp{\tilde{\alpha}_1, \tilde{\alpha}_2}{\gam/2}}{\Tn} 
        - (\tilde{\alpha}_2 - \tilde{\alpha}_1) \leq - \epsilon/4 \right] \\
        &\leq \P \left[ \inf_{0 \leq \tilde{\alpha}_1 \leq \tilde{\alpha}_2 \leq 1} 
        \PC{y_{n+1} \in \PCvp{\tilde{\alpha}_1, \tilde{\alpha}_2}{\gam/2}}{\Tn} 
        - (\tilde{\alpha}_2 - \tilde{\alpha}_1) \leq - \epsilon/4 \right].
    \end{align*}
    By \ref{it:eqi_PacBoundAsymptoticGenJackp}, the last expression in the preceding display vanishes as $n \to \infty$.

    For the first expression we apply \Cref{prop:appRem_finSamEqui} with $\epsilon$ replaced by $\epsilon/4$ and both $\kappa$ and $\gam$ replaced by $\gam/2$ to get
    \begin{align*}
        &\P \left[ \inf_{0 \leq \alpha_1 \leq \alpha_2 \leq 1} \PC{y_{n+1} \in \PCv{\alpha_1, \alpha_2}{\gam}}{\Tn} - \PC{y_{n+1} \in \PCvp{\alpha_1+\epsilon/4, \alpha_2 - \epsilon/4}{\gam/2}}{\Tn} \leq - \epsilon/4 \right] \\
        &\leq \P \left[ \inf_{\alpha_1, \alpha_2 \in \R} \PC{y_{n+1} \in \PCv{\alpha_1, \alpha_2}{\gam}}{\Tn} - \PC{y_{n+1} \in \PCvp{\alpha_1+\epsilon/4, \alpha_2-\epsilon/4}{\gam/2}}{\Tn} \leq - \epsilon/4 \right] \\
        &= \P \left[ \inf_{\alpha_1, \alpha_2 \in \R} \PC{y_{n+1} \in \PCv{\alpha_1-\epsilon/4, \alpha_2+\epsilon/4}{\gam}}{\Tn} - \PC{y_{n+1} \in \PCvp{\alpha_1, \alpha_2}{\gam/2}}{\Tn} \leq - \epsilon/4 \right] \\
        &\leq \dfrac{16}{k\epsilon^2} \sum_{j=1}^k \P( |\pred - \predl{K_j}| > \gam/2).
    \end{align*}
    Since the predictor is stable, this bound vanishes for $n \to \infty$.

    To prove that \ref{it:eqi_PacBoundAsymptoticGenJack} implies \ref{it:eqi_PacBoundAsymptoticGenJackp}
    we can use exactly the same arguments and exchange the roles of CV and CV+.
    
    Next, we prove that \ref{it:eqi_PacBoundAsymptoticConJackp} implies \ref{it:eqi_PacBoundAsymptoticConJack} in the continuous case.
    For any $\alpha_1, \alpha_2 \in \R$ and any $\gam \in \R$
    \begin{align*}
        &\left| \PC{y_{n+1} \in \PCv{\alpha_1, \alpha_2}{\gam}}{\Tn} - \PC{y_{n+1} \in \PCv{\alpha_1, \alpha_2}{0}}{\Tn} \right| \\
        &\leq \EC{ \left| \PC{y_{n+1} \in \PCv{\alpha_1, \alpha_2}{\gam}}{\Tn, x_{n+1}} - \PC{y_{n+1} \in \PCv{\alpha_1, \alpha_2}{0}}{\Tn, x_{n+1}} \right| }{\Tn} \\
        &\leq \EC{ \min(1, 2 |\gam| \supnorm{f_{y_{n+1}\|x_{n+1}}})}{\Tn} = \E( \min(1, 2 |\gam| \supnorm{f_{y_{n+1}\|x_{n+1}}}))
    \end{align*}
    holds almost surely, and the same holds true for the CV+.
    The last expression in the preceding display can be bounded from above by
    \begin{align*}
        2 |\gam| L + \P( \supnorm{f_{y_{n+1}\|x_{n+1}}} > L) \text{ for all } L > 0.
    \end{align*}
    Since $\supnorm{f_{y_{n+1}\|x_{n+1}}}$ is bounded in probability by Assumption \textbf{CC}, with a similar argument as in \cref{eq:appEqi_AssCCimpliesEq1} we can find for every $\epsilon/4 > 0$ a $\gam_\epsilon > 0$ (independent of $n \in \N$) such that $\E( \min(1, 2 |\gam_\epsilon| \supnorm{f_{y_{n+1}\|x_{n+1}}}))$ is bounded from above by $\epsilon/4$ uniformly over $n \in \N$.
    To sum it up, we have shown that 
    \begin{align*}
        \left| \PC{y_{n+1} \in \PCv{\alpha_1, \alpha_2}{\gam_\epsilon}}{\Tn} - \PC{y_{n+1} \in \PCv{\alpha_1, \alpha_2}{0}}{\Tn} \right|
        \leq \epsilon/4
    \end{align*}
    almost surely for every $n \in \N$.

    Since for bounded random variables convergence in $L_1$ coincides with convergence in probability, we only have to show that 
    \begin{align*}
        \limN \P \left( \sup_{0 \leq \alpha_1 \leq \alpha_2 \leq 1} |\PC{y_{n+1} \in \PCv{\alpha_1, \alpha_2}{0}}{\Tn} - (\alpha_2 - \alpha_1)| \geq \epsilon \right) = 0 
        \text{ for all } \epsilon > 0
    \end{align*}
    assuming that the same statement holds true for the CV+ instead of the CV.
    We start with the bound
    \begin{align*}
        \P& \left( \inf_{0 \leq \alpha_1 \leq \alpha_2 \leq 1} \PC{y_{n+1} \in \PCv{\alpha_1, \alpha_2}{0}}{\Tn} - (\alpha_2 - \alpha_1) \leq -\epsilon \right) \\
        \leq&~ \P \left( \inf_{0 \leq \alpha_1 \leq \alpha_2 \leq 1} \PC{y_{n+1} \in \PCv{\alpha_1, \alpha_2}{\gam_\epsilon}}{\Tn} - (\alpha_2 - \alpha_1) 
            \leq -\frac{3}{4} \epsilon \right) \\
        \leq&~ \P \left( \inf_{0 \leq \alpha_1 \leq \alpha_2 \leq 1} \PC{y_{n+1} \in \PCvp{\alpha_1+\epsilon/8, \alpha_2-\epsilon/8}{0}}{\Tn} - (\alpha_2 - \alpha_1 -\epsilon/4) \leq - \epsilon/4 \right) \\
        &+ \P\left( \inf_{0 \leq \alpha_1 \leq \alpha_2 \leq 1} \PC{y_{n+1} \in \PCv{\alpha_1, \alpha_2}{\gam_\epsilon}}{\Tn}
        - \PC{y_{n+1} \in \PCvp{\alpha_1+\epsilon/8, \alpha_2-\epsilon/8}{0}}{\Tn} \leq - \epsilon/4\right).
    \end{align*}
    We treat the two expressions in the upper bound separately.
    By \Cref{prop:appRem_finSamEqui} and the asymptotic stability we conclude that 
    \begin{align*}
        \limN \P\left( \inf_{0 \leq \alpha_1 \leq \alpha_2 \leq 1} \PC{y_{n+1} \in \PCv{\alpha_1, \alpha_2}{\gam_\epsilon}}{\Tn}
        - \PC{y_{n+1} \in \PCvp{\alpha_1+\epsilon/8, \alpha_2-\epsilon/8}{0}}{\Tn} \leq - \epsilon/4\right) = 0.
    \end{align*}
    For the other expression we start with a case distinction: If $(\alpha_2 - \alpha_1) < \epsilon/4$ the term $\PC{y_{n+1} \in \PCvp{\alpha_1+\epsilon/8, \alpha_2-\epsilon/8}{0}}{\Tn} - (\alpha_2 - \alpha_1 -\epsilon/4)$ is nonnegative. In the case $(\alpha_2 - \alpha_1) \geq \epsilon/4$ we have $\alpha_1 + \epsilon/8 \leq \alpha_2 - \epsilon/8$. 
    Setting $\tilde{\alpha}_1 = \alpha_1 + \epsilon/8$ and $\tilde{\alpha}_2 = \alpha_2 - \epsilon/8$, we conclude
    \begin{align*}
        &\P \left( \inf_{0 \leq \alpha_1 \leq \alpha_2 \leq 1} \PC{y_{n+1} \in \PCvp{\alpha_1+\epsilon/8, \alpha_2-\epsilon/8}{0}}{\Tn} - (\alpha_2 - \alpha_1 -\epsilon/4) \leq - \epsilon/4 \right)\\
        &= \P \left( \inf_{\epsilon/8 \leq \tilde{\alpha}_1 \leq \tilde{\alpha}_2 \leq 1-\epsilon/8} \PC{y_{n+1} \in \PCvp{\tilde{\alpha}_1, \tilde{\alpha}_2}{0}}{\Tn} - (\tilde{\alpha}_2 - \tilde{\alpha}_1) \leq - \epsilon/4 \right)\\
        &\leq \P \left( \inf_{0 \leq \alpha_1 \leq \alpha_2 \leq 1} \PC{y_{n+1} \in \PCvp{\alpha_1, \alpha_2}{0}}{\Tn} - (\alpha_2 - \alpha_1) \leq - \epsilon/4 \right)\\
        &\leq \P \left( \sup_{0 \leq \alpha_1 \leq \alpha_2 \leq 1} |\PC{y_{n+1} \in \PCvp{\alpha_1, \alpha_2}{0}}{\Tn} - (\alpha_2 - \alpha_1)| \geq \epsilon/4 \right).
    \end{align*}
    However, the last expression converges to $0$ since CV+ is assumed to be uniformly valid.
    To sum it up, we have shown that 
    \begin{align*}
        \limN \P \left( \inf_{0 \leq \alpha_1 \leq \alpha_2 \leq 1} \PC{y_{n+1} \in \PCv{\alpha_1, \alpha_2}{0}}{\Tn} - (\alpha_2 - \alpha_1) \leq -\epsilon \right) = 0
    \end{align*}
    for all $\epsilon > 0$. In fact, this is enough to prove that
    \begin{align*}
        \limN \P \left( \sup_{0 \leq \alpha_1 \leq \alpha_2 \leq 1} |\PC{y_{n+1} \in \PCv{\alpha_1, \alpha_2}{0}}{\Tn} - (\alpha_2 - \alpha_1)| \geq \epsilon \right) = 0.
    \end{align*}
    To see this, recall that for all $0 \leq \alpha_1 \leq \alpha_2 \leq 1$ we have $\Q{\alpha_1}{\G} \leq \Q{\alpha_2}{\G}$ and therefore
    \begin{align*}
        \R = (-\infty, \Q{\alpha_1}{\G}) \dot{\cup} \PCv{\alpha_1, \alpha_2}{0} \dot{\cup} (\Q{\alpha_2}{\G}, \infty)
    \end{align*}
    where $A \dot{\cup} B$ denotes the union of two \emph{disjoint} sets $A$ and $B$.
    By Assumption \textbf{CC} the distribution of $y_{n+1}$ given $x_{n+1}$ is continuous and therefore any singleton has point mass $0$. Since the three sets are disjoint, we conclude
    \begin{align*}
        \PC{y_{n+1} \in \PCv{\alpha_1, \alpha_2}{0}}{\Tn} \leq 1 - \PC{y_{n+1} \in \PCv{0, \alpha_1}{0}}{\Tn} - \PC{y_{n+1} \in \PCv{\alpha_2,1}{0}}{\Tn}.
    \end{align*}
    Thus, we have
    \begin{align*}
        &\sup_{0 \leq \alpha_1 \leq \alpha_2 \leq 1} \PC{y_{n+1} \in \PCv{\alpha_1, \alpha_2}{0}}{\Tn} - (\alpha_2 - \alpha_1) \\
        &\leq 1 - \inf_{0 \leq \alpha_1 \leq \alpha_2 \leq 1} \left[\PC{y_{n+1} \in \PCv{0, \alpha_1}{0}}{\Tn} - \alpha_1 + \PC{y_{n+1} \in \PCv{\alpha_2,1}{0}}{\Tn} + \alpha_2\right] \\
        &\leq -\inf_{0 \leq \alpha_1\leq 1} \left[\PC{y_{n+1} \in \PCv{0, \alpha_1}{0}}{\Tn} - \alpha_1\right] - \inf_{0 \leq \alpha_2 \leq 1} \left[\PC{y_{n+1} \in \PCv{\alpha_2,1}{0}}{\Tn} - (1 - \alpha_2)\right]\\
        &\leq - 2\inf_{0 \leq \alpha_1 \leq \alpha_2 \leq 1} \left[\PC{y_{n+1} \in \PCv{\alpha_1, \alpha_2}{0}}{\Tn} - (\alpha_2 - \alpha_1)\right].
    \end{align*}
    Hence, we conclude
    \begin{align*}
        &\limN \P \left( \sup_{0 \leq \alpha_1 \leq \alpha_2 \leq 1} \PC{y_{n+1} \in \PCv{\alpha_1, \alpha_2}{0}}{\Tn} - (\alpha_2 - \alpha_1) \geq \epsilon \right)\\
        &\leq \limN \P \left( -2\inf_{0 \leq \alpha_1 \leq \alpha_2 \leq 1} [\PC{y_{n+1} \in \PCv{\alpha_1, \alpha_2}{0}}{\Tn} - (\alpha_2 - \alpha_1)] \geq \epsilon \right)\\
        &= \limN \P \left( \inf_{0 \leq \alpha_1 \leq \alpha_2 \leq 1} \PC{y_{n+1} \in \PCv{\alpha_1, \alpha_2}{0}}{\Tn} - (\alpha_2 - \alpha_1) \leq -\epsilon/2 \right)
        = 0,
    \end{align*}
    which proves \ref{it:eqi_PacBoundAsymptoticConJack}.
    To prove that \ref{it:eqi_PacBoundAsymptoticConJack} implies \ref{it:eqi_PacBoundAsymptoticConJackp}, we can use exactly the same arguments and exchange
    the roles of CV and CV+.
\end{proof}

\begin{proof}[Proof of \Cref{prop:appJa_ShrunkenEqInflated}]
    We first show that \ref{it:appJa_Inflated} implies \ref{it:appJa_Shrunken}.
    For this, we fix $\gam > 0$ and $\epsilon > 0$, $0 \leq \alpha_1 \leq \alpha_2 \leq 1$ and distinguish two cases.
    If $\Q{\alpha_2}{\G} < \Q{\alpha_1}{\G} + 4 \gam$, then $\PCv{\alpha_1, \alpha_2}{-2\gam} = \emptyset$ and therefore the prediction interval has $0$ coverage probability.
    Otherwise, we have $\PCv{\alpha_1, \alpha_2}{-2\gam} \dot{\cup} \PCv{0, \alpha_1}{\gam} \dot{\cup} \PCv{\alpha_2,1}{\gam} \subseteq \R$,
    which implies
    \begin{align*}
        \PC{y_{n+1} \in \PCv{\alpha_1, \alpha_2}{-2\gam}}{\Tn} \leq 1 - \PC{y_{n+1} \in \PCv{0, \alpha_1}{\gam}}{\Tn} - \PC{y_{n+1} \in \PCv{\alpha_2,1}{\gam}}{\Tn}
    \end{align*}
    almost surely. In view of these two cases, we conclude
    \begin{align*}
        &\P \left[ \sup_{0 \leq \alpha_1 \leq \alpha_2 \leq 1} [\PC{y_{n+1} \in \PCv{\alpha_1, \alpha_2}{-2\gam}}{\Tn} - (\alpha_2 - \alpha_1)] > 2\epsilon \right] \\
        &\leq \P \left[ \sup_{0 \leq \alpha_1 \leq \alpha_2 \leq 1} [1 - \PC{y_{n+1} \in \PCv{0, \alpha_1}{\gam}}{\Tn} - \PC{y_{n+1} \in \PCv{\alpha_2, 1}{\gam}}{\Tn} - (\alpha_2 - \alpha_1)] > 2\epsilon \right].
    \end{align*}
    Now the last expression in the preceding display can be bounded from above by
    \begin{align*}
        \P&\left[ \sup_{0 \leq \alpha_2 \leq 1} [(1 - \alpha_2) - \PC{y_{n+1} \in \PCv{\alpha_2, 1}{\gam}}{\Tn}] > \epsilon\right] \\
        &+ \P\left[ \sup_{0 \leq \alpha_1 \leq 1} [\alpha_1 - \PC{y_{n+1} \in \PCv{0, \alpha_1}{\gam}}{\Tn}] > \epsilon\right] \\
        \leq&~ 2 \P \left[ \sup_{0 \leq \alpha_1 \leq \alpha_2 \leq 1} [(\alpha_2 - \alpha_1) - \PC{y_{n+1} \in \PCv{\alpha_1, \alpha_2}{\gam}}{\Tn}] > \epsilon \right].
    \end{align*}
    Note that we have
    \begin{align*}
        &\P \left[ \sup_{0 \leq \alpha_1 \leq \alpha_2 \leq 1} [(\alpha_2 - \alpha_1) - \PC{y_{n+1} \in \PCv{\alpha_1, \alpha_2}{\gam}}{\Tn}] > \epsilon \right]\\
        &= 1 - \P \left[ \inf_{0 \leq \alpha_1 \leq \alpha_2 \leq 1} [\PC{y_{n+1} \in \PCv{\alpha_1, \alpha_2}{\gam}}{\Tn} - (\alpha_2 - \alpha_1)] \geq -\epsilon \right].
    \end{align*}
    Now, \ref{it:appJa_Inflated} implies that the last line in the preceding display converges to $0$ for $n \to \infty$, which proves \ref{it:appJa_Shrunken}.

    Next, we show that \ref{it:appJa_Shrunken} implies \ref{it:appJa_Inflated}. 
    We start with $\gam > 0$ and $\epsilon > 0$ and note that for all $0 \leq \alpha_1 \leq \alpha_2 \leq 1$ we have 
    $\PCv{\alpha_1, \alpha_2}{2\gam} \supseteq \R \backslash (\PCv{0, \alpha_1}{-\gam} \dot{\cup} \PCv{\alpha_2,1}{-\gam})$,
    which implies
    \begin{align*}
        \PC{y_{n+1} \in \PCv{\alpha_1, \alpha_2}{2\gam}}{\Tn} \geq 1 - \PC{y_{n+1} \in \PCv{0, \alpha_1}{-\gam}}{\Tn} - \PC{y_{n+1} \in \PCv{\alpha_2,1}{-\gam}}{\Tn}
    \end{align*}
    almost surely. Thus, we have
    \begin{align*}
        \P& \left[ \inf_{0 \leq \alpha_1 \leq \alpha_2 \leq 1} [\PC{y_{n+1} \in \PCv{\alpha_1, \alpha_2}{2\gam}}{\Tn} - (\alpha_2 - \alpha_1)] < -2\epsilon \right] \\
        \leq& \P \bigg[ \inf_{0 \leq \alpha_1 \leq \alpha_2 \leq 1} [1 - \PC{y_{n+1} \in \PCv{0, \alpha_1}{-\gam}}{\Tn} \\
        &- \PC{y_{n+1} \in \PCv{\alpha_2, 1}{-\gam}}{\Tn} - (\alpha_2 - \alpha_1)] < -2\epsilon \bigg].
    \end{align*}
    With similar arguments as above, the last line in the preceding display can be bounded from above by
    \begin{align*}
        &2 \P \left[ \inf_{0 \leq \alpha_1 \leq \alpha_2 \leq 1} [(\alpha_2 - \alpha_1) - \PC{y_{n+1} \in \PCv{\alpha_1, \alpha_2}{-\gam}}{\Tn}] < -\epsilon \right]\\
        & = 2\left(1 - \P \left( \sup_{0 \leq \alpha_1 \leq \alpha_2 \leq 1} 
                \left[ \PC{ y_{n+1} \in \PCv{\alpha_1, \alpha_2}{-\gam} }{\Tn} - (\alpha_2 - \alpha_1) \right] 
                \leq \epsilon \right) \right),
    \end{align*}
    which proves \ref{it:appJa_Inflated} given \ref{it:appJa_Shrunken}.

    All our arguments can be applied mutatis mutandis if we replace CV with CV+.
\end{proof}

\begin{proof}[Proof of \Cref{cor:appExt_CvPlusValid} and \Cref{cor:eqi_asymptoticValidPiPlus}]
    The statement for the continuous case ($\delta = 0$) and the inflated prediction intervals ($\delta > 0$) are a direct consequence of \Cref{thm:appExt_kFoldCvAlsoValid}.
    For the shrunken prediction intervals, we additionally need the equivalence result of \Cref{prop:appJa_ShrunkenEqInflated}.
\end{proof}

\section{Proofs for \Cref{sec:necessity}}\label{sec:app_NecProofs}
\begin{proof}[Proof of \Cref{prop:nec_expectedInfinityLength}:]
    Since the prediction error is not bounded in probability, there exists a $\delta > 0$ such that for every $K > 0$ we have
    \begin{align*}
        \limsup_{n \in \N} \P(|y_{n+1} -\pred| > K) \geq 2 \delta.
    \end{align*}
    Since the method provides asymptotically marginal conservative prediction intervals, we also have for every $0 \leq \alpha_1 \leq \alpha_2 \leq 1$ fulfilling $\alpha_2 - \alpha_1 \geq 1 - \delta$:
    \begin{align*}
        1 - \delta 
        \leq& \liminf_{n \in \N} \P(y_{n+1} \in PI_{\alpha_1, \alpha_2}^{n,p}(x_{n+1}, \Tn)) \\
        \leq& \liminf_{n \in \N} \P(y_{n+1} \in PI_{\alpha_1, \alpha_2}^{n,p}(x_{n+1}, \Tn), 
        PI_{\alpha_1, \alpha_2}^{n,p}(x_{n+1}, \Tn) \subseteq \pred + [-K, K]) \\
        &+ \limsup_{n \in \N} \P(PI_{\alpha_1, \alpha_2}^{n,p}(x_{n+1}, \Tn) \nsubseteq \pred + [-K, K])\\
        \leq& \liminf_{n \in \N} \P(y_{n+1} \in \pred + [-K, K]) + 
        \limsup_{n \in \N} \P(PI_{\alpha_1, \alpha_2}^{n,p}(x_{n+1}, \Tn) \nsubseteq \pred + [-K, K]).
    \end{align*}
    Combining the last line of the preceding display with the fact that the prediction error is not bounded in probability, we 
    conclude
    \begin{align*}
        \limsup_{n \in \N} \P(PI_{\alpha_1, \alpha_2}^{n,p}(x_{n+1}, \Tn) \nsubseteq \pred + [-K, K]) \geq \delta
    \end{align*}
    for every $K > 0$.
    Note that $0$ is almost surely contained in the set $PI_{\alpha_1, \alpha_2}^{n,p}(x_{n+1}, \Tn) - \pred$. Since $PI_{\alpha_1, \alpha_2}^{n,p}(x_{n+1}, \Tn) - \pred$ is an interval (containing $0$) its length is no less than $K$ whenever $PI_{\alpha_1, \alpha_2}^{n,p}(x_{n+1}, \Tn) - \pred \nsubseteq [-K, K]$.
    Thus, we conclude
    \begin{align*}
        \limsup_{n \in \N} \E(\lambda(PI_{\alpha_1, \alpha_2}^{n,p}(x_{n+1}, \Tn)))
        \geq \limsup_{n \in \N} K \P(PI_{\alpha_1, \alpha_2}^{n,p}(x_{n+1}, \Tn) \nsubseteq \pred + [-K, K])
        \geq \delta K.
    \end{align*}
    Since $K > 0$ can be made arbitrarily large, the statement follows.
\end{proof}

\begin{proof}[Proof of \Cref{lem:nec_equivalentStability}]
    We start by introducing the following notation:
    \begin{align}
        \eta_n &:= |\pred - \predl{n}|\nonumber \\
        \delta_n &:= \E(\pred) - \E(\predl{n})\nonumber\\
        \xi_n &:= \EC{\pred}{\Tnl{n}, x_{n+1}} - \predl{n}.\label{eq:nec_equivalentStabilityEq1}
    \end{align}
    Thus, the condition of \Cref{lem:nec_equivalentStability} that assumes the drift added by a new observation to vanish in $L_2$ can be equivalently written as $\limN \E(\xi_n^2) = 0$. Since the $(2+\xi)$-th absolute moments are uniformly bounded over $n$, the asymptotic stability condition for symmetric prediction algorithms $\pred - \predl{n} \plim 0$ is equivalent to the condition $\limN \E(\eta_n^2) = 0$.
    In other words, we have to show that under the assumptions of \Cref{lem:nec_equivalentStability}, the following statement holds:
    \begin{align*}
        \limN \E(\eta_n^2) = 0 \Longleftrightarrow
        \left[ \limN \sn^2 - \snl^2 = 0 \text{ and } \limN \E(\xi_n^2) = 0 \right].
    \end{align*}

    Before proving the lemma, we need several preliminary considerations:
    By Jensen's inequality, we immediately conclude 
    $\E(\eta_n) \geq \E(|\xi_n|) \geq |\delta_n|$ and $\E(\eta_n^2) \geq \E(\xi_n^2) \geq \delta_n^2$.
    We then have
    \begin{align*}
        \sn^2 
        =&~ \E \left( \left( \pred - \E(\pred) \right)^2 \right) \\
        =&~ \E \left( \left( (\pred - \predl{n}) + (\predl{n} - \E(\pred)) \right)^2 \right)\\
        =&~ \E( \eta_n^2 ) + \E((\predl{n} - \E(\pred))^2) \\
        &+ 2\E( (\pred - \predl{n})(\predl{n} - \E(\pred))) .
    \end{align*}
    Since the mean-squared error can be written as the variance plus the squared bias, we conclude
    \begin{align*}
        \E((\predl{n} - \E(\pred))^2) 
        = \snl^2 + \delta_n^2.
    \end{align*}
    Furthermore, we have
    \begin{align*}
        &\E( (\pred - \predl{n})(\predl{n} - \E(\pred)))\\
        &= \E( \xi_n (\predl{n} - \E(\pred)))\\
        &= \E( \xi_n (\predl{n} - \E(\predl{n})) )
        - \E( \xi_n \delta_n) \\
        &= \E( \xi_n (\predl{n} - \E(\predl{n})) ) - \delta_n^2,
    \end{align*}
    where we used the tower property to conclude $\E(\xi_n) = \delta_n$.
    Putting the pieces together, we have shown that
    \begin{align}\label{eq:main_ineq}
        \sn^2 - \snl^2 = \E(\eta_n^2)
        - \delta_n^2 + 2 \E( \xi_n (\predl{n} - \E(\predl{n})) ).
    \end{align}
    By Cauchy-Schwarz inequality, we get the following bound
    \begin{align}\label{eq:cauchyschwarz}
        |\E( \xi_n (\predl{n} - \E(\predl{n})) )| \leq  \snl \sqrt{\E(\xi_n^2)}.
    \end{align}

    We are now able to start with the proof of \Cref{lem:nec_equivalentStability}: First, assume the prediction algorithm is stable and recall that $\E(\eta_n^2) \geq \E(\xi_n^2) \geq \delta_n^2$ holds. Thus, we conclude that also $\E(\xi_n^2)$ and $\delta_n^2$ converge to $0$ under the asymptotic stability. It remains to show that $\limN \sn^2 - \snl^2 = 0$.
    For this, we stress the fact that $\snl$ is bounded over $n$ and therefore by \cref{eq:cauchyschwarz} we have
    \begin{align*}
        \limN \E( \xi_n (\predl{n} - \E(\predl{n})) ) = 0.
    \end{align*}
    To finish the proof we use \cref{eq:main_ineq} to get
    \begin{align*}
        &\limN \sn^2 - \snl^2 \\
        &= \limN \E(\eta_n^2) - \delta_n^2 + 2 \E( \xi_n (\predl{n} - \E(\predl{n})) ) = 0.
    \end{align*}
    
    It remains to show the reverse direction. For this, assume that the drift added by a new observation vanishes and the variance stabilizes, that is, $\limN \E(\xi_n^2) = 0$ and $\limN \sn^2 - \snl^2 = 0$.
    Then, Jensen's inequality shows that also $\limN \delta_n^2 = 0$ holds true.
    Combining \cref{eq:main_ineq} with \cref{eq:cauchyschwarz} yields
    \begin{align*}
        \E(\eta_n^2) \leq \sn^2 - \snl^2 + \delta_n^2 + 2 \snl \sqrt{\E(\xi_n^2)}.
    \end{align*}
    Recalling the fact that $\snl$ is uniformly bounded over $n$, we conclude that $\E(\eta_n^2)$ vanishes asymptotically, which implies asymptotic stability for the symmetric predictor.
\end{proof}

\begin{proof}[Proof of \Cref{prop:nec_VanishingVariance}]
    Assume the variance stabilizes. Then, by \Cref{lem:nec_equivalentStability} the predictor is stable. Since $y_{n+1}$ and $\pred$ are bounded in $L_2$, the prediction error is bounded in probability. Thus, by \Cref{thm:jack_asymptoticValidPi}, the Jackknife method provides asymptotically uniformly conditionally conservative prediction intervals.

    For the reverse direction, we will assume the variance does not stabilize and proceed in the following two steps:
    Firstly, we will show that $\limsup_{n \to \infty} \E(\ld{\gam}(\Fh, \Ft)) > 0$ for some $\gam > 0$ and, secondly, we will use this to show that indeed there exists a $\epsilon > 0$ such that
    \begin{align*}
        \limsup_{n \to \infty} \P \left(\inf_{0 \leq \alpha_1 \leq \alpha_2 \leq 1} \PC{y_{n+1} \in \PJi{\alpha_1, \alpha_2}{\gam}}{\Tn} - (\alpha_2 - \alpha_1) \leq - \epsilon \right) > 0.
    \end{align*}
    We start with the first part:
    Let $\epsilon > 0$ such that 
    \begin{align}\label{eq:nec_VanishingVariance1}
        8\epsilon \leq \limsup_{n \to \infty} |\sn^2 - \snl^2|.
    \end{align}
    Since the $(2+\xi)$-th absolute moments of $\pred$ and $y_{n+1}$ are uniformly bounded over $n \in \N$ we conclude that $(y_{n+1} - \pred)^2$ is uniformly integrable. 
    Thus, there exists an $M>0$, such that 
    $$\limsup_{n \to \infty} \E((y_{n+1} - \pred)^2 \mathds{1}\{|y_{n+1} - \pred| \geq M\} \leq \epsilon.$$
    With the same argumentation, we can choose $M$ large enough that also 
    $$\limsup_{n \to \infty} \E((y_{n+1} - \predl{n})^2 \mathds{1}\{|y_{n+1} - \predl{n}| \geq M\}) \leq \epsilon.$$
    
    Let $f_M: \R \to \R$ denote the function $f_M(x) = \min(M^2, x^2)$. 
    It is easy to see that
    $$\E(|f_M(y_{n+1} - \pred) - (y_{n+1} - \pred)^2|) \leq \epsilon$$
    and the same holds true if we replace $\pred$ by $\predl{n}$.
    We obtain
    \begin{align*}
        &|\E((y_{n+1} - \pred)^2) - \E((y_{n+1} - \predl{n})^2)|\\
        &= |\E((y_{n+1} - \pred)^2) - \frac{1}{n} \sum_{i=1}^n \E(\hat{u}_i^2)|\\
        &\leq |\E( f_M(y_{n+1} - \pred)) - \frac{1}{n} \sum_{i=1}^n \E(f_M(\hat{u}_i))|
        + 2 \epsilon \\
        &\leq \E \left( \left| \EC{f_M(y_{n+1} - \pred)}{\Tn} - \frac{1}{n} \sum_{i=1}^n f_M(\hat{u}_i) \right| \right) + 2 \epsilon.
    \end{align*}
    We now split $f_M$ into two functions $f_M^+ - f_M^-$ defined as $f_M^+(x) = f_M(x)\mathds{1}\{x \geq 0\}$ and $f_M^-(x) = -f_M(x)\mathds{1}\{x < 0\}$. Then $f_M^+$ and $f_M^-$ are non-decreasing Lipschitz continuous functions with constant $2M$. Thus, 
    we can apply \cref{eq:dis_totalvariationLipschitz} of \Cref{lem:dis_totalVariation} for any $\gam > 0$ to get
    \begin{align*}
        |\EC{f_M(y_{n+1} - \pred)}{\Tn} - \frac{1}{n} \sum_{i=1}^n f_M(\hat{u}_i)| \leq 4M \gam + 2M^2 \ld{\gam}(\Fh, \Ft).
    \end{align*}
    Choosing $\gam = \epsilon/M > 0$ gives the upper bound $4\epsilon + 2M^2 \ld{\epsilon/M}(\Fh, \Ft)$.
    Putting the pieces together, we end up with
    \begin{align}\label{eq:nec_VanishingVariance2}
        \E(\ld{\epsilon/M}(\Fh, \Ft)) \geq \dfrac{|\E((y_{n+1} - \pred)^2) - \E((y_{n+1} - \predl{n})^2)| - 6 \epsilon}{2M^2}.
    \end{align}
    To find a lower bound for the right-hand side, we notice that 
    \begin{align*}
        &|\E((y_{n+1} - \pred)^2) - \E((y_{n+1} - \predl{n})^2)| \\
        &= |\E(\pred^2) - \E((\predl{n})^2) + 2 \E(y_{n+1}(\predl{n} - \pred))|.
    \end{align*}
    Since the one-step-update drift $\xi_n$ defined in \cref{eq:nec_equivalentStabilityEq1} converges to $0$ in $L_2$ and the second moment of $y_{n+1}$ is bounded over $n$ we conclude that 
    \begin{align*}
        \limN \E(y_{n+1}(\predl{n} - \pred)) = -\limN \E(y_{n+1} \xi_n) = 0
    \end{align*}
    as a consequence of the fact that $\predl{n}$ and $y_{n+1}$ are independent of $(y_n, x_n)$.
    Furthermore, we have
    \begin{align*}
        \limN \E(\pred) - \E(\predl{n})
        = \limN \E(\xi_n) = 0.
    \end{align*}
    Recalling that the second moments are uniformly bounded, this also yields
    \begin{align*}
        \limN \E(\pred)^2 - \E(\predl{n})^2 = \limN (\E(\pred) - \E(\predl{n}))(\E(\pred) + \E(\predl{n})) = 0.
    \end{align*}
    Recalling \cref{eq:nec_VanishingVariance1}, we have shown that
    \begin{align*}
        &\limN |\E((y_{n+1} - \pred)^2) - \E((y_{n+1} - \predl{n})^2)| \\
        &\geq \limN |\sn^2 - \snl^2| - \limN |\E(\pred)^2 - \E(\predl{n})^2| 
        \geq 8\epsilon.
    \end{align*}

    Combining this lower bound with \cref{eq:nec_VanishingVariance2} we get
    \begin{align*}
        \limsup_{n \to \infty} \E(\ld{\epsilon/M}(\Fh, \Ft)) \geq \dfrac{\epsilon}{M^2}.
    \end{align*}
    Since $\ld{\epsilon/M}(\Fh, \Ft))$ is bounded from above by $1$ we immediately get
    \begin{align*}
        \limsup_{n \to \infty}  \E(\ld{\epsilon/M}(\Fh, \Ft)) \leq
        \P\left(\ld{\epsilon/M}(\Fh, \Ft) > \frac{\epsilon}{2M^2}\right)
        + \frac{\epsilon}{2M^2}.
    \end{align*}
    In other words, we have found (another) $\tilde{\gam} = \epsilon/(2M) > 0$ and an $\tilde{\epsilon} = \epsilon/(4M^2) > 0$, such that
    $\limsup_{n \to \infty} \P\left(\ld{2\tilde{\gam}}(\Fh, \Ft) > 2\tilde{\epsilon}\right) \geq 2\tilde{\epsilon}$ holds true.
    
    By \Cref{lem:ld_ldDefinitionViaQuantiles} we can rewrite $\ld{2\tilde{\gam}}(\Fh, \Ft)$ as
    \begin{align*}
        \inf\{d \geq 0: \Q{\alpha-d}{\Ft} - 2\tilde{\gam} \leq \Q{\alpha}{\Fh} \leq \Q{\alpha + d}{\Ft} + 2\tilde{\gam} \text{ for all } \alpha \in \R \}.
    \end{align*}
    Thus, we have shown that
    \begin{align*}
        \limsup_{n \to \infty} \P( \exists \alpha \in \R: \Q{\alpha-2\tilde{\epsilon}}{\Ft} - 2\tilde{\gam} > \Q{\alpha}{\Fh} 
        \text{ or } \Q{\alpha + 2\tilde{\epsilon}}{\Ft} + 2\tilde{\gam} < \Q{\alpha}{\Fh}) \geq 2 \tilde{\epsilon},
    \end{align*}
    which yields
    \begin{align*}
        \limsup_{n \to \infty}  \max \left( \P( \exists \alpha \in \R: \Q{\alpha-2\tilde{\epsilon}}{\Ft} - 2\tilde{\gam} > \Q{\alpha}{\Fh}), 
        \P( \exists \alpha \in \R: \Q{\alpha + 2 \tilde{\epsilon}}{\Ft} + 2\tilde{\gam} < \Q{\alpha}{\Fh}) \right) \geq \tilde{\epsilon}.
    \end{align*}
    
    We start with the case $\limsup_{n \to \infty}  \P( \exists \alpha \in \R: \Q{\alpha-2\tilde{\epsilon}}{\Ft} - 2\tilde{\gam} > \Q{\alpha}{\Fh}) \geq \tilde{\epsilon}$. Now, fix training data $\Tn$ such that the inner inequality is fulfilled.
    Since $\Q{\alpha}{F} = -\infty$ whenever $\alpha \leq 0$ and $\Q{\alpha}{F} = +\infty$ whenever $\alpha > 1$ for every cdf $F$, we conclude that
    the inequality can only be fulfilled for $\alpha \in (2 \tilde{\epsilon}, 1]$. Let $\alpha^\ast \in (2 \tilde{\epsilon}, 1]$ denote such an $\alpha$ fulfilling
    $\Q{\alpha^\ast-2\tilde{\epsilon}}{\Ft} - 2\tilde{\gam} > \Q{\alpha^\ast}{\Fh}$.
    We then have for $\alpha_1 = 0$ and $\alpha_2 = \alpha^\ast$
    \begin{align*}
        &\PC{y_{n+1} \in \PJi{0, \alpha^\ast}{\tilde{\gam}}}{\Tn} = \PC{y_{n+1} - \pred \leq \Q{\alpha^\ast}{\Fh} + \tilde{\gam}}{\Tn} \\
        &\leq \PC{y_{n+1} - \pred \leq \Q{\alpha^\ast-2\tilde{\epsilon}}{\Ft} - \tilde{\gam}}{\Tn} 
        = \Ft(\Q{\alpha^\ast-2\tilde{\epsilon}}{\Ft} - \tilde{\gam}) \leq \alpha^\ast - 2 \tilde{\epsilon} = \alpha_2 - \alpha_1 - 2 \tilde{\epsilon}.
    \end{align*}
    Thus, in the first case, we get
    \begin{align*}
        \limsup_{n \to \infty} \P \left( 
        \inf_{0 \leq \alpha_1 \leq \alpha_2 \leq 1} \PC{y_{n+1} \in \PJi{\alpha_1, \alpha_2}{\tilde{\gam}}}{\Tn} - (\alpha_2 - \alpha_1) \leq - 2\tilde{\epsilon} \right) \geq \tilde{\epsilon}. 
    \end{align*}

    In the second case, we proceed similarly with the exception that we define $\alpha^\ast \in (0, 1-2\tilde{\epsilon}]$ such that
    $\Q{\alpha^\ast + 2\tilde{\epsilon}}{\Ft} + 2\tilde{\gam} < \Q{\alpha^\ast}{\Fh}$ holds true. Then, for $\alpha_1 = \alpha^\ast$ and $\alpha_2 = 1$ we have
    \begin{align*}
        \PC{y_{n+1} \notin \PJi{\alpha^\ast, 1}{\tilde{\gam}}}{\Tn} 
        &\geq \PC{y_{n+1} - \pred < \Q{\alpha^\ast}{\Fh} - \tilde{\gam}}{\Tn} \\
        &\geq \PC{y_{n+1} - \pred < \Q{\alpha^\ast + 2\tilde{\epsilon}}{\Ft} +\tilde{\gam}}{\Tn}
        \geq \Ft(\Q{\alpha^\ast + 2\tilde{\epsilon}}{\Ft})
        \geq \alpha^\ast + 2 \tilde{\epsilon},     
    \end{align*}
    which yields $\PC{y_{n+1} \in \PJi{\alpha^\ast, 1}{\tilde{\gam}}}{\Tn} \leq 1 - \alpha^\ast - 2\tilde{\epsilon} = \alpha_2 - \alpha_1 - 2 \tilde{\epsilon}$ and, subsequently,
    \begin{align*}
        \limsup_{n \to \infty} \P \left( 
        \inf_{0 \leq \alpha_1 \leq \alpha_2 \leq 1} \PC{y_{n+1} \in \PJi{\alpha_1, \alpha_2}{\tilde{\gam}}}{\Tn} - (\alpha_2 - \alpha_1) \leq - 2\tilde{\epsilon} \right) \geq \tilde{\epsilon}. 
    \end{align*}
\end{proof}

\begin{proof}[Proof of \Cref{prop:nec_OneStepUpdateDrift}]
    We start with the observation that by assumption
    \begin{align*}
        &\limsup_{n \to \infty} \E \left| \EC{y_{n+1} - \pred}{\Tn} - \avg{n}\hat{u}_i \right|\\
        &\geq \limsup_{n \to \infty} \left| \E(y_{n+1} - \pred) - \avg{n} \E(\hat{u}_i) \right|\\
        &= \limsup_{n \to \infty} | \E(y_{n+1} - \pred) - \E(y_{n+1} - \predl{n}) | > 0. 
    \end{align*}
    Let $M > 0$ and define functions $f_M^+, f_M^-$ and $f_M^0$ as
    \begin{align*}
        f_M^+(x) &= \max(0, \min(x, M)),\\
        f_M^-(x) &= \min(0, \max(x, -M)),\\
        f_M^0 &= x - f_M^+(x) - f_M^-(x).
    \end{align*}
    In particular, we have $x = f_M^+(x) + f_M^-(x) + f_M^0(x)$ for all $x \in \R$.
    Since the $(1+\xi)$-th absolute moments of $y_{n+1}$ and $\pred$ are uniformly bounded over $n \in \N$,
     $y_{n+1}$ and $\pred$ are uniformly integrable and hence so is the prediction error. The same holds true for $\pred - \predl{n}$.\footnote{This shows that indeed we could weaken the assumptions of this Lemma to $(1+\xi)$-th absolute moments instead of $(2+\xi)$-th absolute moments. We have decided to use the stronger condition to be consistent with the assumptions of the other two results in the subsection.}
    Thus, we conclude
    \begin{align*}
        \lim_{M \to \infty} \sup_{n \in \N} \E(|f_M^0(y_{n+1} - \pred)|) = 0, 
    \end{align*}
    which shows the existence of an $M > 0$, such that 
    \begin{align*}
        &\limsup_{n \to \infty} 
        \E \left| \EC{f_M^+(y_{n+1} - \pred)}{\Tn} - \avg{n}f_M^+(\hat{u}_i) \right| \\
        &+ \E \left| \EC{f_M^-(y_{n+1} - \pred)}{\Tn} - \avg{n}f_M^-(\hat{u}_i) \right|
        > 0.
    \end{align*}
    Let $C_M > 0$ denote the limes superior in the preceding display.
    Since $f_M^+$ and $-f_M^-$ are nondecreasing $1$-Lipschitz continuous functions, we can apply \cref{eq:dis_totalvariationLipschitz} of \Cref{lem:dis_totalVariation} and get for any arbitrary $\gam > 0$
    \begin{align*}
        \limsup_{n \to \infty} 2\gam + 2M \E(\ld{\gam}(\Fh, \Ft)) \geq C_M.
    \end{align*}
    By choosing $0 < 2\gam < C_M$ we conclude $\limsup_{n \to \infty} \E(\ld{\gam}(\Fh, \Ft)) > 0$.
    Now, we can proceed as in the proof of \Cref{prop:nec_VanishingVariance} to show that this indeed implies 
    \begin{align*}
        \limsup_{n \to \infty} \P \left( 
        \inf_{0 \leq \alpha_1 \leq \alpha_2 \leq 1} \PC{y_{n+1} \in \PJi{\alpha_1, \alpha_2}{\gam}}{\Tn} - (\alpha_2 - \alpha_1) \leq - 2\tilde{\epsilon} \right) \geq \tilde{\epsilon}. 
    \end{align*}
    for some $\gam > 0$ and $\tilde{\epsilon} > 0$.
\end{proof}

\section{Proofs for \Cref{sec:appRemarks}}\label{sec:app_RemProofs}
\begin{proof}[Proof of \Cref{prop:appExt_symmetrizedJackknife}]  
    Define 
    \begin{align*}
        \Ft^s = \PC{|y_{n+1} - \pred| \leq t}{\Tn}
    \end{align*}
    and $\Fh^s$ as the empirical distribution function of the absolute values $|\hat{u}_i|$ of the leave-one-out residuals.
    The main difference between the classical Jackknife and its symmetric version is now to use the empirical distribution function $\Fh^s$ 
    to estimate the distribution function $\Ft^s$ of the \emph{absolute} value $|y_{n+1} - \pred|$ of the prediction error. As we will see, the upper bounds for the \ldname between $\Fh$ and $\Ft$ and the \ldname between its symmetric counterparts coincide, and therefore \Cref{thm:jack_asymptoticValidPi} also holds for the symmetric Jackknife.

    To be more precise, \Cref{lem:appJa_modifyingSteinLeebsLemma} holds if we replace $\Ft$, $\Fh$ with $\Ft^s$ and $\Fh^s$, respectively and $\hat{u}_i$, $y_{n+1} - \pred$ and $y_{n+1} - \predl{i}$ with their absolute values. Furthermore, using $||y_{n+1} - \pred| - |y_{n+1} - \predl{i}|| \leq |\pred - \predl{i}|$ in \cref{eq:appJa_boundingLdMiddlePartEq1} we conclude that \Cref{lem:appJa_boundingLdMiddlePart} also holds if we replace $\Ft$, $\Fh$ with their symmetric variants and $y_{n+1} - \pred$ as well as $\pred - \predl{i}$ with their absolute values in the proof.
    \Cref{lem:appJa_ultimateBoundJackknife} also remains true with $\Fh^s$ and $\Ft^s$ instead of $\Fh$ and $\Ft$, if we replace $y_{n+1} - \pred - \mu$ and $\pred - \predl{i}$ with their absolute values in the proof.
    To put it in other words, \Cref{lem:appJa_ultimateBoundJackknife} yields
    the same upper bound for the \ldname between $\Fh^s$ and $\Ft^s$ as for the \ldname between $\Fh$ and $\Ft$ and therefore \Cref{prop:appRem_finSamCV} holds for the symmetric Jackknife.
    Furthermore, we conclude that under the assumptions of \Cref{thm:jack_asymptoticValidPi} also $\limN \E(\ld{\gam}(\Fh^s, \Ft^s)) = 0$ for all $\gam > 0$.
    As a direct consequence, \Cref{thm:jack_asymptoticValidPi} holds for the symmetric Jackknife.   

    With a similar argument, \Cref{thm:eqi_PacBoundAsymptotic} can be shown to hold true for the symmetric versions of the Jackknife and the Jackknife+:
    For every fixed $x_{n+1}$, $\Tn$ and $\gam \geq 0$ \Cref{cor:ld_BoundForEcdfs} allows us to bound the \ldname between $G_s$ and $G_s^+$ as follows
    \begin{align*}
        \ld{\gam}(G_s, G_s^+) \leq \avg{n} \mathds{1}\{|\pred - \predl{i}| > \gam\},
    \end{align*}
    which is exactly the same bound we derived for their asymmetric versions. 
    Therefore, \Cref{lem:appEqi_CondCoverageEquivalence} holds if we replace the original Jackknife and Jackknife+ with their symmetric variants. As a direct consequence, \Cref{prop:appRem_finSamEqui} and therefore also \Cref{thm:eqi_PacBoundAsymptotic} remain true for the symmetric variants.
    
    Combining the statements of \Cref{thm:jack_asymptoticValidPi} and \Cref{thm:eqi_PacBoundAsymptotic} for the symmetric Jackknife and Jackknife+ immediately leads to the symmetric analogue of \Cref{cor:eqi_asymptoticValidPiPlus}.
\end{proof}

\begin{proof}[Proof of \Cref{prop:appExt_equivalenceForFittedValues}:]
    Since we are comparing the fitted values to the Jackknife, we have $k_n = n$ and $K_j = \{j\}$ in this proof.
    
    Fix $x_{n+1}$, training data $\Tn$ and $\gam \geq 0$. Then the \ldname between $H$ and $G$ can be bounded from above using \Cref{cor:ld_BoundForEcdfs} as follows:
    \begin{align*}
        \ld{\gam}(G,H) \leq \avg{n} \mathds{1}\{|w_i - v_i| > \gam\} = \avg{n} \mathds{1}\{|\hat{y}_i - \hat{y}_i^{\backslash i}| > \gam\}.
    \end{align*}
    We would like to emphasize that the upper bound does not depend on $x_{n+1}$. Thus, for $\epsilon > 0$ we have
    \begin{align}\label{eq:appExt_equivalenceForFittedValuesEq1}
        \PC{\ld{\gam}(G,H) > \epsilon}{\Tn} \leq \mathds{1}\left\{\avg{n} \mathds{1}\{|\hat{y}_i - \hat{y}_i^{\backslash i}| > \gam\} > \epsilon \right\}.
    \end{align}
    almost surely. 

    We now proceed similarly to the proof of \Cref{lem:appEqi_CondCoverageEquivalence}:
    Applying \Cref{lem:ld_ldDefinitionViaQuantiles} yields
    \begin{align*}
        \PC{y_{n+1} \in PI^{fv}_{\alpha_1-\epsilon, \alpha_2+\epsilon}(\kappa + \gam)}{\Tn, x_{n+1}}
        \geq&~ \PC{y_{n+1} \in \PJi{\alpha_1, \alpha_2}{\kappa}}{\Tn, x_{n+1}} \\
        &- \mathds{1}\{\ld{\gam}(G,H) > \epsilon\}
    \end{align*}
    almost surely. Taking the expectation with respect to $x_{n+1}$ and using \cref{eq:appExt_equivalenceForFittedValuesEq1} yields
    \begin{nalign}\label{eq:appExt_equivalenceForFittedValuesEq2}
        \PC{y_{n+1} \in PI^{fv}_{\alpha_1-\epsilon, \alpha_2+\epsilon}(\kappa + \gam)}{\Tn}
        \geq&~ \PC{y_{n+1} \in \PJi{\alpha_1, \alpha_2}{\kappa}}{\Tn} \\
        &- \mathds{1}\left\{\avg{n} \mathds{1}\{|\hat{y}_i - \hat{y}_i^{\backslash i}| > \gam\} > \epsilon \right\},
    \end{nalign}
    which gives a result similar to \cref{eq:appEqi_CondCoverageEquivalenceEq1} in \Cref{lem:appEqi_CondCoverageEquivalence} with the difference, that the term $\frac{1}{\epsilon n} \sum_{i=1}^n \P(|\pred - \predl{i}| > \gam\|\Tn)$ is replaced by $\mathds{1}\{\frac{1}{n} \sum_{i=1}^n \mathds{1}\{|\hat{y}_i - \hat{y}_i^{\backslash i}| > \gam\} > \epsilon \}$.
    With the same arguments, we can also show that \cref{eq:appExt_equivalenceForFittedValuesEq2} holds if we exchange the roles of the prediction intervals. 
    We then have
    \begin{align*}
        \P\left( \mathds{1}\left\{\frac{1}{n} \sum_{i=1}^n \mathds{1}\{|\hat{y}_i - \hat{y}_i^{\backslash i}| > \gam\} > \epsilon \right\} \geq \epsilon \right)
        &\leq \P\left( \frac{1}{n} \sum_{i=1}^n \mathds{1}\{|\hat{y}_i - \hat{y}_i^{\backslash i}| > \gam \} > \epsilon \right)\\
        &\leq \frac{1}{n \epsilon} \sum_{i=1}^n \P(|\hat{y}_i - \hat{y}_i^{\backslash i}| > \gam),
    \end{align*}
    where we used Markov's inequality for the last inequality.
    In particular, we can proceed as in the proof of \Cref{prop:appRem_finSamEqui} to get a similar result with the difference that the term $\frac{1}{\epsilon^2 n} \P(|\pred - \predl{i}| > \gam)$ is replaced by $\frac{1}{n \epsilon} \sum_{i=1}^n \P(|\hat{y}_i - \hat{y}_i^{\backslash i}| > \gam)$. In particular, there is only an $\epsilon$ instead of $\epsilon^2$ in the denominator because the upper bound for $\ld{\gam}(G,H)$ is independent of $x_{n+1}$.
    For the asymptotic statement, we can now proceed as in the proof of \Cref{thm:eqi_PacBoundAsymptotic}.
\end{proof}

\begin{proof}[Proof of \Cref{prop:ext_lossFunctions}]
    The result directly follows from \Cref{lem:dis_totalVariation}.
    To see this, fix $\Tn$ and define $X$ as the prediction error's absolute value conditional on $\Tn$ and let $Y$ be uniformly distributed on the absolute values of the leave-one-out residuals in the sense that its distribution function is given by $F_Y(t) = \avg{n} \mathds{1}\{|\hat{u}_i| \leq t\}.$
    For any $M > 0$ we now apply \cref{eq:dis_totalvariationLower1} and \cref{eq:dis_totalvariationUpper1} of \Cref{lem:dis_totalVariation} to the function $f_M(x) = \max(-M, \min(M, \ell(x)))$ to get
    \begin{align*}
        &\P \left(\avg{n} f_M(|\hat{u}_i| - \epsilon) - \epsilon \leq \EC{f_M(|y_{n+1} - \pred|)}{\Tn} \leq \avg{n} f_M(|\hat{u}_i| + \epsilon) + \epsilon \right) \\
        &\geq \P( \ld{\epsilon}(F_X, F_Y) \leq \dfrac{\epsilon}{2M}).
    \end{align*}
    Note that for every $M \geq -\ell(\epsilon)$ we have
    $\ell(|y_{n+1} - \predl{n}| + \epsilon) \geq \ell(\epsilon) \geq -M$ and therefore 
    $\frac{1}{n} \sum_{i=1}^n f_M(|\hat{u}_i| + \epsilon) \leq \frac{1}{n} \sum_{i=1}^n \ell(|\hat{u}_i| + \epsilon)$ holds almost surely for such an $M > 0$.
    
    Furthermore, by the monotonicity of $\ell$ we have $\ell(-\epsilon) \leq \ell(|y_{n+1} - \predl{n}| - \epsilon) \leq \ell(|y_{n+1} - \predl{n}|)$, 
    and therefore the uniform integrability of $\ell(|y_{n+1} - \predl{n}| - \epsilon)$ follows from that of $\ell(|y_{n+1} - \predl{n}|)$.
    Thus, for any $\delta > 0$ this yields
    \begin{align*}
        &\lim_{M \to \infty} \limN \P\left(\left| \avg{n} \ell(|\hat{u}_i| - \epsilon) - f_M(|\hat{u}_i| - \epsilon)\right| \geq \delta\right) \\
        &\leq \lim_{M \to \infty} \limN \dfrac{1}{\delta} \E \left| \avg{n} \ell(|\hat{u}_i| - \epsilon) - f_M(|\hat{u}_i| - \epsilon)\right| \\
        &\leq \lim_{M \to \infty} \limN \dfrac{1}{n \delta} \sum_{i=1}^n \E \left| \ell(|\hat{u}_i| - \epsilon) - f_M(|\hat{u}_i| - \epsilon)\right| \\
        &= \lim_{M \to \infty} \limN \dfrac{1}{\delta} \E \left|\ell(|y_{n+1} - \predl{n}| - \epsilon) - f_M(|y_{n+1} - \predl{n}| - \epsilon)\right|
        = 0.
    \end{align*}
    Using a similar argument for $\EC{f_M(|y_{n+1} - \pred|)}{\Tn}$, we conclude that to prove \Cref{prop:ext_lossFunctions} 
    it is enough to show that
    \begin{align*}
        \lim_{M \to \infty} \limN \P \Bigg(\avg{n} f_M(|\hat{u}_i| - \epsilon) - \epsilon 
        &\leq \EC{f_M(|y_{n+1} - \pred|)}{\Tn} \\
        &\leq \avg{n} f_M(|\hat{u}_i| + \epsilon) + \epsilon \Bigg) = 1.
    \end{align*}
    
    Note that the distribution functions $F_X$ and $F_Y$ of $X$ and $Y$, respectively, coincide with $\Ft^s$ and $\Fh^s$ in the proof of \Cref{prop:appExt_symmetrizedJackknife}.
    Proceeding similar to the proof of \Cref{prop:appExt_symmetrizedJackknife}, we can show that $\limN \E( \ld{\epsilon}(F_X, F_Y)) = 0$ whenever $\limN \E(\ld{\epsilon}(\Fh, \Ft)) = 0$,
    where the latter can be derived from the asymptotic stability and the stochastic boundedness of the prediction error.    
\end{proof}

\begin{proof}[Proof of \Cref{cor:ext_MseConsistency}]
    Since the $(2+\xi)$-th absolute moments are uniformly bounded, $(y_{n+1} - \pred)^2$ and $(y_{n+1} - \predl{n})^2$ are uniformly integrable. Thus, we can apply \Cref{prop:ext_lossFunctions} to the function $\ell(x) = (\max(0,x))^2$ to get
    \begin{align*}
        \limN \P \left( \avg{n} \ell(|\hat{u}_i| - \epsilon) - \epsilon \leq \EC{ (y_{n+1} - \predl{n})^2}{\Tn} \leq \avg{n} (|\hat{u}_i| + \epsilon)^2 + \epsilon \right) = 1
    \end{align*}
    for all $\epsilon > 0$.
    Using standard arguments combining the inequality from above with the uniform integrability yields
    \begin{align*}
        \limsup_{n \to \infty} \E \left| \EC{ (y_{n+1} - \predl{n})^2}{\Tn} - \avg{n} (\hat{u}_i)^2 \right| \leq \epsilon + \epsilon^2 + 2\epsilon \limsup_{n \to \infty} \E(|y_{n+1} - \predl{n}|).
    \end{align*}
    Since $\epsilon > 0$ can be made arbitrarily small, the result follows.
\end{proof}

\begin{proof}[Proof of \Cref{cor:appExt_misclassificationError}]
    The result follows upon application of \Cref{prop:ext_lossFunctions}. For this, we start with the observation that the prediction error is bounded by $K$. 
    Furthermore, the predictor is stable by assumption.
    Defining the function $\ell(x) = \mathds{1}\{x \geq \frac{1}{2}\}$ we can now apply \Cref{prop:ext_lossFunctions} to get 
    \begin{align*}
        \P\Bigg[ \avg{n} \mathds{1}\left\{|\hat{u}_i| - \epsilon \geq \frac{1}{2}\right\} - \epsilon
        &\leq \PC{ |y_{n+1} - \pred| \geq \frac{1}{2}}{\Tn} \\
        &\leq \avg{n} \mathds{1}\left\{|\hat{u}_i| + \epsilon \geq \frac{1}{2}\right\} + \epsilon \Bigg] \underset{n \to \infty}{\longrightarrow} 1
    \end{align*}
    for every $\epsilon > 0$. Setting $\epsilon < \frac{1}{2}$ and noting that $\hat{u}_i$ and the prediction error are integer-valued we conclude that
    \begin{align*}
        \limN \P\left( \avg{n} \mathds{1}\{\hat{u}_i \neq 0\} - \epsilon \leq \PC{ y_{n+1} \neq \pred}{\Tn}
        \leq \avg{n} \mathds{1}\{\hat{u}_i \neq 0\} + \epsilon\right) = 1.
    \end{align*}
    Since $0 < \epsilon < \frac{1}{2}$ can be made arbitrarily small, we conclude $\avg{n} \mathds{1}\{\hat{u}_i \neq 0\} - \PC{ y_{n+1} \neq \pred}{\Tn} \plim 0$.
    By boundedness, this also implies the convergence to $0$ in $L_1$.
\end{proof}

\begin{proof}[Proof of \Cref{lem:appRem_counterExMStability}]
    For all $m,n \in \N$ let $\mathcal{P}_n^x$ be any $p_n$-dimensional distribution of $x_{n+m}^{(n)}$ such that its first component $x_{n+m,1}^{(n)}$ is Dirac distributed at the point $n$.
    Let $k \in \N$ and denote any training data set consisting of $k$ data points by $\mathcal{T}_k$.
    Define an algorithm $\mathcal{A}$ as follows:
    \begin{align*}
        \mathcal{A}_{p,k}(x, \mathcal{T}_k) = \mathcal{A}_{1,k}(x_1, \mathcal{T}_k)
        = \begin{cases}
            0 & \text{ if } x_1 \geq k \text{ and}\\
            L & \text{ else,}
        \end{cases}
    \end{align*}
    where $x_1$ denotes the first component of a vector $x \in \R^p$.
    Now, let $n \in \N$, $k \in \N$ and $x_{n+m}^{(n)} \sim \mathcal{P}_n^x$. We then have for every $k \in \N$
    \begin{align*}
        \mathcal{A}_{p_n,k}(x_{n+m}^{(n)}, \mathcal{T}_k) = \mathcal{A}_{1,k}(x_{n+m,1}^{(n)}, \mathcal{T}_k) = L\mathds{1}\{n < k\}.
    \end{align*}
    In particular, this yields that
    $\mathcal{A}_{p_n,n-1}(x_{n+m}^{(n)}, \Tnl{n}) = \mathcal{A}_{p_n,n}(x_{n+m}^{(n)}, \Tn) = 0$
    while for all $m \geq 2$ we have $\mathcal{A}_{p_n,n+m-1}(x_{n+m}^{(n)}, \mathcal{T}_{n+m-1}^{(n)}) = L$.
    Now, we conclude
    \begin{align*}
        \beta_{m,n-1}^{out}(\mathcal{A},\mathcal{P}_n) 
        :&= \E
        (|\mathcal{A}_{p_n,n+m-1}(x_{n+m}^{(n)}, \mathcal{T}_{n+m-1}^{(n)}) - \mathcal{A}_{p_n,n-1}(x_{n+m}^{(n)}, \Tnl{n})|) \\
        &= |\mathcal{A}_{1,n+m-1}(n, \mathcal{T}_{n+m-1}^{(n)}) - \mathcal{A}_{1,n-1}(n, \Tnl{n})|
        = L\mathds{1}\{m \geq 2\}.
    \end{align*}
    With the same argument one can show that $\widetilde{\beta}^{out}_{m,n-1}(\mathcal{A},\mathcal{P}_n)(L) = \mathds{1}\{m \geq 2\}$. 
    Since the predictor coincides with $\mathcal{A}_{p_n,n}(x_{n+1}^{(n)}, \Tn) = \mathcal{A}_{p_n,n}(x_{n+m}^{(n)}, \Tn) = 0$, it is bounded (almost surely/in probability). 
    In particular, the prediction error coincides with $y_{n+1}^{(n)}$ and is therefore bounded in probability whenever $y_{n+1}^{(n)}$ is.

\end{proof}

\begin{proof}[Proof of \Cref{lem:appRem_mStabilityGrowsWithSqM}:]
    We start proving inequalities in finite samples. For this, we suppress the dependence on $n$, that is, we write $m$ for $m_n$, $p$ for $p_n$ and $\beta$ for $\beta_n$.
    Recalling that $x_{n+m+1} \sim \mathcal{N}(0, I_{p})$ and the first absolute moment of a mean-zero Gaussian random variable with variance $\sigma^2$ is $\sigma \sqrt{2/\pi}$, we conclude
    \begin{align*}
        \E(|x_{n+m+1}'(\bhl{n} - \bhl{n+m})|) = \sqrt{\frac{2}{\pi}} \E( \norm{\bhl{n} - \bhl{n+m} }).
    \end{align*}
    Let $A = (x_1, \ldots, x_n)' \in \R^{n \times p}$, $B = (x_{n+1}, \ldots, x_{m+n})' \in \R^{m \times p}$, $Y = (y_1, \ldots, y_n)' \in \R^n$ and $Z = (y_{n+1}, \ldots, y_{n+m}) \in \R^m$.
    We then have
    \begin{nalign}\label{eq:appappRem_mStabilityGrowsWithSqMEq1}
        \bhl{n} - \bhl{n+m}
        &= \left( A'A + B'B + \lambda (n + m) I_p \right)^{-1} \left( A'Y + (B'B + \lambda m I_p)\bhl{n} - (A'Y + B'Z) \right)\\
        &= \left( A'A + B'B + \lambda (n + m) I_p \right)^{-1} \left( (B'B + \lambda m I_p)\bhl{n} - B'Z \right).
    \end{nalign}
    The idea of the proof is now the following: Based on \cref{eq:appappRem_mStabilityGrowsWithSqMEq1}, we will show that 
    $\liminf_{n \to \infty} \frac{n+m}{\sqrt{mp}} \E( \norm{\bhl{n} - \bhl{n+m} }) \geq C_1$ and $\limsup_{n \to \infty} \frac{n}{\sqrt{p}} \E( \norm{\bhl{n} - \bhl{n+1} }) \leq C_2$ for constants $C_1$ and $C_2$ depending only on $\rho, \lambda$ and $\sigma$, which proves \Cref{lem:appRem_mStabilityGrowsWithSqM}. Indeed, standard results of random matrix theory show that under the assumptions of \Cref{lem:appRem_mStabilityGrowsWithSqM} the eigenvalues of $A'A + B'B + \lambda (n + m) I_p$ are of the order $n+m$ almost surely (cf. \cite{bai1993limit}) and the main part of this proof will be concerned with finding bounds for $\norm{(B'B + \lambda m I_p)\bhl{n} - B'Z}$.

    For this, let 
    \begin{align*}
        g(A,B) = \mathds{1}\{a_- \leq \dfrac{\sigma_i(A)}{\sqrt{n}} \leq a_+, 
        b_- \leq \dfrac{\sigma_i(B)}{\sqrt{m}} \leq b_+ \text{ for all } 1 \leq i \leq p\},
    \end{align*}
    where $\sigma_1(D) \geq \ldots \geq \sigma_p(D)$ denote the singular values of a matrix $D \in \R^{l \times p}$ with $p \leq l$ and $a_-, a_+, b_-, b_+$ will be strictly positive values defined later.\footnote{
        Recall that we have $p \leq \rho m \leq \rho n$ for a $\rho < 1$ and therefore the matrices $A$ and $B$ have fewer columns than rows.
    }
    We then have
    \begin{align*}
        &\E( \norm{ \bhl{n} - \bhl{n+m} } )
        \geq \E( \norm{ \bhl{n} - \bhl{n+m} } g(A,B)) \\
        &\geq \E \left( \left[ \sigma_1 \left(A'A + B'B + \lambda(n+m) I_p \right) \right]^{-1} g(A,B)
        \norm{B'Z - (B'B + \lambda m I_p) \bhl{n}} \right)\\
        &\geq \left( \lambda(n+m) + n a_+^2 + m b_+^2 \right)^{-1} 
        \E( g(A,B) \norm{B'Z - (B'B + \lambda m I_p) \bhl{n}}).
    \end{align*}
    Using the abbreviation $c := (\lambda(n+m) + n a_+^2 + m b_+^2)^{-1}$ we have for every $\delta \geq 0$ the inequality
    \begin{align*}
        &c \E( g(A,B) \norm{B'Z - (B'B + \lambda m I_p) \bhl{n}}) \\
        &\geq c \E_{A,B,Y}\left( g(A,B) \delta \P \left( \norm{B'Z - (B'B + \lambda m I_p) \bhl{n}} \geq \delta \bigg\| A,B,Y \right) \right).
    \end{align*}
    Our goal is now to find a lower bound for $\PC{\norm{B'Z - (B'B + \lambda m I_p) \bhl{n}} \geq \delta}{A,B,Y}$. For this we rewrite $Z$ as $B \beta + U$, where $U = (u_{n+1}, \ldots, u_{n+m})'$ is an $m$-dimensional multivariate centered Gaussian random vector with covariance matrix $\sigma^2 I_m$ which is independent of $A,B$ and $Y$. Thus, we have almost surely
    \begin{align*}
        &\PC{\norm{B'Z - (B'B + \lambda m I_p) \bhl{n}} \geq \delta}{A,B,Y} \\
        &= 1 - \PC{\norm{B'B \beta + B'U  - (B'B + \lambda m I_p) \bhl{n}} < \delta}{A,B,Y} \\
        &\geq 1 - \supnorm{f_{B'U}} leb(B_\delta^p((B'B + \lambda m I_p) \bhl{n} - B'B \beta)),
    \end{align*}
    where $B_\delta^p(x)$ denotes the $p$-dimensional ball with radius $\delta$ and center $x$,
    $f_{B'U}$ the density of the random variable $B'U$ and $leb$ denotes here the $p$-dimensional Lebesgue measure.
    Now, the volume of the ball is independent of its center and we can therefore replace it by $leb(B_\delta^p(\mathbf{0}))$, where $\mathbf{0}$ denotes the $p$-dimensional zero vector.
    Since we condition on $B$, we conclude that
    $\supnorm{f_{B'U}}$ conditional on $B$ equals $det(B'B \sigma^2)^{-\frac{1}{2}} (2 \pi)^{-\frac{p}{2}}$.
    Furthermore, on the event $g(A,B) \neq 0$ we have
    $det(B'B \sigma^2)^{-\frac{1}{2}} \leq (\sigma b_- \sqrt{m})^{-p}$.
    In other words, we have
    \begin{align*}
        &c \E_{A,B}( g(A,B)) \delta \left(1 - leb(B_\delta^p(\mathbf{0})) (\sigma^2 b_-^2 m 2 \pi)^{-\frac{p}{2}} \right) \\
        &= c \delta \left(1 - leb(B_\delta^p(\mathbf{0})) (\sigma^2 b_-^2 m 2 \pi)^{-\frac{p}{2}} \right) \times \\
        &\P\left(a_- \leq \dfrac{\sigma_i(A)}{\sqrt{n}} \leq a_+, 
        b_- \leq \dfrac{\sigma_i(B)}{\sqrt{m}} \leq b_+ \text{ for all } 1 \leq i \leq p\right).
    \end{align*}
    We now claim that 
    $leb(B_\delta^p(\mathbf{0})) \leq \left( \dfrac{\delta^2 2 \pi e}{p} \right)^{\frac{p}{2}}$.
    To see this, we start with the fact that 
    $leb(B_\delta^p(\mathbf{0})) = \dfrac{\delta^p \pi^{\frac{p}{2}}}{\Gamma(\frac{p}{2} + 1)}$,
    where $\Gamma(x) = \int_0^\infty t^{x-1} e^{-t} dt$ denotes the gamma function.
    We now claim that $\Gamma(1 + p/2) \geq (\frac{p}{2e})^{\frac{p}{2}}$ holds true for all $p \in \N$. 
    Since the Gamma-function is nondecreasing on $[\frac{3}{2}, \infty)$, we have for $p \geq 3$
    \begin{align*}
        \Gamma\left(1 + \frac{p}{2}\right) \geq \int_{\frac{p}{2}}^{1 + \frac{p}{2}} \Gamma(x) dx.
    \end{align*}
    Because $[\frac{p}{2}, 1 + \frac{p}{2}]$ has length $1$, we can use Jensen's inequality to get that
    \begin{align*}
        \ln\left(\Gamma\left(1 + \frac{p}{2}\right)\right) \geq \ln \left(\int_{\frac{p}{2}}^{1 + \frac{p}{2}} \Gamma(x) dx \right)
        \geq \int_{\frac{p}{2}}^{1 + \frac{p}{2}} \ln(\Gamma(x)) dx.
    \end{align*}
    By Raabe's formula (cf. subsection $6.441$ in \cite{gradshteyn2014table}), we have
    \begin{align*}
        \int_{\frac{p}{2}}^{1 + \frac{p}{2}} \ln(\Gamma(x)) dx = \dfrac{1}{2}\ln(2 \pi) + \dfrac{p}{2}\ln\left(\frac{p}{2}\right) - \frac{p}{2},
    \end{align*}
    which gives
    \begin{align*}
        \Gamma\left(1 + \frac{p}{2}\right) \geq \sqrt{2\pi} \left(\dfrac{p}{2e}\right)^{\frac{p}{2}} \geq \left(\dfrac{p}{2e}\right)^{\frac{p}{2}}.
    \end{align*}
    The case $p=1$ can be seen from the fact that $\Gamma(\frac{3}{2}) = \frac{\sqrt{\pi}}{2} \geq \frac{1}{\sqrt{2e}}$, and the case $p=2$ follows from
    $\Gamma(2) = 1 \geq \frac{1}{e}$.    

    Thus, we have
    \begin{align*}
        leb(B_\delta^p(\mathbf{0})) \leq \left( \dfrac{\delta^2 2 \pi e}{p} \right)^{\frac{p}{2}}.
    \end{align*}
    Putting the pieces together, we have 
    \begin{align*}
        leb(B_\delta^p(\mathbf{0})) (\sigma^2 b_-^2 m 2 \pi)^{-\frac{p}{2}}
        \leq \left( \dfrac{\delta^2 e}{m p b_-^2 \sigma^2} \right)^{\frac{p}{2}}.
    \end{align*}
    Now, defining $\delta^* = \frac{ \sigma b_- \sqrt{mp}}{e}$ yields
    \begin{align*}
        leb(B_{\delta^*}^p(\mathbf{0})) (\sigma^2 b_-^2 m 2 \pi)^{-\frac{p}{2}}
        \leq e^{-\frac{p}{2}} \leq \dfrac{1}{\sqrt{e}}.
    \end{align*}
    Thus, we have
    \begin{align*}
        \E( \norm{ \bhl{n} - \bhl{n+m} } ) 
        &\geq 
        (\lambda(n+m) + n a_+^2 + m b_+^2)^{-1} \dfrac{ \sigma b_- \sqrt{mp}}{e} \left(1 - \dfrac{1}{\sqrt{e}} \right) \times \\
        &\P(a_- \leq \dfrac{\sigma_i(A)}{\sqrt{n}} \leq a_+, 
        b_- \leq \dfrac{\sigma_i(B)}{\sqrt{m}} \leq b_+ \text{ for all } 1 \leq i \leq p).
    \end{align*}
    Since $p \leq \rho m \leq \rho n$, standard results for random matrices (cf. \cite{bai1993limit}) show that the probability of all eigenvalues of $A'A/n$ and $B'B/m$ being contained in the interval $[(1-\sqrt{\rho})^2 - \epsilon, (1+\sqrt{\rho})^2 + \epsilon]$ with $\epsilon > 0$ converges to $1$ for $n \to \infty$. Since $\rho < 1$ we can find positive values $a_-, a_+, b_-, b_+$, such that
    \begin{align*}
        \limN \P(a_- \leq \dfrac{\sigma_i(A)}{\sqrt{n}} \leq a_+, 
        b_- \leq \dfrac{\sigma_i(B)}{\sqrt{m}} \leq b_+ \text{ for all } 1 \leq i \leq p_n) = 1.
    \end{align*}
    Thus, there exists a constant $C_1 > 0$ depending on $\rho, \lambda$ and $\sigma$, such that
    \begin{align*}
        \liminf_{n \to \infty} \dfrac{n+m_n}{\sqrt{m_n p_n}} \E( \norm{ \bhl{n} - \bhl{n+m_n} } ) \geq C_1.
    \end{align*}
    This proves a lower bound which crucially relies on the fact that $p_n \leq \rho m_n < m_n$.

    We now prove an upper bound for $m = 1$:
    Recalling \cref{eq:appappRem_mStabilityGrowsWithSqMEq1} we have
    \begin{align*}
        \bhl{n} - \bhl{n+1}
        = \left( A'A + x_{n+1}x_{n+1}' + \lambda (n + 1) I_p \right)^{-1} \left( (x_{n+1}x_{n+1}' + \lambda I_p)\bhl{n} - x_{n+1}y_{n+1} \right).
    \end{align*}
    Thus, we get the following upper bound:
    \begin{align*}
        \E( \norm{\bhl{n} - \bhl{n+1}} ) \leq 
        \dfrac{1}{\lambda (n+1)} \left( \E\left(\norm{x_{n+1}} |y_{n+1}| + \norm{ (x_{n+1}x_{n+1}' + \lambda I_p)\bhl{n} } \right) \right).
    \end{align*}
    However, we have
    \begin{align*}
        \E(\norm{x_{n+1}} |y_{n+1}| ) \leq \left( \E( \norm{x_{n+1}}^2 ) \E( y_{n+1}^2) \right)^{\frac{1}{2}}
        = \sqrt{p} \sqrt{\sigma^2 + \norm{\beta_n}^2}.
    \end{align*}
    We can bound the other term from above as follows:
    \begin{align*}
        \E(\norm{ (x_{n+1} x_{n+1}' + \lambda I_p)\bhl{n} })
        &\leq \E( \norm{x_{n+1}} |x_{n+1}'\bhl{n} | + \lambda \norm{\bhl{n}} ) \\
        &\leq \left( \E( \norm{x_{n+1}}^2 ) \E(|x_{n+1}'\bhl{n}|^2) \right)^{\frac{1}{2}} + \lambda \E(\norm{\bhl{n}}).
    \end{align*}
    Now, we have 
    $E( \norm{x_{n+1}}^2 ) = p$, $\E(|x_{n+1}'\bhl{n}|^2) = \E(\norm{\bhl{n}}^2)$.
    Let denote $s_i$ the $i$-th singular value of $A/\sqrt{n}$. 
    Then the singular values of $(A'A + \lambda n I_p)^{-1}A'$ are given by $\frac{s_i}{\sqrt{n}(s_i^2 + \lambda)}$ and can therefore be bounded from above by $\dfrac{1}{2 \sqrt{n \lambda}}$.
    Recalling the fact that 
    \begin{align*}
        \bhl{n} = (A'A + \lambda n I_p)^{-1}A'Y
    \end{align*}
    this yields
    \begin{align*}
        \E( \norm{\bhl{n}}^2 ) \leq \dfrac{1}{4n \lambda} \E(\norm{Y}^2)
        = \dfrac{\E(y_1)^2}{4 \lambda}.
    \end{align*}
    Putting the pieces together, we conclude that
    \begin{align*}
        \E(\norm{ (x_{n+1}x_{n+1}' + \lambda I_p)\bhl{n} }) \leq \tilde{C}_2 \sqrt{p_n},
    \end{align*}
    where $\tilde{C}_2 > 0$ depends on $M, \sigma$ and $\lambda$.
    In particular, this shows that
    \begin{align*}
        \E( \norm{\bhl{n} - \bhl{n+1}} ) \leq C_2 \dfrac{\sqrt{p_n}}{n},
    \end{align*}
    where $C_2 > 0$ is another constant depending only on $M, \sigma$ and $\lambda$. 
    Thus, there exists a constant $C > 0$ depending on $M, \sigma, \lambda$ and $\rho$, such that
    \begin{align*}
        \liminf_{n \to \infty} \dfrac{\beta_{m,n}^{out}(\mathcal{A}, \mathcal{P}_n)}{\sqrt{m_n} \beta_{1,n}^{out}(\mathcal{A}, \mathcal{P}_n)} \geq C,
    \end{align*}
    which was exactly what we had to prove.
\end{proof}

    \ifthenelse{\equal{\version}{arxive}}{
        \section{Proofs for \Cref{sec:ld}}\label{sec:appLd}
\begin{proof}[Proof of \Cref{lem:ld_basicProperties}:]
	To prove \Cref{it:ld_bpInfimumAttained} we choose a sequence $(\epsilon_n)_{n \in \N}$ converging from above to $\LD$. 
	We then have
	\begin{align*}
		F(t-\gam) - \epsilon_n \leq G(t) \leq F(t+\gam) + \epsilon_n \fatir 
	\end{align*}
	for each $n \in \N$.
	Taking the limits yields
	\begin{align*}
		F(t-\gam) - \LD & = \limN F(t-\gam) - \epsilon_n \leq G(t)                    \\
		                & \leq \limN F(t+\gam) + \epsilon_n = F(t+\gam) + \LD \fatir. 
	\end{align*}

	It is easy to see that $\LD$ can be equivalently written as
	\begin{align}\label{eq:appLd_SymmetricAltDefOfLd}
		\LD = \inf \{\epsilon \geq 0: G(t) \leq F(t+\gam)+\epsilon \text{ and } F(t) \leq G(t+\gam)+\epsilon \fatir \} 
	\end{align}
	or as
	\begin{align*}
		\LD = \inf \{\epsilon \geq 0: G(t - \gam) - \epsilon \leq F(t) \text{ and } F(t - \gam) - \epsilon \leq G(t) \fatir \}. 
	\end{align*}
	Since the definition in \cref{eq:appLd_SymmetricAltDefOfLd} is symmetric in $F$ and $G$, we conclude \Cref{it:ld_bpSymmetry}.

    For \Cref{it:ld_bpMonoCont}, the monotonicity of $F$ and \Cref{it:ld_bpInfimumAttained} yield for any $\epsilon \geq \gam \geq 0$:
	\begin{align*}
		F(t - \epsilon) - \LD
		  & \leq F(t - \gam) - \LD \leq G(t) \\
		  & \leq F(t + \gam) + \LD           
		\leq F(t + \epsilon) + \LD \fatir,
	\end{align*}
	which proves $\ld{\epsilon}(F,G) \leq \LD$. Thus, the function $\gam \mapsto \LD$ is nonincreasing.
	For the continuity from the right, we start with a sequence $(\gam_n)_{n \in \N}$ converging from above to $\gam \geq 0$.
	By the monotonicity and the fact that $\ld{\gam_n}(F,G)$ is bounded by $1$, we conclude that $\limN \ld{\gam_n}(F,G)$ exists.
	Moreover, we have $\ld{\gam_n}(F,G) \leq \LD$ for all $n \in \N$. 
	Thus, it suffices to show that $\LD \leq \limN \ld{\gam_n}(F,G)$.
	By \cref{eq:appLd_SymmetricAltDefOfLd} we have for all $n \in \N$ and $t \in \R$
	\begin{align*}
		G(t) \leq F(t + \gam_n) + \ld{\gam_n}(F,G) \text{ and } 
		F(t) \leq G(t + \gam_n) + \ld{\gam_n}(F,G).             
	\end{align*}
	By the continuity from the right of $F$ and $G$, this yields
	\begin{align*}
		G(t) & \leq \limN F(t + \gam_n) + \ld{\gam_n}(F,G) = F(t + \gam) + \limN \ld{\gam_n}(F,G) \text{ and } \\
		F(t) & \leq \limN G(t + \gam_n) + \ld{\gam_n}(F,G) = G(t + \gam) + \limN \ld{\gam_n}(F,G),             
	\end{align*}
	which implies $\LD \leq \limN \ld{\gam_n}(F,G)$. Thus, the function $\gam \mapsto \LD$ is continuous from the right.
	
	In order to prove \Cref{it:ld_bpAltDef} we write $\mathcal{S}$ as a shorthand for 
	$\{\epsilon \geq 0: F(t-\gam) - \epsilon \leq G(t) \leq F(t+\gam) + \epsilon \fatir \}$.
	It is easy to see that
	\begin{align*}
		\mathcal{S} 
		  & = \{\epsilon \geq 0: F(t) \leq G(t + \gam) + \epsilon, G(t) \leq F(t+\gam) + \epsilon \fatir \}  \\
		  & = \{\epsilon \geq 0: F(t) - G(t + \gam) \leq \epsilon, G(t) - F(t+\gam) \leq \epsilon \fatir \}, 
	\end{align*}
	which can be rewritten as
	\begin{align*}
		  & \{\epsilon \geq 0: \max(F(t) - G(t + \gam), G(t) - F(t+\gam)) \leq \epsilon \fatir \}                      \\
		  & = \{\epsilon \geq 0: \sup_{t \in \R} \max(F(t) - G(t + \gam), G(t) - F(t+\gam)) \leq \epsilon \}           \\
		  & = \left[ \max \left(0, \sup_{t \in \R} \max(F(t) - G(t + \gam), G(t) - F(t+\gam)) \right), \infty \right). 
	\end{align*}
	Recalling the definition of $\LD$ we conclude
	\begin{align*}
		\LD = \inf( \mathcal{S})                                                                     
		= \max \left(0, \sup_{t \in \R} \max(F(t) - G(t + \gam), G(t) - F(t+\gam)) \right). 
	\end{align*}
	However, as $F$ and $G$ are distribution functions we have $\lim_{x \to \infty} F(x) - G(x+\gam) = 0$. Thus, 
	$\sup_{t \in \R} \max(F(t) - G(t + \gam), G(t) - F(t+\gam))$ is nonnegative and we can remove the outer maximum.

	\Cref{it:ld_bpTriangle} is a direct consequence of \Cref{it:ld_bpInfimumAttained} and \Cref{it:ld_bpAltDef}
	as for every $t \in \R$ we have
	\begin{align*}
		  & H(t - \gam_1 - \gam_2) - \ld{\gam_2}(G,H) - \ld{\gam_1}(F,G) 
		\leq G(t - \gam_1) - \ld{\gam_1}(F,G)
		\leq F(t)\\
		  & \leq G(t + \gam_1) + \ld{\gam_1}(F,G)                        
		\leq H(t + \gam_1 + \gam_2) + \ld{\gam_2}(G,H) + \ld{\gam_1}(F,G).
	\end{align*}

	We next prove \ref{it:ld_bpLevy}. By the monotonicity of $F$ and \cref{eq:ld_ldAttainsMinimum} we have for every $t \in \R$:
	\begin{align*}
		  & F(t - \max(\gam, \LD)) - \max(\gam, \LD)       
		\leq F(t - \gam) - \LD \\
		  & \leq G(t)                                      
		\leq F(t + \gam) + \LD \\
		  & \leq F(t + \max(\gam, \LD)) + \max(\gam, \LD). 
	\end{align*}
	Thus, $L(F,G) \leq \max(\gam, \LD)$.
	Now suppose that $\gam \geq L(F,G)$. 
	For every $t \in \R$ the monotonicity of $F$ yields
	\begin{align*}
		  & F(t - \gam) - L(F,G)        
		\leq F(t - L(F,G)) - L(F,G) 
		\leq G(t) \\
		  & \leq F(t + L(F,G)) + L(F,G) 
		\leq F(t + \gam) + L(F,G),
	\end{align*}
	which implies $\LD \leq L(F,G)$ by the definition of $\LD$.
	
	It remains to show $\min(\gam, \LD) \leq L(F,G)$.
	For this, we distinguish two cases: If $\gam \leq L(F,G)$, the inequality is trivially fulfilled. If, otherwise, $\gam > L(F,G)$ holds true, then $\LD \leq L(F,G)$ as shown before.   
	
	\Cref{it:ld_bpScaling} is a direct consequence of \Cref{it:ld_bpAltDef}.
\end{proof}

\begin{proof}[Proof of \Cref{cor:ld_ldWeakConvergence}:]
	If $F_n$ converges weakly to $F$, then $\limN L(F_n, F) = 0$. Applying \Cref{it:ld_bpLevy} this implies
	\begin{align*}
		\limsup_{n \to \infty} \min(\gam, \ld{\gam}(F_n, F)) 
		\leq \limsup_{n \to \infty} L(F_n, F) = 0            
	\end{align*}
	for every $\gam \geq 0$ which shows that $\limN \ld{\gam}(F_n, F) = 0$ for all positive $\gam$.
	
	For the other direction, we start with an arbitrary $\gam > 0$. 
	For every such $\gam$ we have $\limN \ld{\gam}(F_n, F) = 0$ by assumption.
	Now, \Cref{it:ld_bpLevy} implies
	\begin{align*}
		0 \leq \limsup_{n \to \infty} L(F_n, F)                           
		\leq \limsup_{n \to \infty} \max(\gam, \ld{\gam}(F_n, F)) = \gam. 
	\end{align*}
	As $\gam > 0$ can be made arbitrarily small, we conclude $\limN L(F_n, F) = 0$, which proves the weak convergence.
\end{proof}

\begin{proof}[Proof of \Cref{prop:ld_quantileInequality}:]
	It suffices to show $\Q{\alpha}{G} \leq \Q{\alpha + \LD}{F} + \gam$ for every $\alpha \in \R$ because we can exchange the roles of $F$ and $G$.
	For this, we distinguish three cases:
	If $\Q{\alpha + \LD}{F} = \infty$ the inequality is trivially fulfilled.
	If $\Q{\alpha + \LD}{F} = - \infty$ we conclude $\alpha + \LD \leq 0$. 
	In particular, this shows that $\alpha \leq 0$ and therefore we also have $\Q{\alpha}{G} = -\infty$,
	which proves the inequality.
	In the last case we have $\Q{\alpha + \LD}{F} \in \R$ and \Cref{it:ld_bpInfimumAttained} yields
	\begin{align*}
		G(\Q{\alpha + \LD}{F} + \gam)
		  & \geq F(\Q{\alpha + \LD}{F}) - \LD 
		\geq \alpha.
	\end{align*}
	Thus, $\Q{\alpha}{G} = \inf\{x \in \R: G(x) \geq \alpha\} \leq \Q{\alpha + \LD}{F} + \gam$, which proves the claim.
\end{proof}

\begin{proof}[Proof of \Cref{lem:ld_ldDefinitionViaQuantiles}:]
	Assume $\LD \leq \epsilon$. 
	Then the monotonicity of the quantiles together with \Cref{prop:ld_quantileInequality} implies
	\begin{align*}
		\Q{\alpha-\epsilon}{G} - \gam 
		  & \leq \Q{\alpha-\LD}{G} - \gam 
		\leq  \Q{\alpha}{F} 
		\leq \Q{\alpha+\LD}{G} + \gam 
		\leq \Q{\alpha+\epsilon}{G} + \gam.
	\end{align*}
	For the other direction, we assume that $\LD > \epsilon \geq 0$.
	In view of \cref{eq:appLd_SymmetricAltDefOfLd} this implies the existence of $t^\ast \in \R$, such that
	either $F(t^\ast) > G(t^\ast + \gam) + \epsilon$ or $G(t^\ast) > F(t^\ast + \gam) + \epsilon$ holds true.
	W.l.o.g. we assume $F(t^\ast) > G(t^\ast + \gam) + \epsilon$.
	Defining $\alpha = F(t^\ast)$ we conclude $\alpha \in (\epsilon, 1]$. 
	Thus, we have $\Q{\alpha}{F} \leq t^\ast < \infty$ by definition. Furthermore, we conclude $\Q{\alpha}{F} > -\infty$ because $\alpha > 0$.
	In other words, we have $\Q{\alpha}{F} \in \R$, which yields
	\begin{align*}
		G(\Q{\alpha}{F} + \gam) \leq G(t^\ast + \gam) < F(t^\ast) - \epsilon 
		= \alpha - \epsilon \leq G(\Q{\alpha-\epsilon}{G}).                  
	\end{align*}
	Hence, we conclude $\Q{\alpha}{F} + \gam < \Q{\alpha-\epsilon}{G}$.   
\end{proof}

\begin{proof}[Proof of \Cref{cor:ld_generalPiInequality}]
	\Cref{cor:ld_generalPiInequality} is a direct consequence of \Cref{prop:ld_quantileInequality}.
\end{proof}

\begin{proof}[Proof of \Cref{lem:ld_WassersteinTypeBound}:]
	Let $X \sim \mu_F$, $Y \sim \mu_G$ and $\P \in \Gamma(\mu_F, \mu_G)$.
	We then have for every $t \in \R$:
	\begin{align*}
		F(t) - G(t + \gam) 
		  & = \mu_F(\{X \leq t\}) - \mu_G(\{Y \leq t + \gam\}) 
		\leq \P(X \leq t < Y - \gam) \\
		  & \leq \P(Y - X > \gam) \leq \P(|Y-X| > \gam).       
	\end{align*}
	By exchanging the roles of $F$ and $G$, we also get the inequality
	$G(t) - F(t + \gam) \leq \P(|Y-X| > \gam)$. 
	Recalling \cref{it:ld_bpAltDef} and the fact that $t \in \R$ was arbitrary we conclude
	\begin{align*}
		\LD = \sup_{t \in \R} \max(F(t) - G(t + \gam), G(t) - F(t + \gam)) \leq \P(|Y-X| > \gam). 
	\end{align*}
	Since $X \sim \mu_F$, $Y \sim \mu_G$ and $\P \in \Gamma(\mu_F, \mu_G)$ was arbitrary, we conclude
	\begin{align*}
		\LD \leq \inf_{\substack{X \sim \mu_F, Y \sim \mu_G, \\ \P \in \Gamma(\mu_F,\mu_G)}} \P( |X-Y| > \gam).
	\end{align*}
	To finish the proof, we apply Markov's inequality.
\end{proof}

\begin{proof}[Proof of \Cref{cor:ld_BoundForEcdfs}:]
	This can be seen from \Cref{lem:ld_WassersteinTypeBound}.
\end{proof}
    }{}
\end{document}